\ifx\shlhetal\undefinedcontrolsequenc\let\shlhetal\relax\fi

\documentstyle[11pt]{article}

\input amssym.def
\input amssym.tex

\newcommand{\ZFC}{{\rm ZFC}}
\newcommand{\GCH}{{\rm GCH}}
\newcommand{\ccc}{{\rm c.c.c.}}
\newcommand{\cc}{{\rm c.c.}}
\newcommand{\acc}{{\rm acc}}
\newcommand{\nacc}{{\rm nacc}}
\newcommand{\com}{{\rm com}}
\newcommand{\Card}{{\rm Card}}
\newcommand{\cel}{{\rm c}}
\newcommand{\cf}{{\rm cf}\/} 
\newcommand{\pcf}{{\rm pcf}\/}
\newcommand{\Reg}{{\rm Reg}}
\newcommand{\lh}{{\rm lg}\/}
\newcommand{\Levy}{{\rm Levy}} 
\newcommand{\id}{{\rm id}}
\newcommand{\bd}{{\rm bd}}
\newcommand{\tcf}{{\rm tcf}}
\newcommand{\rest}{{\restriction}} 
\newcommand{\suc}{{\rm succ}} 
\newcommand{\dom}{{\rm dom}} 
\newcommand{\rng}{{\rm rng}}
\newcommand{\red}{{\rm red}}
\newcommand{\green}{{\rm green}}
\newcommand{\forces}{\Vdash} 
\newcommand{\otp}{{\rm otp}}
\newcommand{\wmk}{{\rm wmk}}
\newcommand{\swmk}{{\rm smk}}

\newcommand{\QED}{\hfill\vrule width 6pt height 6pt depth 0pt 
\vspace{0.1in}} 
\newcommand{\Proof}{\noindent {\sc Proof} \hspace{0.2in}} 
\newcommand{\gt}{{\frak t}}
\newcommand{\bz}{{\bf 0}}
\newcommand{\bo}{{\bf 1}}
\newcommand{\A}{{\cal A}}  
\newcommand{\ba}{{\Bbb B}}
\newcommand{\C}{{\cal C}}
\newcommand{\F}{{\cal F}} 
\newcommand{\cH}{{\cal H}}
\newcommand{\K}{{\cal K}}
\newcommand{\q}{{\Bbb Q}}
\renewcommand{\P}{{\cal P}}
\newcommand{\p}{{\Bbb P}}
\renewcommand{\S}{{\cal S}}
\newcommand{\bV}{{\bf V}} 
\newcommand{\U}{{\cal U}}

\newtheorem{theorem}{Theorem}[section] 
\newtheorem{claim}{Claim}[theorem]

\newtheorem{lemma}[theorem]{Lemma} 
\newtheorem{proposition}[theorem]{Proposition} 
 
\newtheorem{problem}[theorem]{Problem} 
\newtheorem{thema}[theorem]{Thema} 
\newtheorem{pws}[theorem]{The Problem We Address} 
\newtheorem{pms}[theorem]{Question (mostly solved)}
\newtheorem{definition}[theorem]{Definition}
\newtheorem{notation}[theorem]{Notation}
\newtheorem{remark}[theorem]{Remark}
\newtheorem{conclusion}[theorem]{Conclusion}
\newtheorem{mainconc}[theorem]{Main Conclusion}
\newtheorem{comment}[theorem]{Comment}
\newtheorem{conrem}[theorem]{Concluding Remarks}

\title{Cellularity of free products of Boolean algebras (or topologies)} 
\author{
{\bf Saharon Shelah}\thanks{The research partially supported by ``The Israel
Science Foundation'' administered by The Israel Academy of Sciences and
Humanities. Publication 575. We thank Andrzej Roslanowski for writing sections
1 - 5 from lectures, 6 - 7 from notes}\\
Institute of Mathematics\\
The Hebrew University of Jerusalem\\
91904 Jerusalem, Israel\\
and\\
Department of Mathematics\\
Rutgers University\\
New Brunswick, NJ 08854, USA
}

\date\today 

\begin{document} 
\baselineskip13.14 truept

\nocite{EK}
\nocite{HJSh:249}
\nocite{HaJuSz}
\nocite{Jn}
\nocite{Ju}
\nocite{KnMg78}
\nocite{LMSh:198}
\nocite{MgSh:204}
\nocite{M1}
\nocite{M2}
\nocite{Sh:g}
\nocite{Sh:92}
\nocite{Sh:93}
\nocite{Sh:126}
\nocite{Sh:345}
\nocite{Sh:345a}
\nocite{Sh:355}
\nocite{Sh:410}
\nocite{Sh:420}
\nocite{Sh:430}
\nocite{Sh:460}
\nocite{Sh:481}
\nocite{Sh:513}
\nocite{Sh:572}
\nocite{Sh:576}

\nocite{So74}
\setcounter{section}{-1}
\setcounter{page}{0}
\maketitle

\begin{abstract}
The aim of this paper is to present an answer to Problem 1 of \cite{M1},
\cite{M2}. We do this by proving in particular that 
\begin{quotation}
\noindent if {\em $\mu$ is a strong limit singular cardinal,
$\theta=(2^{\cf(\mu)})^+$, $2^\mu=\mu^+$}

\noindent then {\em there are Boolean algebras $\ba_1,\ba_2$ such that
$\cel(\ba_1)=\mu$, 

$\cel(\ba_2)<\theta$ but $\cel(\ba_1*\ba_2)=\mu^+$. 
}
\end{quotation}
Further we improve this result, deal with the method and the
neccessity of the assumptions. 
\end{abstract}
\eject

\section{Introduction} 
\begin{notation}
\label{notation}
\begin{enumerate}
\item In the present paper all cardinals are infinite so we will not repeat
this additional demand. Cardinals will be denoted by $\lambda$, $\mu$,
$\theta$ (with possible indexes) while ordinal numbers will be called
$\alpha$, $\beta$, $\zeta$, $\xi$, $\varepsilon$, $i$, $j$. Usually $\delta$
will stand for a limit ordinal (we may forget to repeat this assumption). 
\item Sequences of ordinals will be called $\eta$, $\nu$, $\rho$ (with
possible indexes). For sequences $\eta_1,\eta_2$ their longest common initial
segment is denoted by $\eta_1\wedge\eta_2$. The length of the sequence $\eta$
is $\lh(\eta)$.
\item Ideals are supposed to be proper and contain all singletons. For a limit
ordinal $\delta$ the ideal of bounded subsets of $\delta$ is denoted by
$J^{\bd}_\delta$. If $I$ is an ideal on a set $X$ then $I^+$ is the family of
$I$-large sets, i.e. 
\[a\in I^+\quad\mbox{ if and only if }\quad a\subseteq X\ \&\ a\notin I\]
and $I^c$ is the dual filter of sets with the complements in $I$.
\end{enumerate}
\end{notation}

\begin{notation}
\label{morenotation}
\begin{enumerate}
\item In a Boolean algebra we denote the Boolean operations by $\cap$ (and
$\bigcap$\/), $\cup$ (and $\bigcup$), $-$. The distinguished elements are
$\bz$ and $\bo$. In the cases which may be confusing we will add indexes to
underline in which Boolean algebra the operation (the element) is considered,
but generally we will not do it. 
\item For a Boolean algebra $\ba$ and an element $x\in\ba$ we denote:
\[x^0=x\quad\mbox{ and }\quad x^1=-x.\]
\item The free product of Boolean algebras $\ba_1$, $\ba_2$ is denoted by
$\ba_1*\ba_2$. We will use $\bigstar$ to denote the free product of a family
of Boolean algebras. For $\bigstar_{i<\sigma} \ba_i$ (wlog $\ba_i$
pairwise disjoint) is the Boolean algebra generated by the formal
intersection $\bar b= \bigcap\limits_{i<\sigma} b_i$, $b_i\in \ba$
freely except $\tau(\bar b^1, \ldots, \bar b^n)=0_{\ba}$ iff some
$u\in \sigma$ finite non empty for simplicity
$$
\bigstar_{i\in u} \ba_i \models \tau(\bigcap_{i\in u} b^1_i\ldots)\neq
0
$$
\end{enumerate}
\end{notation}

\begin{definition}
\begin{enumerate}
\item A Boolean algebra $\ba$ satisfies the $\lambda$-cc if there is no
family $\F\subseteq\ba^+\stackrel{\rm def}{=}\ba\setminus\{\bz\}$ such
that $|\F|=\lambda$ and any two members of $\F$ are disjoint (i.e. their meet
in $\ba$ is $\bz$).
\item The cellularity of the algebra $\ba$ is
\[\cel(\ba)=\sup\{|\F|: \F\subseteq\ba^+\ \&\ (\forall x,y\in\F)(x\neq
y\ \Rightarrow x\cap y=\bz)\},\]
\[\cel^+(\ba)=\sup\{|\F|^+: \F\subseteq\ba^+\ \&\ (\forall x,y\in\F)(x\neq
y\ \Rightarrow x\cap y=\bz)\}.\]
\item For a topological space $(X,\tau)$:
\[\begin{array}{ll}
\cel(X,\tau)=\sup\{|\U|:&\U\mbox{ is a family of pairwise disjoint}\\
\ &\mbox{nonempty open sets}\}.\\
\end{array}\]
\end{enumerate}
\end{definition}
The problem can be posed in each of the three ways ({\em $\lambda$-cc} is the
way of forcing, {\em the cellularity of Boolean algebras} is the approach of
Boolean algebraists, and {\em the cellularity of a topological space} is the
way of general topologists). It is well known that the three are equivalent,
though (1) makes the attainment problem more explicite. We use the second
approach. 

A stronger property then $\lambda$-cc is the $\lambda$-Knaster property. This
property behaves nicely in free products -- it is productive. We will use it
in our construction.

\begin{definition}
\label{knaster}
A Boolean algebra $\ba$ has the $\lambda$-Knaster property if for every
sequence $\langle z_\varepsilon: \varepsilon<\lambda\rangle\subseteq\ba^+$
there is $A\in [\lambda]^{\textstyle \lambda}$ such that
\[\varepsilon_1,\varepsilon_2\in A\ \ \ \Rightarrow\ \ \ z_{\varepsilon_1}\cap
z_{\varepsilon_2}\neq\bz.\] \end{definition}

We are interested in the behaviour of the cellularity of Boolean algebras when
the free product of them is considered.

\begin{thema}
\label{thema}
When, for Boolean algebras $\ba_1$, $\ba_2$
\[\cel^+(\ba_1)\leq\lambda_1\ \&\ \cel^+(\ba_2)\leq\lambda_2\ \ \
\Rightarrow\ \ \ \cel^+(\ba_1*\ba_2)\leq\lambda_1+\lambda_2 ?\]
\end{thema}

There are a lot of results about it, particularly if $\lambda_1=\lambda_2$
(see \cite{Sh:g} or \cite{M1}, more \cite{Sh:572}). It is well know
that if 
\[(\lambda_1^+ +\lambda_2^+)\longrightarrow (\lambda_1^+,\lambda_2^+)^2\]
then the answer is ``yes''. These are exactly the cases for which ``yes''
answer is known. Under {\bf GCH} the only problem which remained open was the
one presented below: 

\begin{pws}
\label{pwa}
\ \\
{\bf (posed by D. Monk as Problem 1 in \cite{M1}, \cite{M2} under $\GCH$)}\\
Are there Boolean algebras $\ba_1$, $\ba_2$ and cardinals $\mu,\theta$ such
that 
\begin{enumerate}
\item $\lambda_1=\mu$ is singular, $\mu>\lambda_2=\theta>\cf(\mu)$ and
\item $\cel(\ba_1)=\mu$, $\cel(\ba_2)\leq\theta$ but $\cel(\ba_1*\ba_2)>\mu$?
\end{enumerate}
\end{pws}

We will answer this question proving in particular the following result (see
\ref{main}): 
\begin{quotation}
If {\em $\mu$ is a strong limit singular cardinal,
$\theta=(2^{\cf(\mu)})^+$, $2^\mu=\mu^+$}

then {\em there are Boolean algebras $\ba_1,\ba_2$ such that
$\cel(\ba_1)=\mu$, 

$\cel(\ba_2)<\theta$ but $\cel(\ba_1*\ba_2)=\mu^+$. 
}
\end{quotation}
Later we deal with better results by refining the method.

\begin{remark}
{\em 
On products of many Boolean algebras and square bracket arrows see
\cite{Sh:345}, 1.2A, 1.3B.

If $\lambda\longrightarrow [\mu]^2$, $[\tau<\sigma\ \Rightarrow\
2^\tau<\theta]$, the cardinals $\theta,\sigma$ are possibly finite,
$\ba_i$ (for $i<\theta$) are Boolean algebras such that for each $j<\theta$
the free product $\mathop{\bigstar}\limits_{i\in\theta\setminus\{j\}}\ba_i$
satisfies the $\mu$-cc {\em then} the algebra $\ba=\mathop{\bigstar}
\limits_{i<\theta}\ba_i$ satisfies the $\lambda$-cc.  
}
\end{remark}
[Why? if $\langle a^\zeta_i: i<\theta\rangle\in
\prod\limits_{i<\theta}{\ba}^+_i$ for $\zeta<\lambda$ such that for
every $\zeta<\xi<\lambda$ for some $i$, $\ba_i \models$
``$a^\zeta_i\cap a^\xi_i=\bz$'', let $\i=\i(\zeta, \xi)$, and we can
find $A\in [\lambda]^\mu$, $j<\theta$ such that for $\zeta<\xi$ from
$A$, $\i(\zeta, \xi)=j$, so $\langle a^\zeta_i: i<\theta, i\neq
i^*\rangle$ for $\zeta\in A$ exemplifies $\bigstar_{i\in
\theta\setminus\{i^*\}}{\ba}_i$ fails the $\mu-\cc$. We can deal also
with ultraproducts and other products similarly.]

\section{Preliminaries: products of ideals}
\begin{notation}
\label{quantifiers}
For an ideal $J$ on $\delta$ the quantifier $(\forall^J i<\delta)$ means ``for
all $i<\delta$ except a set from the ideal'', i.e.
\[(\forall^J i<\delta)\varphi(i)\ \equiv\ \{i<\delta: \neg\varphi(i)\}\in J.\]
The dual quantifier $(\exists^J i<\delta)$ means ``for a $J$-positive set of
$i<\delta$''. 
\end{notation}

\begin{proposition}
\label{prefub}
Assume that $\lambda^{0}>\lambda^{1}>\ldots>\lambda^{n-1}$ are cardinals,
$I^\ell$ are ideals on $\lambda^\ell$ (for $\ell<n$) and
$B\subseteq\prod\limits_{\ell<n}\lambda^\ell$. Further suppose that  
\begin{description}
\item[$(\alpha)$] $(\exists^{I^{0}}\gamma_{0})\ldots (\exists^{I^{n-1}}
\gamma_{n-1})(\langle\gamma_\ell:\ell<n\rangle\in B)$
\item[$(\beta)$] the ideal $I^\ell$ is $(2^{\lambda^{\ell+1}})^+$-complete
(for $\ell+1<n$).
\end{description}
{\em Then} there are sets $X_\ell\subseteq\lambda^\ell$, $X_\ell\notin I^\ell$
such that $\prod\limits_{\ell<n} X_\ell\subseteq B$.

\noindent [Note that this translates the situation to arity 1; it is a kind of
polarized $(1,\ldots,1)$-partitions with ideals.]
\end{proposition}

\Proof We show it by induction on $n$. Define
\[\begin{array}{ll}
E_0\stackrel{\rm def}{=}&\{(\gamma',\gamma''): \gamma',\gamma''<\lambda^0
\quad\mbox{ and}\\ 
\ &\mbox{for all }\gamma_1<\lambda^1,\ldots,\gamma_{n-1}<\lambda^{n-1}\ \mbox{
we have}\\ 
\ &(\langle\gamma',\gamma_1,\ldots,\gamma_{n-1}\rangle\in B\ \ \
\Leftrightarrow\ \ \ \langle\gamma'',\gamma_1,\ldots,\gamma_{n-1}\rangle
\in B)\}. 
  \end{array}\]
Clearly $E_0$ is an equivalence relation on $\lambda^0$ with $\leq
2^{\prod_{0<m<n}\lambda^m}=2^{\lambda^1}$ equivalence classes. Hence the set
\[A_0\stackrel{\rm def}{=}\bigcup\{A: A\mbox{ is an }E_0\mbox{-equivalence
class}, A\in I^0\}\]
is in the ideal $I^0$. Let 
\[A_0^*\stackrel{\rm def}{=}\{\gamma_0<\lambda^0:(\exists^{I^{1}}\gamma_{1})
\ldots (\exists^{I^{n-1}}\gamma_{n-1})(\langle\gamma_0,\gamma_1,\ldots,
\gamma_{n-1}\rangle\in B).\]
The assumption $(\alpha)$ implies that $A_0^*\notin I^0$ and hence we may
choose $\gamma^*_0\in A_0^*\setminus A_0$. Let 
\[B_1\stackrel{\rm def}{=}\{\bar{\gamma}\in\prod_{k=1}^{n-1}\lambda^k:
\langle\gamma^*_0 \rangle\hat{\ }\bar{\gamma}\in B\}.\]
Since $\gamma^*_0\in A_0^*$ we are sure that 
\[(\exists^{I^{1}}\gamma_{1})\ldots (\exists^{I^{n-1}}\gamma_{n-1})(\langle
\gamma_1,\ldots, \gamma_{n-1}\rangle\in B_1).\]
Hence we may apply the inductive hypothesis for $n-1$ and $B_1$ and we find
sets $X_1\in (I^1)^+,\ldots,X_{n-1}\in (I^{n-1})^+$ such that
$\prod\limits_{\ell=1}^{n-1}X_\ell\subseteq B_1$, so then
\[(\forall\gamma_1\in X_1)\ldots(\forall\gamma_{n-1}\in X_{n-1})(\langle
\gamma^*_0,\gamma_1,\ldots,\gamma_{n-1}\rangle\in B).\] 
Take $X_0$ to be the $E_0$-equivalence class of $\gamma^*_0$ (so $X_0\in
(I^0)^+$ as $\gamma_0^*\notin A_0$). By the definition of the relation $E_0$
and the choice of the sets $X_\ell$ we have that for each $\gamma_0\in X_0$
\[(\forall\gamma_1\in X_1)\ldots(\forall\gamma_{n-1}\in X_{n-1})(\langle
\gamma_0,\gamma_1,\ldots,\gamma_{n-1}\rangle\in B)\] 
what means that $\prod\limits_{\ell<n}X_\ell\subseteq B$. The proposition is
proved. \QED

\begin{proposition}
\label{weakfub}
Assume that $\lambda_0>\lambda_1>\ldots>\lambda_{n-1}\geq\sigma$ are
cardinals, $I_\ell$ are ideals on $\lambda_\ell$ (for $\ell<n$) and
$B\subseteq\prod\limits_{\ell<n}\lambda_\ell$. Further suppose that
\begin{description}
\item[$(\alpha)$] $(\exists^{I_0}\gamma_0)\ldots(\exists^{I_{n-1}}
\gamma_{n-1})(\langle\gamma_\ell:\ell<n\rangle\in B)$
\item[$(\beta)$]  for each $\ell<n-1$ the ideal $I_\ell$ is
$((\lambda_{\ell+1})^\sigma)^+$-complete, $[\lambda_{n-1}]^{<\sigma}\subseteq
I_{n-1}$. 
\end{description}
Then there are sets $X_\ell\in [\lambda_\ell]^\sigma$ such that
$\prod\limits_{\ell<n} X_\ell\subseteq B$. 
\end{proposition}

\Proof The proof is by induction on $n$. If $n=1$ then there is nothing to do
as $I_{n-1}$ contains all subsets of $\lambda_{n-1}$ of size
$<\sigma$ and $\lambda_{n_1}\geq \sigma$ so every $A\in I^+_{n_1}$ has
cardinality $\geq \sigma$. 

\noindent Let $n>1$ and let
\[a_0\stackrel{\rm def}{=}\{\gamma\in\lambda_0: (\exists^{I_1}\gamma_1)\ldots
(\exists^{I_{n-1}}\gamma_{n-1})(\langle\gamma,\gamma_1,\ldots,\gamma_{n-1}
\rangle\in B)\}.\]
By our assumptions we know that $a_0\in (I_0)^+$. For each $\gamma\in
a_0$ we may apply the inductive hypothesis to the set
\[B_\gamma\stackrel{\rm def}{=}\{\langle\gamma_1,\ldots,\gamma_{n-1}\rangle\in
\prod_{0<\ell<n}\lambda_\ell:\langle\gamma,\gamma_1,\ldots,\gamma_{n-1}
\rangle\in B\}\]
and we get sets $X^\gamma_1\in [\lambda_1]^\sigma,\ldots,X^\gamma_{n-1}\in
[\lambda_{n-1}]^\sigma$ such that 
\[\prod\limits_{0<\ell<n}X^\gamma_\ell\subseteq B_\gamma.\]
There is at most $(\lambda_1)^\sigma$ possible sequences $\langle X^\gamma_1,
\ldots,X^\gamma_{n-1}\rangle$, the ideal $I_0$ is
$((\lambda_1)^\sigma)^+$--complete  so for some sequence $\langle X_1,\ldots,
X_{n-1}\rangle$ and a set $a^*\subseteq a_0$, $a^*\in (I_0)^+$ we have
\[(\forall\gamma\in a^*)(X^\gamma_1=X_1\ \&\ \ldots\ \&\ X_{n-1}^\gamma=
X_{n-1}).\]
Choose $X_0\in [a^*]^\sigma$ (remember that $I_0$ contains singletons and it
is complete enough to make sure that $\sigma\leq |a^*|$). Clearly
$\prod\limits_{\ell<n} X_\ell\subseteq B$. \QED

\begin{remark}
{\em 
We can use $\sigma_0\geq \sigma_1\geq \ldots\geq\sigma_{n-1}$, $I_\ell$
is $(\lambda^{\sigma_{\ell+1}}_{\ell+1})^+$-complete,
$[\lambda_\ell]^{<\sigma_\ell}\subseteq I_\ell$.
}
\end{remark}

\begin{proposition}
\label{fubini}
Assume that $n<\omega$ and $\lambda^m_\ell$, $\chi^m_\ell$, $P^m_\ell$,
$I^m_\ell$, $I^m$ and $B$ are such that for $\ell,m\leq n$:
\begin{description}
\item[$(\alpha)$] $I^m_\ell$ is a $\chi^m_\ell$-complete ideal on
$\lambda^m_\ell$ (for $\ell,m\leq n$)
\item[$(\beta)$]  $P^m_\ell\subseteq\P(\lambda^m_\ell)$ is a family dense in
$(I^m_\ell)^+$ in the sense that:
\[(\forall X\in (I^m_\ell)^+)(\exists a\in P^m_\ell)(a\subseteq X)\]
\item[$(\gamma)$] $I^m=\{X\!\subseteq\!\prod\limits_{\ell\leq n}
\lambda^m_\ell: \neg(\exists^{I^m_0}\gamma_0)\ldots(\exists^{I^m_n}\gamma_n)
(\langle\gamma_0,\ldots,\gamma_n\rangle\in X)\}$

\noindent [thus $I^m$ is the ideal on $\prod\limits_{\ell\leq
n}\lambda^m_\ell$ such that the dual filter $(I^m)^c$ is the Fubini product of
filters $(I^m_0)^c,\ldots,(I^m_n)^c$]
\item[$(\delta)$] $\chi^m_{n-m}>\sum\limits_{\ell=m+1}^{n}(|P^\ell_{n-\ell}|+
\sum\limits_{k=0}^{n-\ell}\lambda^\ell_k)$ 
\item[$(\varepsilon)$] $B\subseteq\prod\limits_{m\leq n}\prod\limits_{\ell\leq
n}\lambda^m_\ell$ is a set satisfying
\[(\exists^{I^0}\eta_0)(\exists^{I^1}\eta_1)\ldots (\exists^{I^n}\eta_n)
(\langle\eta_0,\eta_1,\ldots,\eta_n\rangle\in B).\]
\end{description}
{\em Then} there are sets $X_0,\ldots,X_n$ such that for $m\leq n$:
\begin{description}
\item[(a)] $X_m\subseteq\prod\limits_{\ell\leq n-m}\lambda^m_\ell$
\item[(b)] if $\eta,\nu\in X_m$, $\eta\neq\nu$ then
\begin{description}
\item[(i)\ ] $\eta\rest (n-m)=\nu\rest (n-m)$
\item[(ii) ] $\eta(n-m)\neq\nu (n-m)$
\end{description}
\item[(c)] $\{\eta(n-m):\eta\in X_m\}\in P^m_{n-m}$\ \ \ \ \ \ \ \ \ \ and
\item[(d)] for each $\langle\eta_0,\ldots,\eta_n\rangle\in\prod\limits_{m\leq
n} X_m$ there is $\langle\eta^*_0,\ldots,\eta^*_n\rangle\in B$ such that  

$(\forall m\leq n)(\eta_m\trianglelefteq\eta^*_m)$.
\end{description}
\end{proposition}

\noindent{\bf Remark}~\ref{fubini}.A:
\begin{enumerate}
\item Note that the sets $X_m$ in the assertion of \ref{fubini} may be thought
of as sets of the form $X_m=\{\nu_m\hat{\ }\langle\alpha\rangle: \alpha\in
a_m\}$ for some $\nu_m\in\prod\limits_{\ell<n-m}\lambda^m_\ell$ and $a_m\in
P^m_{n-m}$. 
\item We will apply this proposition with
$\lambda^m_\ell=\lambda_\ell$, $I^m_\ell=I_\ell$ and\\
$\lambda_\ell>\chi_\ell>\sum_{k<\ell}\lambda_k$.
\item In the assumption $(\delta)$ of \ref{fubini} we may have that the last
sum on the right hand side of the inequality ranges from $k=0$ to $n-\ell-1$.
We did not formulate that assumption in this way as with $n-\ell$ there it is
easier to handle the induction step and this change is not important for our
applications. 
\item In the assertion {\bf (d)} of \ref{fubini} we can make $\eta^*_\ell$
depending on $\langle\eta_0,\ldots,\eta_\ell\rangle$ only.
\end{enumerate}
\medskip

\Proof The proof is by induction on $n$. For $n=0$ there is nothing to do. Let
us describe the induction step.

Suppose $0<n<\omega$ and $\lambda^m_\ell$, $\chi^m_\ell$, $P^m_\ell$,
$I^m_\ell$, $I^m$ (for $\ell,m\leq n$) and $B$ satisfy the assumptions
$(\alpha)$--$(\varepsilon)$. Let
\[B^*\stackrel{\rm def}{=}\{\langle\eta_0,\eta_1\rest n,\ldots,\eta_n\rest
n\rangle: \eta_m\in\prod_{\ell\leq n}\lambda^m_\ell\mbox{ (for $m\leq n$) and
} \langle\eta_0,\eta_1,\ldots,\eta_n\rangle\in B\}\]
and for $\eta_0\in\prod\limits_{\ell\leq n}\lambda^0_\ell$ let
\[B^*_{\eta_0}\stackrel{\rm def}{=}\{\langle\nu_1,\ldots,\nu_n\rangle\in
\prod_{m=1}^n\prod_{\ell=0}^{n-1}\lambda^m_\ell: \langle\eta_0,\nu_1,\ldots,
\nu_n\rangle\in B^*\}.\]
Let $J^m$ (for $1\leq m\leq n$) be the ideal on $\prod\limits_{\ell=0}^{n-1}
\lambda^m_\ell$ coming from the ideals $I^m_\ell$, i.e. a set $X\subseteq
\prod\limits_{\ell<n}\lambda^m_\ell$ is in $J^m$ if and only if
\[\neg (\exists^{I^m_0}\gamma_0)\ldots(\exists^{I^m_{n-1}}\gamma_{n-1})
(\langle\gamma_0,\ldots,\gamma_{n-1}\rangle\in X).\]
Let us call the set $B^*_{\eta_0}$ {\em big} if
\[(\exists^{J^1}\nu_1)\ldots(\exists^{J^n}\nu_n)(\langle\nu_1,\ldots,\nu_n
\rangle\in B^*_{\eta_0}).\]
We may write more explicitely what the bigness means: the above condition is
equivalent to
\[\begin{array}{ll}
\ &(\exists^{I^1_0}\gamma^1_0)\ldots (\exists^{I^1_{n-1}}\gamma^1_{n-1})
\ldots\\  
\ldots &(\exists^{I^n_0}\gamma^n_0)\ldots (\exists^{I^n_{n-1}}\gamma^n_{n-1})
(\langle\langle\gamma^1_0,\ldots,\gamma^1_{n-1}\rangle,\ldots \langle
\gamma^n_0,\ldots,\gamma^n_{n-1}\rangle \rangle\in B^*_{\eta_0})\\
  \end{array}\]
which means 
\[\begin{array}{ll}
\ &(\exists^{I^1_0}\gamma^1_0)\ldots\ldots(\exists^{I^n_{n-1}}\gamma^n_{n-1})\\
\ &(\exists\gamma^1_n)\ldots(\exists\gamma^n_n)(\langle\eta_0,\langle
\gamma^1_0,\ldots,\gamma^1_n\rangle,\ldots,\langle\gamma^n_0,\ldots,
\gamma^n_n\rangle\rangle\in B).\\
  \end{array}\]
By the assumptions $(\gamma)$ and $(\varepsilon)$ we know that
\[\begin{array}{ll}
\ &(\exists^{I^0_0}\gamma^0_0)\ldots(\exists^{I^0_n}\gamma^0_n)
(\exists^{I^1_0}\gamma^1_0)\ldots(\exists^{I^1_n}\gamma^1_n)\ldots\\
\ldots &(\exists^{I^n_0}\gamma^n_0)\ldots(\exists^{I^n_{n}}\gamma^n_{n})
(\langle\langle\gamma_0^0,\ldots,\gamma^0_n\rangle,\langle
\gamma^1_0,\ldots,\gamma^1_n\rangle,\ldots,\langle\gamma^n_0,\ldots,
\gamma^n_n\rangle\rangle\in B).\\
  \end{array}\]
Obviously any quantifier $(\exists^{I^m_\ell}\gamma^m_\ell)$ above may be
replaced by $(\exists\gamma^m_\ell)$ and then ``moved'' right as much as we
want. 
Consequently we get
\[\begin{array}{l}
(\exists\gamma^0_0)\ldots(\exists\gamma^0_{n-1})(\exists^{I^0_n}
\gamma^0_n)(\exists^{I^1_0}\gamma^1_0)\ldots(\exists^{I^1_{n-1}}
\gamma^1_{n-1})\ldots\ldots(\exists^{I^n_0}\gamma^n_0)\ldots
(\exists^{I^n_{n-1}}\gamma^n_{n-1})\\
(\exists\gamma^1_n)\ldots(\exists\gamma^n_n)(\langle\langle\gamma_0^0,
\ldots,\gamma^0_n\rangle,\langle\gamma^1_0,\ldots,\gamma^1_n\rangle,\ldots,
\langle\gamma^n_0,\ldots,\gamma^n_n\rangle\rangle\in B)\\
  \end{array}\]
which means that
\[(\exists\gamma^0_0)\ldots(\exists\gamma^0_{n-1})(\exists^{I^0_n}
\gamma^0_n)(B^*_{\langle\gamma^0_0,\ldots,\gamma^n_n\rangle}\mbox{ is big}).\]
Hence we find $\gamma^0_0,\ldots,\gamma^0_{n-1}$ and a set $a\in (I^0_n)^+$
such that
\[(\forall\gamma\in a)(B^*_{\langle\gamma^0_0,\ldots,\gamma^n_n\rangle}\mbox{
is big}).\]
Note that the assumptions of the proposition are such that if we know that
$B^*_{\eta_0}$ is big then we may apply the inductive hypothesis to
\[\lambda^m_\ell,\chi^m_\ell,P^m_\ell,I^m_\ell, J^m\ \mbox{ (for $1\leq m\leq
n$, $\ell\leq n-1$) and } B^*_{\eta_0}.\]
Consequently for each $\gamma\in a$ we find sets $X^\gamma_1,\ldots,
X^\gamma_n$ such that for $1\leq m\leq n$:
\begin{description}
\item[(a)$^*$] $X^\gamma_m\subseteq\prod\limits_{\ell\leq n-m}\lambda^m_\ell$
\item[(b)$^*$] if $\eta,\nu\in X^\gamma_m$, $\eta\neq\nu$ then
\begin{description}
\item[(i)\ ] $\eta\rest (n-m)=\nu\rest (n-m)$
\item[(ii) ] $\eta(n-m)\neq\nu (n-m)$
\end{description}
\item[(c)$^*$] $\{\eta(n-m):\eta\in X^\gamma_m\}\in P^m_{n-m}$\ \ \ \ \ \ \ \
\ \ and 
\item[(d)$^*$] for all $\langle\eta_0,\ldots,\eta_n\rangle\in\prod
\limits_{m\leq n} X^\gamma_m$ we have
\[(\exists\langle\eta^*_0,\ldots,\eta^*_n\rangle
\in B^*_{\langle\gamma^0_0,\ldots,\gamma^0_{n-1},\gamma\rangle})(\forall
1\leq m\leq n)(\nu_m\trianglelefteq\nu^*_m).\]
\end{description}
Now we may ask how mane possibilities for $X^\gamma_m$ do we have: not too
many. If we fix the common initial segment (see {\bf (b)$^*$}) the only
freedom we have is in choosing an element of $P^m_{n-m}$ (see {\bf (c)$^*$}).
Consequently there are at most $|P^m_{n-m}|+\sum\limits_{\ell\leq n-m}
\lambda^m_\ell$ possible values for $X^\gamma_m$ and hence there are at most
\[\sum_{m=1}^n(|P^m_{n-m}|+\sum_{\ell\leq n-m}\lambda^m_\ell)<\chi^0_n\]
possible values for the sequence $\langle
X^\gamma_1,\ldots,X^\gamma_n\rangle$. Since the ideal $I^0_n$ is
$\chi^0_n$-complete we find a sequence $\langle X_1,\ldots,X_n\rangle$ and a
set $b\subseteq a$, $b\in (I^0_n)^+$ such that 
\[(\forall\gamma\in b)(\langle X^\gamma_1,\ldots,X^\gamma_n\rangle = \langle
X_1,\ldots, X_n\rangle).\]
Next choose $b^0_n\in P^0_n$ such that $b^0_n\subseteq b$ and put
\[X_0=\{\langle\gamma^0_0,\ldots,\gamma^0_{n-1},\gamma\rangle: \gamma\in
b^0_n\}.\] 
Now it is a routine to check that the sets $X_0,X_1,\ldots,X_n$ are as
required (i.e. they satisfy clauses {\bf (a)}--{\bf (d)}). \QED

\section{Cofinal sequences in trees} 
\begin{notation}
\label{evenmorenot}
\begin{enumerate}
\item For a tree $T\subseteq {}^{\delta{>}}\mu$ the set of $\delta$-branches
through $T$ is
\[{\lim}_\delta(T)\stackrel{\rm def}{=}\{\eta\in{}^\delta\mu: (\forall\alpha
<\delta)(\eta\rest\alpha\in T)\}.\]
The $i$-th level (for $i<\delta$) of the tree $T$ is 
\[T_i\stackrel{\rm def}{=} T\cap {}^i\mu\]
and $T_{{<}i}\stackrel{\rm def}{=}\bigcup\limits_{j<i}T_j$.

If $\eta\in T$ then the set of immediate successors of $\eta$ in $T$ is
\[\suc_T\stackrel{\rm def}{=}\{\nu\in T: \eta\vartriangleleft\nu\ \&\
\lh(\nu)=\lh(\eta)+1\}.\] 
We shall not distinguish strictly between $\suc_T(\eta)$ and $\{\alpha:
\eta\hat{\ }\langle\alpha\rangle \in T\}$.
\end{enumerate}
\end{notation}

\begin{definition}
\label{cofinal}
\begin{enumerate}
\item $\K_{\mu,\delta}$ is the family of all pairs $(T,\bar{\lambda})$ such
that $T\subseteq{}^{\delta{>}}\mu$ is a tree with $\delta$ levels and
$\bar{\lambda}=\langle\lambda_\eta:\eta\in T\rangle$ is a sequence of
cardinals such that for each $\eta\in T$ we have $\suc_T(\eta)=\lambda_\eta$
(compare the previous remark about not distinguishing $\suc_T(\eta)$ and
$\{\alpha:\eta\hat{\ }\langle\alpha\rangle\in T\}$).

\item For a limit ordinal $\delta$ and a cardinal $\mu$  we let
\[\begin{array}{ll}
\K^{\id}_{\mu,\delta}\stackrel{\rm def}{=}\{(T,\bar{\lambda},\bar{I}):
&(T,\bar{\lambda})\in \K_{\mu,\delta}, \ \bar{I}=\langle I_\eta:\eta\in
T\rangle\\ 
\ &\mbox{each }I_\eta\mbox{ is an ideal on }\lambda_\eta=\suc_T(\eta)\}.\\
\end{array}\]

Let $(T,\bar{\lambda},\bar{I})\in\K_{\mu,\delta}^{\id}$ and let $J$ be an
ideal on $\delta$ (including $J^{\bd}_\delta$ if we do not say otherwise).
Further let
$\bar{\eta}=\langle\eta_\alpha:\alpha<\lambda\rangle\subseteq{\lim}_\delta(T)$ 
be a  sequence of $\delta$-branches through $T$.

\item We say that $\bar{\eta}$ is {\em $J$-cofinal in} $(T,\bar{\lambda},
\bar{I})$ if 
\begin{description}
\item[(a)] $\eta_\alpha\neq\eta_\beta$ for distinct $\alpha,\beta<\lambda$
\item[(b)] for every sequence $\bar{A}{=}\langle A_\eta\!:\eta\in T\rangle
\in\!\!\prod\limits_{\eta\in T} I_\eta$ there is $\alpha^*{<}\lambda$ such
that   
\[\alpha^*\leq\alpha<\lambda\ \ \ \Rightarrow\ \ \ (\forall^J i<\delta)
(\eta_\alpha\rest (i+1)\notin A_{\eta_\alpha\rest i}).\]
\end{description}

\item If $I$ is an ideal on $\lambda$ then we say that
$(\bar{\eta},I)$ is {\em a $J$-cofinal pair} for $(T,\bar{\lambda},
\bar{I})$ if  
\begin{description}
\item[(a)] $\eta_\alpha\neq\eta_\beta$ for distinct $\alpha,\beta<\lambda$
\item[(b)] for every sequence $\bar{A}=\langle A_\eta:\eta\in T\rangle\in
\prod\limits_{\eta\in T} I_\eta$ there is $A\in I$ such that 
\[\alpha\in\lambda\setminus A\ \ \ \Rightarrow\ \ \ (\forall^J i<\delta)
(\eta_\alpha\rest (i+1)\notin A_{\eta_\alpha\rest i}).\]
\end{description}

\item The sequence $\bar{\eta}$ is {\em strongly $J$-cofinal in}
$(T,\bar{\lambda}, \bar{I})$ if 
\begin{description}
\item[(a)] $\eta_\alpha\neq\eta_\beta$ for distinct $\alpha,\beta<\lambda$
\item[(b)] for every $n<\omega$ and functions $F_0,\ldots,F_{n}$ there is
$\alpha^*<\lambda$ such that  

\noindent {\em if} $m\leq n$, $\alpha_0<\ldots<\alpha_n<\lambda$,
$\alpha^*\leq\alpha_m$ 

\noindent {\em then} the set
\[\hspace{-1.3cm}\begin{array}{ll}
\{i<\delta:&\mbox{(i)\ \ } (\forall \ell<m)(\lambda_{\eta_{\alpha_\ell}\rest
i}< \lambda_{\eta_{\alpha_m}\rest i})\mbox{ and}\\
\ &\mbox{(ii) }F_m(\eta_{\alpha_0}\rest (i{+}1),\ldots,
\eta_{\alpha_{m-1}}\rest (i{+}1),\eta_{\alpha_m}\rest i ,\ldots,
\eta_{\alpha_n}\rest i)\in I_{\eta_{\alpha_m}\rest i}\\ 
\ &\ \ \mbox{(and well defined) but }\\
\ &\ \ \eta_{\alpha_m}\rest(i{+}1)\in F_m(\eta_{\alpha_0}\rest
(i{+}1),\ldots, \eta_{\alpha_{m-1}}\rest (i{+}1), \eta_{\alpha_m}\rest
i,\ldots, \eta_{\alpha_n}\rest i)\}\\
\end{array}\] 
is in the ideal $J$.
\end{description}
[Note: in {\bf (b)} above we may have $\alpha^*<\alpha_0$, this causes no real
change.] 
\item The sequence $\bar{\eta}$ is {\em stronger $J$-cofinal in}
$(T,\bar{\lambda}, \bar{I})$ if 
\begin{description}
\item[(a)] $\eta_\alpha\neq\eta_\beta$ for distinct $\alpha,\beta<\lambda$
\item[(b)] for every $n<\omega$ and functions $F_0,\ldots,F_{n}$ there
is $\alpha^*<\lambda$ such that  

\noindent {\em if} $m\leq n$, $\alpha_0<\ldots<\alpha_n<\lambda$,
$\alpha^*\leq\alpha_m$ 

\noindent {\em then} the set
\[\hspace{-1.3cm}\begin{array}{ll}
\{i<\delta:&\mbox{(ii) }F_m(\eta_{\alpha_0}\rest (i{+}1),\ldots,
\eta_{\alpha_{m-1}}\rest (i{+}1),\eta_{\alpha_m}\rest i ,\ldots,
\eta_{\alpha_n}\rest i)\in I_{\eta_{\alpha_m}\rest i}\\ 
\ &\ \ \mbox{(and well defined) but }\\
\ &\ \ \eta_{\alpha_m}\rest(i{+}1)\in F_m(\eta_{\alpha_0}\rest
(i{+}1),\ldots, \eta_{\alpha_{m-1}}\rest (i{+}1), \eta_{\alpha_m}\rest
i,\ldots,\eta_{\alpha_n}\rest i)\} 
   \end{array}\]
is in the ideal $J$.
\end{description}

\item The sequence $\bar{\eta}$ is {\em strongest $J$-cofinal in}
$(T,\bar{\lambda}, \bar{I})$ if 
\begin{description}
\item[(a)] $\eta_\alpha\neq\eta_\beta$ for distinct $\alpha,\beta<\lambda$
\item[(b)] for every $n<\omega$ and functions $F_0,\dots,F_n$ there is
$\alpha^*<\lambda$ such that  

\noindent {\em if} $m\leq n$, $\alpha_0<\ldots<\alpha_n<\lambda$,
$\alpha^*\leq\alpha_m$ 

\noindent {\em then} the set
\[\hspace{-1.8cm}\begin{array}{ll}
\{i<\delta:&\mbox{(i')\ \ } (\exists \ell<m)(\lambda_{\eta_{\alpha_\ell}\rest
i}\geq \lambda_{\eta_{\alpha_m}\rest i})\mbox{ or}\\
\ &\mbox{(ii) }F_m(\eta_{\alpha_0}\rest (i{+}1),\ldots,
\eta_{\alpha_{m-1}}\rest (i{+}1),\eta_{\alpha_m}\rest i ,\ldots,
\eta_{\alpha_n}\rest i)\in I_{\eta_{\alpha_m}\rest i}\\ 
\ &\ \ \mbox{(and well defined) but }\\
\ &\ \ \eta_{\alpha_m}\rest(i{+}1)\in F_m(\eta_{\alpha_0}\rest
(i{+}1),\ldots, \eta_{\alpha_{m-1}}\rest (i{+}1), \eta_{\alpha_m}\rest
i,\ldots, \eta_{\alpha_n}\rest i)\}\\
 \end{array}\] 
is in the ideal $J$.
\item The sequence $\bar \theta$ is {\em big$^*$ $J$-cofinal} in
$(T, \bar \lambda, \bar I)$ if
\begin{description}
\item[(a)] $\eta_\alpha\neq \eta_\beta$ for distinct $\alpha$,
$\beta<\lambda$
\item[(b)] for every $\eta$ and functions $F_0, \ldots, F_n$ there is
$\alpha^*$ such that 

\noindent if $\alpha_0<\ldots <\alpha_n$ and $\alpha^*\leq \alpha_m$ then for
$m\leq n$ the set
\[\hspace{-1.8cm}\begin{array}{ll}
\{i<\delta: & \mbox{ if }\nu_\ell=\left\{
			\begin{array}{ll}
	\eta_{\alpha_\ell}\rest (i{+}1) & \mbox{ if }
		\lambda_{\eta_{\alpha_\ell}\restriction i}=
		\lambda_{\eta_{\alpha-m}\rest i}\ \mbox{ or}\\
	\ 		&\lambda_{\eta_{\alpha_\ell}\rest i}=
		\lambda_{\eta_{\alpha_n}\rest i}\ \mbox{ and }\
		\eta_{\alpha_\ell}(i)< \eta_{\alpha_m}(i)\\
	\eta_{\alpha_\ell}\restriction i & \mbox{ if not}\\
			\end{array}
	\right.\\
\ &\mbox{ then we have }\\
\ &\eta_{\alpha_m}(i)\in F_m(\nu_\ell)\in I_{\eta_{\alpha_m}\restriction
i}\}\\  
\end{array}\]
is in the ideal $J$.
\end{description}
\item We define ``strongly$^*$ $J$-cofinal'', ``stronger$^*$ $J$-cofinal'' and
``strongest$^*$ big $J$-cofinal'' in almost the same way, replacing the
requirement that $\alpha^*\leq \alpha_m$ in 5{\bf (b)}, 6{\bf (b)}, 7{\bf (b)}
above (respectively) by $\alpha^*\leq\alpha_0$.
\end{description}
\end{enumerate}
\end{definition}

\begin{remark}
\label{impcof}
{\em
\begin{enumerate}
\item Note that ``strongest $J$-cofinal'' implies ``stronger $J$-cofinal''
and this implies ``strongly $J$-cofinal''. ``Stronger $J$-cofinal'' 
implies ``$J$-cofinal''. Also ``bigger'' $\Rightarrow$ ``big''
$\Rightarrow$ ``cofinal'', ``big'' $\Rightarrow$ ``strongly''.
\item The different notions of ``strong $J$-cofinality'' (the conditions {\em
(i)} and {\em (i')}) are to allow us to carry some diagonalization arguments. 
\item The difference between ``strongly $J$-cofinal'' and ``strongly$^*$
$J$-cofinal'' etc is, in our context, immaterial. we may in all places in this
paper replace the respective notion with its version with ``$*$'' and no harm
will be done. 
\end{enumerate}
}
\end{remark}

\begin{remark}
{\em
\begin{enumerate}
\item Remind {\bf pcf}:

An important case is when $\langle\lambda_i: i<\delta\rangle$ is an
increasing sequence of regular cardinals, $\lambda_i>\prod\limits_{j<i}
\lambda_j$, $\lambda_\eta=\lambda_{\lh(\eta)}$,
$I_\eta=J^{\bd}_{\lambda_\eta}$ and $\lambda=
\tcf(\prod\limits_{i<\delta}\lambda_i/J)$. 

\item Moreover we are interested in more complicated $I_\eta$'s (as in
\cite{Sh:430}, section 5), connected to our problem, so ``the existence of the
true cofinality'' is less clear. But the assumption $2^\mu=\mu^+$ will rescue
us. 

\item There are natural stronger demands of cofinality since here we are not
interested just in $x_\alpha$'s but also in Boolean combinations. Thus
naturally we are interested in behaviours of large sets of $n$-tuples, see
\ref{super}. 
\end{enumerate}
}
\end{remark}

\begin{proposition}
\label{strimpstr}
Suppose that $(T,\bar{\lambda},\bar{I})\in\K^{\id}_{\mu,\delta}$,
$\bar{\eta}=\langle\eta_\alpha:\alpha<\lambda\rangle\subseteq{\lim}_\delta
(T)$ and $J$ is an ideal on $\delta$, $J\supseteq J^{\bd}_\delta$. 
\begin{enumerate}
\item Assume that 
\begin{description}
\item[$(\circledcirc)$] if $\alpha<\beta<\lambda$ then $(\forall^J i<\delta)
(\lambda_{\eta_\alpha\rest i}<\lambda_{\eta_\beta\rest i})$.
\end{description}
Then the following are equivalent
\[\mbox{``}\bar{\eta}\mbox{ is strongly $J$-cofinal for }(T,\bar{\lambda},
\bar{I})\mbox{''}\] 
\[\mbox{``}\bar{\eta}\mbox{ is stronger $J$-cofinal for }(T,\bar{\lambda},
\bar{I})\mbox{''}\] 
\[\mbox{``}\bar{\eta}\mbox{ is strongest $J$-cofinal for }(T,\bar{\lambda},
\bar{I})\mbox{''.}\]
\[\mbox{``}\bar{\eta}\mbox{ is big $J$-cofinal for }(T,\bar{\lambda},
\bar{I})\mbox{''.}\]
\item If $I_\nu\supseteq J^{\bd}_{\lambda_\nu}$ and $\lambda_\nu
=\lambda_{\lh(\nu)}$ for each $\nu\in T$ and the
sequence $\bar{\eta}$ is stronger $J$-cofinal for $(T,\bar{\lambda},\bar{I})$
then for some $\alpha^*<\lambda$ the sequence $\langle\eta_\alpha:\alpha^*\leq
\alpha<\lambda\rangle$ is $<_J$-increasing. 
\item If $\eta\in T_i \Rightarrow \lambda_\eta=\lambda_i$ and $\bar
\eta <_J$-increasing in $\prod\limits_{i<\delta}$ {\em then}
``big'' is equivalent to ``stronger''. \QED
\end{enumerate}
\end{proposition}

\begin{proposition}
\label{getcofinal}
Suppose that
\begin{enumerate}
\item $\langle\lambda_i: i<\delta\rangle$ is an increasing sequence of regular
cardinals, $\delta<\lambda_0$ is a limit ordinal, 
\item $T=\bigcup\limits_{i<\delta}\prod\limits_{j<i}\lambda_j$,
$I_\eta=I_{\lh(\eta)}=J^{\bd}_{\lambda_{\lh(\eta)}}$,
$\lambda_\eta=\lambda_{\lh(\eta)}$, 
\item $J$ is an ideal on $\delta$,
$\lambda=\tcf(\prod\limits_{i<\delta}\lambda_i/J)$ and it is examplified 
by a sequence $\bar{\eta}=\langle\eta_\alpha:\alpha<\lambda\rangle\subseteq
\prod\limits_{i<\delta}\lambda_i$
\item for each $i<\delta$
\[|\{\eta_\alpha\rest i: \alpha<\lambda\}|<\lambda_i\]
(so e.g. $\lambda_i>\prod\limits_{j<i}\lambda_j$ suffices).
\end{enumerate}
Then the sequence $\bar{\eta}$ is $J$-cofinal in $(T,\bar{\lambda},\bar{I})$.
\end{proposition}

\Proof First note that our assumptions imply that each ideal
$I_\eta=I_{\lh(\eta)}$ is $|\{\eta_\alpha\rest \lh(\eta):
\alpha<\lambda\}|^+$-complete. Hence for  
each sequence $\bar{A}=\langle A_\eta:\eta\in T\rangle\in
\prod\limits_{\eta\in T} I_\eta$ and $i<\delta$ the set
\[A_i\stackrel{\rm def}{=}\bigcup\{A_{\eta_\alpha\rest i}: \alpha<\lambda\}\]
is in the the ideal $I_i$, i.e. it is bounded in $\lambda_i$ (for $i<\delta$).
(We should remind here our convention which says in this case that we do not
distinguish $\lambda_i$ and $\suc_T(\eta)$ if $\lh(\eta)=i$,
see \ref{evenmorenot}.) Take $\eta^*\in\prod\limits_{i<\delta}\lambda_i$
such that for each $i<\delta$ we have $A_i\subseteq\eta^*(i)$. As the sequence
$\bar{\eta}$ realizes the true cofinality of
$\prod\limits_{i<\delta}\lambda_i/J$ we find $\alpha^*<\lambda$ such that 
\[\alpha^*\leq\alpha<\lambda\ \ \ \Rightarrow\ \ \ \{i<\delta: \eta_\alpha(i)
<\eta^*(i)\}\in J\]
which allows us to finish the proof. \QED
\medskip

It follows from the above proposition that the notion of $J$-cofinal sequences
is not empty. Of course it is better to have ``strongly (or even: stronger)
$J$-cofinal'' sequences $\bar{\eta}$. So it is nice to have that sometimes the 
weaker notion implies the stronger one.

\begin{proposition}
\label{cofimpstrong}
Assume that $\delta$ is a limit ordinal, $\mu$ is a cardinal,
$(T,\bar{\lambda},\bar{I})\in\K^{\id}_{\mu,\delta}$. Let $J$ be an ideal on
$\delta$ such that $J\supseteq J^{\bd}_\delta$ (which is our standard
hypothesis). Further suppose that
\begin{description}
\item[$(\circledast)$] if $\eta\in T_i$ then the ideal $I_\eta$ is
$(|T_i|+\sum\{\lambda_\nu\!: \nu{\in} T_i\ \&\ \lambda_\nu{<}\lambda_\eta
\})^+$--complete.   
\end{description}
{\em Then} each $J$-cofinal sequence $\bar{\eta}$ for
$(T,\bar{\lambda},\bar{I})$ is strongly $J$-cofinal for
$(T,\bar{\lambda},\bar{I})$. 

If in addition $\eta\neq\nu\in T_i\ \Rightarrow\ \lambda_\eta\neq\lambda_\nu$
then $\bar \eta$ is big $J$-cofinal for $(T,\bar{\lambda},\bar{I})$. Also if
in addition  
\[\eta\in T_i\ \ \Rightarrow\ \ (\exists^{!1}\nu\in
T_i)(\lambda_\nu=\lambda_\eta)\ \vee\ [(\exists^{\leq \lambda_\eta} 
\nu\in T_i)(\lambda_\nu =\lambda_\eta)\ \&\ I_\eta\mbox{ normal}]\]
then $\bar \eta$ is big $J$-cofinal.
\end{proposition}

\Proof Let $n<\omega$ and $F_0,\ldots,F_n$ be $(n+1)$-place functions. First
we define a sequence $\bar{A}=\langle A_\eta: \eta\in T\rangle$. For $m\leq n$
and a sequence $\langle\eta_m,\ldots,\eta_n\rangle\subseteq T_i$ we put 
\[\begin{array}{ll}
A^m_{\langle\eta_m,\ldots,\eta_n\rangle}{=}\bigcup\{& F_m(\nu_0,\ldots,
\nu_{m-1},\eta_m,\ldots,\eta_n):\nu_0,\ldots,\nu_{m-1}\in T_{i+1},\\
\ &\quad\quad(\nu_0,\ldots,\nu_{m-1},\eta_m,\ldots,\eta_n)\!\in\!\dom(F),\\ 
\ &\quad\quad\lambda_{\nu_0\rest i}<\lambda_\eta,\ldots,\lambda_{\nu_{m-1}
\rest i}<\lambda_{\eta_m}\\
\ &\quad\quad\mbox{and }F(\nu_0,\ldots,\nu_{m-1},\eta_m,\ldots,\eta_n)\in
I_{\eta_m}\} 
  \end{array}\]
and next for $\eta\in T_i$ let 
\[A_\eta=\bigcup\{A^m_{\langle\eta,\eta_{m+1},\ldots,\eta_n\rangle}: m\leq n\
\&\ \eta_{m+1},\ldots,\eta_n\in T_i\}.\]
Note that the assumption $(\circledast)$ was set up so that
$A^m_{\langle\eta_m,\ldots,\eta_n\rangle}\in I_{\eta_m}$ and the sets $A_\eta$
are in $I_\eta$ (for $\eta\in T$).

By the $J$-cofinality of $\bar{\eta}$, for some $\alpha^*<\lambda$ we have
\[\alpha^*\leq\alpha<\lambda\ \ \ \Rightarrow\ \ \ (\forall^J i<\delta)
(\eta_\alpha\rest(i+1)\notin A_{\eta_\alpha\rest i}).\]
We are going to prove that this $\alpha^*$ is as required in the definition of
strongly $J$-cofinal sequences. So suppose that $m\leq n$, $\alpha_0<\ldots
<\alpha_n<\lambda$ and $\alpha^*\leq\alpha_m$. By the choice of $\alpha^*$ we
have that the set $A\stackrel{\rm def}{=}\{i<\delta: \eta_{\alpha_m}
\rest(i+1)\in A_{\eta_{\alpha_m}\rest i}\}$ is in the ideal $J$. But if
$i<\delta$ is such 
that  
\begin{quotation}
\noindent $(\forall \ell<m)(\lambda_{\eta_{\alpha_\ell}\rest i}<
\lambda_{\eta_{\alpha_m}\rest i})$ and

\noindent $F(\eta_{\alpha_0}\rest(i+1),\ldots,
\eta_{\alpha_{m-1}}\rest (i+1),\eta_{\alpha_m}\rest i,\ldots,\eta_{\alpha_n}
\rest i)\in I_{\eta_{\alpha_m}\rest i}$ but

\noindent $\eta_{\alpha_m}\rest (i+1)\in F(\eta_{\alpha_0}\rest(i+1),\ldots,
\eta_{\alpha_{m-1}}\rest (i+1),\eta_{\alpha_m}\rest i,\ldots,\eta_{\alpha_n}
\rest i)$ 
\end{quotation}
then clearly $\eta_{\alpha_m}\rest (i+1)\in A^m_{\langle\eta_{\alpha_m}\rest
i,\ldots,\eta_{\alpha_n}\rest i\rangle}$ and so $i\in A$. This finishes the
proof. 

The ``big'' verion should be clear too. \QED 

\begin{proposition}
\label{stronger}
Assume that $\mu$ is a strong limit uncountable cardinal and
$\langle\mu_i:i<\delta \rangle$ is an increasing sequence of cardinals with
limit $\mu$. Further suppose that $(T,\bar{\lambda},\bar{I})\in
\K^{\id}_{\mu,\delta}$, $|T_i|\leq\mu_i$ (for $i<\delta$), $\lambda_\eta<\mu$
and each $I_\eta$ is $\mu_{\lh(\eta)}^+$-complete and contains all singletons
(for $\eta\in T$). Finally assume 
$2^\mu=\mu^+$ and let $J$ be an ideal on $\delta$, $J\supseteq
J^{\bd}_\delta$.  

\noindent Then there exists a stronger $J$-cofinal sequence $\bar{\eta}$ for
$(T,\bar{\lambda},\bar{I})$ of the length $\mu^+$ (even for
$J=J^{\bd}_\delta$). 

\noindent We can get ``big'' if
$$
\rho\neq \eta\in T_i\ \&\ \lambda_\rho=\lambda_\eta \Rightarrow
(\exists^{\leq \lambda_\eta}\nu\in T_i)(\lambda_\nu=\lambda_\eta)\ \&\
I_\eta\mbox{ normal.}
$$
\end{proposition}
\Proof This is a straight diagonal argument. Put
\[\begin{array}{ll}
Y\stackrel{\rm def}{=}\{\langle F_0,\ldots,F_n\rangle: &n<\omega \mbox{ and
each }F_\ell\mbox{ is a function with}\\  
\ &\dom(F)\subseteq T^{n+1},
\rng(F)\subseteq\bigcup\limits_{\eta\in T} I_\eta\}.\\ 
  \end{array}\] 
Since $|Y|=\mu^\mu=\mu^+$ (remember that $\mu$ is strong limit and
$\lambda_\eta<\mu$ for $\eta\in T$) we may choose an enumeration $Y=\{\langle
F^\xi_0,\ldots,F^\xi_{n_\xi}\rangle: \xi<\mu^+\}$. For each $\zeta<\mu^+$
choose an increasing sequence $\langle\A^\zeta_i: i<\delta\rangle$ such that
$|\A^\zeta_i|\leq\mu_i$ and $\zeta=\bigcup\limits_{i<\delta}\A^\zeta_i$. Now
we choose by induction on $\zeta<\mu^+$ branches $\eta_\zeta$ such that for
each $\zeta$ the restriction $\eta_\zeta\rest i$ is defined by induction on
$\zeta$ as follows.  

\noindent If $i=0$ or $i$ is limit then there is nothing to do.

\noindent Suppose now that we have defined $\eta_\zeta\rest i$ and $\eta_\xi$
for $\xi<\zeta$. We find $\eta_\zeta(i)$ such that 
\begin{description}
\item[$(\alpha)$] $\eta_\zeta(i)\in\lambda_{\eta_\zeta\rest i}$
\item[$(\beta)$] if $\varepsilon\in\A^\zeta_i$, $m\leq n_\varepsilon$,
$\alpha_0,\ldots,\alpha_{m-1}\in \A^\zeta_i$ (hence $\alpha_\ell<\zeta$ so
$\eta_{\alpha_\ell}$ are defined already), $\nu_{m+1},\ldots,\nu_n\in T_i$ and
\[\hspace{-1cm}F^\varepsilon_m(\eta_{\alpha_0}\rest (i+1),\ldots,
\eta_{\alpha_{m-1}}\rest (i+1), \eta_\zeta\rest i, \nu_{m+1},\ldots,\nu_n) \in
I_{\eta_\zeta\rest i}\quad\mbox{and well defined}\] 
then 
\[\eta_\zeta\rest (i+1)\notin F_m^\varepsilon(\eta_{\alpha_0}\rest
(i+1),\ldots, \eta_{\alpha_{m-1}}\rest (i+1),\eta_\zeta\rest i,
\nu_{m+1},\ldots,\nu_n)\] 
\item[$(\gamma)$] $\eta_\zeta\rest(i+1)\notin\{\eta_\varepsilon\rest (i+1):
\varepsilon\in\A^\zeta_i\}$.
\end{description}
{\em Why it is possible?\/} Note that there is $\leq\aleph_0+|\A^\zeta_i|+
|\A^\zeta_i|^{\textstyle {<}\aleph_0}+|T_i|\leq\mu_i$ negative demands 
and each of them says that $\eta_\zeta\rest (i+1)$ is not in some set from
$I_{\eta_\zeta\rest i}$ (remember that we have assumed that the ideals
$I_{\eta_\zeta\rest i}$ contain singletons). Consequently using the
completeness of the ideal we may satisfy the requirements
$(\alpha)$--$(\gamma)$ above. 

Now of course $\eta_\zeta\in{\lim}_\delta(T)$. Moreover if
$\varepsilon<\zeta<\mu^+$ then $(\exists i<\delta)(\varepsilon\in \A^\zeta_i)$
which implies $(\exists i<\delta)(\eta_\varepsilon\rest (i+1)\neq
\eta_\zeta\rest (i+1))$. Consequently
\[\varepsilon<\zeta<\mu^+\ \ \ \Rightarrow\ \ \
\eta_\varepsilon\neq\eta_\zeta.\] 
Checking the demand {\bf (b)} of ``stronger $J$-cofinal'' is straightforward:
for functions $F_0,\ldots,F_n$ (and $n\in\omega$) take $\varepsilon$ such that
\[\langle F_0,\ldots,F_n\rangle = \langle F^\varepsilon_0,\ldots,
F^\varepsilon_{n_\varepsilon}\rangle\]
and put $\alpha^*=\varepsilon+1$. Suppose now that $m\leq n$,
$\alpha_0<\ldots<\alpha_n<\lambda$, $\alpha^*\leq\alpha_m$. Let $i^*<\delta$
be such that for $i>i^*$ we have
\[\varepsilon,\alpha_0,\ldots,\alpha_{m-1}\in \A^{\alpha_m}_i.\]
Then by the choice of $\eta_{\alpha_m}\rest(i+1)$ we have that for each
$i>i^*$ 
\smallskip

if $F_m^\varepsilon(\eta_{\alpha_0}\rest (i+1),\ldots,
\eta_{\alpha_{m-1}}\rest (i+1),\eta_\zeta\rest i, \eta_{\alpha_{m+1}}\rest i,
\ldots,\eta_{\alpha_n}\rest i)\in I_{\eta_{\alpha_m}\rest i}$

then $\eta_{\alpha_m}\rest i\notin F_m^\varepsilon(\eta_{\alpha_0}\rest
(i+1),\ldots, \eta_{\alpha_{m-1}}\rest (i+1),\eta_\zeta\rest i,
\eta_{\alpha_{m+1}}\rest i, \ldots,\eta_{\alpha_n}\rest i)$.
\smallskip

\noindent This finishes the proof. \QED 

\begin{remark}
{\em
The proof above can be carried for functions $F$ which depend on
$(\eta_{\alpha_0},\ldots,\eta_{\alpha_{m-1}},\eta_{\alpha_m}\rest i,\ldots,
\eta_{\alpha_n}\rest i)$. This will be natural later.
}
\end{remark}

Let us note that if the ideals $I_\eta$ are sufficiently complete then
$J$-cofinal sequences cannot be too short.

\begin{proposition}
Suppose that $(T,\bar{\lambda},\bar{I})\in\K^{\id}_{\mu,\delta}$ is such that
for each $\eta\in T_i$, $i<\delta$ the ideal $I_\eta$ is
$(\kappa_i)^+$--complete (enough if $[\lambda_\eta]^{\kappa_i}
\subseteq I_\eta$). Let $J\supseteq J^{\bd}_\delta$ be an ideal on $\delta$
and let $\bar{\eta}=\langle\eta_\alpha:\alpha<\delta^*\rangle$ be a
$J$-cofinal sequence for $(T,\bar{\lambda},\bar{I})$. Then  
\[\delta^*>\lim\sup_J \kappa_i\]
and consequently 
\[\cf(\delta^*)>\lim\sup_J \kappa_i.\]
\end{proposition}

\Proof Fix an enumeration $\delta^*=\{\alpha_\varepsilon: \varepsilon<
|\delta^*|\}$ and for $\alpha<\delta^*$ let $\zeta(\alpha)$ be the unique
$\zeta$ such that $\alpha=\alpha_\zeta$. 

\noindent For $\eta\in T_i$, $i<\delta$ put
\[A_\eta\stackrel{\rm def}{=}\{\nu\in\suc_T(\eta): (\exists\varepsilon\leq
\kappa_i)(\nu\vartriangleleft\eta_\varepsilon)\}.\] 
Clearly $|A_\eta|\leq \kappa_i$ and hence $A_\eta\in I_\eta$. Apply the
$J$-cofinality of $\bar{\eta}$ to the sequence $\bar{A}=\langle A_\eta:
\eta\in T\rangle$. Thus there is $\alpha^*<\delta^*$ such that for each
$\alpha\in [\alpha^*,\delta^*)$ we have
\[(\forall^J i<\delta)(\eta_\alpha\rest (i+1)\notin A_{\eta_\alpha\rest i})\]
and hence
\[(\forall^J i<\delta)(\zeta(\alpha)>\kappa_i)\]
and consequently
\[\zeta(\alpha)\geq \lim\sup_J \kappa_i.\]
Hence we conclude that $|\delta^*|>\lim\sup_J \kappa_i$.

\noindent For the part ``consequently'' of the proposition note that if
$\langle\eta_\alpha: \alpha<\delta^*\rangle$ is $J$-cofinal (in
$(T,\bar{\lambda},\bar{I})$) and $A\subseteq\delta^*$ is cofinal in $\delta^*$
then $\langle\eta_\alpha: \alpha\in A\rangle$ is $J$-cofinal too. \QED

\begin{remark}
{\em
\begin{enumerate}
\item So if we have a $J$-cofinal sequence of the length $\delta^*$ then we
also have one of the length $\cf(\delta^*)$. Thus assuming regularity of the
length is natural.
\item Moreover the assumption that the length of the sequence is above
$|\delta|+|T|$ is very natural and in most cases it will follow from the
$J$-cofinality (and completeness assumptions). However we will try to state
this condition in the assumptions whenever it is used in the proof (even if
it can be concluded from the other assumptions).
\end{enumerate}
}
\end{remark}

\section{Getting $(\kappa,{\rm not}\lambda)$-Knaster algebras}

\begin{proposition}
\label{algide}
Let $\lambda,\sigma$ be cardinals such that $(\forall\alpha<\sigma)
(2^{|\alpha|}<\lambda)$, $\sigma$ is regular. Then there are a Boolean algebra
$\ba$, a sequence $\langle y_\alpha:\alpha<\lambda\rangle\subseteq\ba^+$ and
an ideal $I$ on $\lambda$ such that
\begin{description}
\item[(a)] if $X\subseteq\lambda$, $X\notin I$ then $(\exists
\alpha,\beta\in X)(\ba\models y_\alpha\cap y_\beta=\bz)$
\item[(b)] the ideal $I$ is $\sigma$-complete
\item[(c)] the algebra $\ba$ satisfies the $\mu$-Knaster condition for any
regular uncountable $\mu$ (really $\ba$ is free).
\end{description}
\end{proposition}

\Proof Let $\ba$ be the Boolean algebra freely generated by
$\{z_\alpha:\alpha<\lambda\}$ (so the demand {\bf (c)} is satisfied). Let
$A=\{(\alpha,\beta):\alpha<\beta<\lambda\}$ and
$y_{(\alpha,\beta)}=z_\alpha-z_\beta(\neq\bz)$ (for $(\alpha,\beta)\in A$).
The ideal $I$ of subsets of $A$ defined by
\begin{quotation}
\noindent a set $X\subseteq A$ is in $I$ if and only if

\noindent there are $\zeta<\sigma$, $X_\varepsilon\subseteq A$ (for
$\varepsilon<\zeta$) such that $X\subseteq\bigcup\limits_{\varepsilon<\zeta}
X_\varepsilon$ and for every $\varepsilon<\zeta$ no two
$y_{(\alpha_1,\beta_1)}, y_{(\alpha_2,\beta_2)}\in X_\varepsilon$ are disjoint
in $\ba$.
\end{quotation}
First note that
\begin{claim}
$A\notin I$.
\end{claim}

\noindent{\em Why?}\ \ \ If not then we have witnesses $\zeta<\sigma$ and
$X_\varepsilon$ (for $\varepsilon<\zeta$) for it. 
So $A=\bigcup\limits_{\varepsilon<\zeta}X_\varepsilon$ and hence for
$(\alpha,\beta)\in A$ we have $\varepsilon(\alpha,\beta)$ such that
$y_{(\alpha,\beta)}\in X_{\varepsilon(\alpha,\beta)}$. So
$\varepsilon(\cdot,\cdot)$ is actually a function from $[\lambda]^2$ to
$\zeta<\sigma$. By the Erd\"os--Rado theorem we find $\alpha<\beta<\gamma<
\lambda$ such that $\varepsilon(\alpha,\beta)=\varepsilon(\beta,\gamma)$. But
\[y_{(\alpha,\beta)}\cap y_{(\beta,\gamma)}=(z_\alpha-z_\beta)\cap
(z_\beta-z_\gamma)=\bz\]
so $(\alpha,\beta)$, $(\beta,\gamma)$ cannot be in the same $X_\varepsilon$
-- a contradiction. 
\medskip

To finish the proof note that $I$ is $\sigma$-complete (as $\sigma$ is
regular), if $X\notin I$ then, by the definition of $I$, there are two
disjoint elements in $\{y_{(\alpha,\beta)}: (\alpha,\beta)\in X\}$. Finally
$|A|=\lambda$. \QED 

\begin{definition}
\label{well1}
\begin{description}
\item[(a)] A pair $(\ba,\bar{y})$ is called {\em a $\lambda$-marked Boolean
algebra} if $\ba$ is a Boolean algebra and $\bar{y}=\langle
y_\alpha:\alpha<\lambda\rangle$ is a sequence of non-zero elements of $\ba$.
\item[(b)] A triple $(\ba,\bar{y},I)$ is called a {\em $(\lambda,\chi)$-well
marked Boolean algebra} if $(\ba,\bar{y})$ is a $\lambda$-marked Boolean
algebra, $\chi$ is a regular cardinal and $I$ is a (proper) $\chi$-complete
ideal on $\lambda$ such that 
\[\{A\subseteq\lambda: (\forall \alpha,\beta\in A)(\ba\models y_\alpha\cap
y_\beta\neq \bz)\}\subseteq I.\]
By {\em $\lambda$-well marked Boolean algebra} we will mean
$(\lambda,\aleph_0)$-well marked one. As in the above situation $\lambda$ can
be read from $\bar{y}$ (as $\lambda=\lh(\bar{y})$) we may omit it and then we
may speak just about well marked Boolean algebras. 
\end{description}
\end{definition}

\begin{remark}
\label{well2}
{\em
Thus proposition \ref{algide} says that if $\lambda,\sigma$ are regular
cardinals and
\[(\forall \alpha<\sigma)(2^{|\alpha|}<\lambda)\]
then there exists a $(\lambda,\sigma)$-well marked Boolean algebra
$(\ba,\bar{y},I)$ such that $\ba$ satisfies the $\kappa$-Knaster
property for every $\kappa$. 
}
\end{remark}

\begin{definition}
\label{constructor}
\begin{description}
\item[(a)] For cardinals $\mu$ and $\lambda$ and a limit ordinal $\delta$,
{\em a $(\delta,\mu,\lambda)$-constructor} is a system 
\[\C=(T, \bar{\lambda}, \bar{\eta},\langle (\ba_\eta,\bar{y}_\eta):\eta\in
T\rangle)\]   
such that 
\begin{enumerate}
\item $(T,\bar{\lambda})\in \K_{\mu,\delta}$, 
\item $\bar{\eta}= \langle \eta_i: i\in \lambda\rangle$ where
$\eta_i\in\lim_\delta(T)$ (for $i<\lambda$) are distinct $\delta$-branches
through $T$ and 
\item for each $\eta\in T$:\ \ \ $(\ba_\eta,\bar{y}_\eta)$ is a
$\lambda_\eta$-marked Boolean algebra, i.e. $\bar{y}_\eta=\langle
y_{\eta\hat{\ }\langle\alpha\rangle}:\alpha< \lambda_\eta\rangle\subseteq
\ba_\eta^+$ (usually this will be an enumeration of $\ba_\eta^+$). 
\end{enumerate}
\item[(b)] Let $\C$ be a constructor (as above). We define Boolean algebras
$\ba_2=\ba^{\red}=\ba^{\red}(\C)$ and $\ba_1=\ba^{\green}=\ba^{\green}(\C)$
by:
\medskip

\noindent $\ba^{\red}$ is the Boolean algebra freely generated by
$\{x_i:i<\lambda\}$ {\em except that}
\begin{quotation}
\noindent {\em if} $i_0,\ldots,i_{n-1}<\lambda$, $\nu=\eta_{i_0}\rest
\zeta =\eta_{i_1}\rest\zeta=\ldots=\eta_{i_{n-1}}\rest\zeta$ and
$\ba_\nu\models\bigcap\limits_{\ell<n} y_{\eta_{i_\ell}\rest(\zeta+1)}
=\bz$ 

\noindent {\em then} $\bigcap\limits_{\ell<n} x_{i_\ell}=\bz$
\end{quotation}
[Note: we may demand that the sequence $\langle\eta_{i_\ell}(\zeta):\ell<n
\rangle$ is strictly increasing, this will cause no difference.]
\medskip

\noindent $\ba^{\green}$ is the Boolean algebra freely generated by
$\{x_i:i<\lambda\}$ {\em except that}
\begin{quotation}
\noindent {\em if} $\nu=\eta_i\rest\zeta =\eta_j\rest\zeta$,
$\eta_i(\zeta)\neq\eta_j(\zeta)$ and 

\noindent $\ba_\nu\models y_{\eta_i\rest(\zeta+1)}\cap
y_{\eta_j\rest(\zeta+1)}\neq\bz$  

\noindent {\em then} $x_i\cap x_j=\bz$.
\end{quotation}
\end{description}
\end{definition}

\begin{remark}
{\em
\begin{enumerate}
\item The equations for the green case can look strange but they have to be
dual to the ones of the red case.
\item ``Freely generated except $\ldots$'' means that a Boolean combination is
non-zero except when some (finitely many) conditions implies it. For this it
is enough to look at elements of the form 
\[x_{i_0}^{\gt_0}\cap\ldots \cap x_{i_{n-1}}^{\gt_{n-1}}\]
where $\gt_\ell\in\{0,1\}$.
\item Working in the free product $\ba^{\red}*\ba^{\green}$ we will use
the same notation for elements (e.g. generators) of $\ba^{\red}$ as for
elements of $\ba^{\green}$. Thus $x_i$ may stay either for the respective
generator in $\ba^{\red}$ or $\ba^{\green}$. We hope that this will not be
confusing, as one can easily decide in which algebra the element is considered
from the place of it (if $x\in\ba^{\red}$, $y\in\ba^{\green}$ then $(x,y)$
will stay for the element $x\cap_{\ba^{\red}*\ba^{\green}} y\in\ba^{\red}*
\ba^{\green}$). In particular we may write $(x_i,x_i)$ for an element
which could be denoted $x_i^{\red}\cap x_i^{\green}$.
\end{enumerate}
}
\end{remark}

\begin{remark}
{\em
If the pair $(\ba^{\red},\ba^{\green})$ is a counterexample with the free
product $\ba^{\red} *\ba^{\green}$ failing the $\lambda$-cc but each of the
algebras satisfying that condition then each of the algebras fails the
$\lambda$-Knaster condition. But $\ba^{\red}$ is supposed to have $\kappa$-cc
($\kappa$ smaller than $\lambda$). This is known to restrict $\lambda$. 
}
\end{remark}

\begin{proposition}
\label{notknaster}
Assume that $\C=(T,\bar{\lambda},\bar{\eta},\langle(\ba_\eta,\bar{y}_\eta):
\eta\in T\rangle)$ is a $(\delta,\mu,\lambda)$-constructor and $J\supseteq
J^{\bd}_\delta$ is an ideal on $\delta$ such that 
\begin{description}
\item[(a)] $\bar{\eta}=\langle \eta_i: i\in T\rangle$ is $J$-cofinal for
$(T,\bar{\lambda},\bar{I})$ 
\item[(b)] if $X\in I^+_\eta$ then 
\[(\exists\alpha,\beta\in X)(\ba_\eta\models y_{\eta\hat{\
}\langle\alpha\rangle}\cap y_{\eta\hat{\ }\langle\beta\rangle}=\bz).\]  
\end{description}
Then the sequence $\langle x^{\red}_\alpha:\alpha<\lambda\rangle$ examplifies
that $\ba^{\red}(\C)$ fails the $\lambda$-Knaster condition.
\end{proposition}

\noindent{\bf Explanation:}\ \ \ The above proposition is not just something
in the direction of Problem~\ref{pwa}. The tuple $(\ba^{\red},\bar{x},
J^{\bd}_\lambda)$ is like $(\ba_\eta,\bar{y}_\eta, I_\eta)$ but
$J^{\bd}_\lambda$ is nicer than ideals given by previous results. Using such
objects makes building examples for Problem~\ref{pwa} much easier.  
\medskip

\Proof It is enough to show that
\begin{quotation}
\noindent for each $Y\in [\lambda]^{\textstyle\lambda}$ one can find
$\varepsilon,\zeta\in Y$ such that
\[\ba_{\eta_\varepsilon\rest i}\models y_{\eta_\varepsilon\rest (i+1)}\cap
y_{\eta_\zeta\rest (i+1)}=\bz\]
where $i=\lh(\eta_\varepsilon\wedge\eta_\zeta)$.
\end{quotation}
For this, for each $\nu\in T$ we put
\[A_\nu\stackrel{\rm def}{=}\{\alpha<\lambda_\nu: (\exists\varepsilon\in Y)
(\nu\hat{\ }\langle\alpha\rangle\vartriangleleft\eta_\varepsilon)\}.\]

\begin{claim}
There is $\nu\in T$ such that $A_\nu\notin I_\nu$.
\end{claim}

\noindent{\em Why?}\ \ \ First note that by the definition of $A_\nu$, for
each $\varepsilon\in Y$ we have
\[(\forall i<\delta)(\eta_\varepsilon \hat{\ }\langle i\rangle \in
A_{\eta_\varepsilon\rest i}).\] 
Now, if we had that $A_\nu\in I_\nu$ for all $\nu\in T$ then we could apply
the assumption that $\bar{\eta}$ is $J$-cofinal for
$(T,\bar{\lambda},\bar{I})$ to the sequence $\langle A_\nu: \nu\in T\rangle$.
Thus we would find $\alpha^*<\lambda$ such that 
\[\alpha^*\leq\alpha<\lambda\ \ \ \Rightarrow\ \ \ \{i<\delta: \eta_\alpha(i)
\notin A_{\eta_\alpha \rest i}\}\in J\]
which contradicts our previous remark (remember $|Y|=\lambda$). The claim is
proved. 
\medskip

Due to the claim we find $\nu\in T$ such that $A_\nu\notin I_\nu$. By the part
{\bf (b)} of our assumptions we find $\alpha,\beta\in A_\nu$ such that 
\[\ba_\nu\models y_{\nu\hat{\ }\langle\alpha\rangle}\cap y_{\nu\hat{\
}\langle\beta\rangle}=\bz.\]
Choose $\varepsilon,\zeta\in Y$ such that $\nu\hat{\ }\langle\alpha\rangle
\vartriangleleft\eta_\varepsilon$, $\nu\hat{\ }\langle\beta\rangle
\vartriangleleft\eta_\zeta$ (see the definition of $A_\nu$). Then
$\nu=\eta_\varepsilon\wedge\eta_\zeta$ and 
\[\ba_\nu\models y_{\eta_\varepsilon\rest (i+1)}\cap y_{\eta_\zeta\rest
(i+1)}=\bz\]
(where $i=\lh(\nu)$), finishing the proof of the proposition. \QED

\begin{lemma}
\label{consknaster}
Let $\C=(T, \bar{\lambda}, \bar{\eta},\langle(\ba_\eta,\bar{y}_\eta):
\eta\in T\rangle)$ be a $(\delta,\mu,\lambda)$-constructor such that 
\begin{description}
\item[$(\bigstar)$] the Boolean algebras $\ba_\eta$ satisfy  the
$(2^{|\delta|})^+$--Knaster condition.
\end{description}
{\em Then} the Boolean algebra $\ba^{\red}(\C)$ satisfies the
$(2^{|\delta|})^+$--Knaster condition.  
\noindent In fact we may replace $(2^{|\delta|})^+$ above by any regular
cardinal $\theta$ such that
\[(\forall\alpha<\theta)(|\alpha|^{|\delta|}<\theta).\]  

\noindent To get that $\ba^{\red}(\C)$ satisfies the $(2^{|\delta|})^+$--cc it
is enough if instead of $(\bigstar)$ we assume
\begin{description}
\item[$(\bigstar\bigstar)$] every free product of finitely many of the Boolean
algebras $\ba_\eta$ satisfies the $(2^{|\delta|})^+$--cc. 
\end{description}
\end{lemma}

\noindent{\bf Remark}: \ \ \ 1. Usually we will have $\delta=\cf(\mu)$.

\noindent 2. Later we will get more (e.g. $|\delta|^+$-Knaster if
$(T,\bar{\eta})$ is hereditarily free, see \ref{free}, \ref{morekna}). 
\medskip

\Proof Let $\theta=(2^{|\delta|})^+$ and assume $(\bigstar)$ (the other cases
have the same proofs). Suppose that
$z_\varepsilon\in\ba^{\red}\setminus\{\bz\}$ (for $\varepsilon<\theta$). We 
start with a series of reductions which we describe fully here but later, in
similar situations, we will state what is the result of the procedure only.
\medskip

\noindent{\bf Standard cleaning}:\ \ \ Each $z_\varepsilon$ is a Boolean
combination of some generators $x_{i_0},\ldots,x_{i_{n-1}}$. But, as we want
to find a subsequence with non-zero intersections, we may replace
$z_\varepsilon$ by any non-zero $z\leq z_\varepsilon$. Consequently we may
assume that each $z_\varepsilon$ is an intersection of some generators or
their complements. Further, as $\cf(\theta)=\theta>\aleph_0$ we may assume that
the number of generators needed for this representation does not depend on
$\varepsilon$ and is equal to, say, $n^*$. Thus we have two functions
\[i:\theta\times n^*\longrightarrow\lambda\quad\mbox{ and }\quad
\gt:\theta\times n^*\longrightarrow 2\] such that for each
$\varepsilon<\theta$: \[z_\varepsilon=\bigcap_{\ell<n^*}
(x_{i(\varepsilon,\ell)})^{\gt(\varepsilon,\ell)}\] and there is no repetition
in $\langle i(\varepsilon,\ell):\ell<n^*\rangle$.  Moreover we may assume that
$\gt(\varepsilon,\ell)$ does not depend on $\varepsilon$, i.e.
$\gt(\varepsilon,\ell)=\gt(\ell)$.  Due to the $\Delta$-lemma for finite sets
we may assume that $\langle \langle i(\varepsilon,\ell):\ell<n^*\rangle:
\varepsilon<\theta\rangle$ is a $\Delta$-system of sequences, i.e.:
\begin{description} \item[$(*)_1$]
$i(\varepsilon,\ell_1)=i(\varepsilon,\ell_2)\ \ \Rightarrow\ \ \ell_1=\ell_2$\
\ \ \ and \item[$(*)_2$] for some $w\subseteq n^*$ we have
\[\hspace{-0.7cm}(\exists\varepsilon_1\!<\!\varepsilon_2<\theta)
(i(\varepsilon_1,\ell)\!=\! i(\varepsilon_2,\ell))\ \mbox{ iff }\
(\forall\varepsilon_1,\varepsilon_2\!<\!\theta)(i(\varepsilon_1,\ell)\!=\!
i(\varepsilon_2,\ell))\ \mbox{ iff }\ \ell\in w\] \end{description} Now note
that, by the definition of the algebra $\ba^{\red}$, \begin{description}
\item[$(*)_3$] $z_{\varepsilon_1}\cap z_{\varepsilon_2}=\bz$\ \ \ \ if and
only if
\[\hspace{-1cm}\bigcap\{x^{\gt(\ell)}_{i(\varepsilon_1,\ell)}:
\ell<n^*, \gt(\ell)=0\}
\cap\bigcap\{x^{\gt(\ell)}_{i(\varepsilon_2,\ell)}: \ell<n^*,
\gt(\ell)=0\}=\bz.\] \end{description} Consequently we may
assume that
\[(\forall\ell<n^*)(\forall\varepsilon<\theta)(\gt(\ell)=0).\]
\medskip

\noindent{\bf Explanation of what we are going to do now:}\ \ \ We want to
replace the sequence $\langle z_\varepsilon:\varepsilon<\theta\rangle$ by a
large subsequence such that the places of splitting between two branches used
in two different $z_\varepsilon$'s will be uniform. Then we will be able to
translate our $\theta$--cc problem to the one on the algebras $\ba_\eta$.
\medskip

Let 
\[A_\varepsilon\stackrel{\rm def}{=}\{\nu\in {}^{\delta{>}}\mu: (\exists
j<\varepsilon)(\exists \ell<n^*)(\nu\vartriangleleft\eta_{i(j,\ell)})\}\]
and let $B_\varepsilon$ be the closure of $A_\varepsilon$:
\[\begin{array}{ll}
B_\varepsilon\stackrel{\rm def}{=}\{\rho\in {}^{\delta{\geq}}\mu: &\rho\in
A_\varepsilon \mbox{ or } \lh(\rho)\mbox{ is a limit ordinal and}\\
\ &(\forall\zeta<\lh(\rho))(\rho\rest\zeta\in A_\varepsilon)\}\\
  \end{array}\]
Note that $|A_\varepsilon|\leq |\varepsilon|\cdot |\delta|$ and hence
$|B_\varepsilon|\leq |A_\varepsilon|^{{\leq}|\delta|}<\theta$. Next we define
(for $\varepsilon<\theta$, $\ell<n^*$): 
\[\zeta(\varepsilon,\ell)\stackrel{\rm def}{=}\sup\{\zeta<\delta:
\eta_{i(\varepsilon,\ell)} \rest\zeta\in B_\varepsilon\}.\]
Thus $\zeta(\varepsilon,\ell)\leq\lh(\eta_{i(\varepsilon,\ell)})=\delta$. Let
$\S=\{\varepsilon<\theta:\cf(\varepsilon)>|\delta|\}$. For each
$\varepsilon\in\S$ we neccessarily have
\[\eta_{i(\varepsilon,\ell)}\rest\zeta(\varepsilon,\ell)\in
B_\varepsilon\quad\mbox{ and }\quad B_\varepsilon=\bigcup_{\xi<\varepsilon}
B_\xi\] 
(remember that $\cf(\varepsilon)>|\delta|$ and for limit $\varepsilon$ we have
$A_\varepsilon=\bigcup\limits_{\xi<\varepsilon} A_\xi$) and hence
\[\eta_{i(\varepsilon,\ell)}\rest\zeta(\varepsilon,\ell)\in
B_{\xi(\varepsilon,\ell)},\quad\quad\mbox{ for some }\xi(\varepsilon,\ell)<
\varepsilon.\] 
Let $\xi(\varepsilon)=\max\{\xi(\varepsilon,\ell):\ell<n^*\}$. By the Fodor
lemma we find $\xi^*<\theta$ such that the set
\[\S_1\stackrel{\rm def}{=}\{\varepsilon\in\S: \xi(\varepsilon)=\xi^*\}\]
is stationary. Thus $\eta_{i(\varepsilon,\ell)}\rest
\zeta(\varepsilon,\ell)\in B_{\xi^*}$ for each $\varepsilon\in\S_1$,
$\ell<n^*$. Since $|B_{\xi^*}|,|\delta|<\theta$ we find
$\nu_0,\ldots,\nu_{n^*-1}\in B_{\xi^*}$ and $\alpha(\ell_1,\ell_2)\leq\delta$
(for $\ell_1\leq\ell_2<n^*$) such that the set
\[\begin{array}{ll}
\S_2\stackrel{\rm def}{=}\{\varepsilon\in\S_1:& (\forall
\ell< n^*)(\eta_{i(\varepsilon,\ell)}{\rest}\zeta(\varepsilon,\ell)
=\nu_\ell)\ \&\\
\ &\&\ (\forall\ell_1\leq\ell_2<n^*)(\lh(\eta_{i(\varepsilon,\ell_1)}\wedge
\eta_{i(\varepsilon,\ell_2)})=\alpha(\ell_1,\ell_2))\}\\
\end{array}\]  
is stationary. Further, applying the $\Delta$-lemma we find a set $\S_3\in
[\S_2]^{\textstyle \theta}$ such that  
\[\{\langle\eta_{i(\varepsilon,\ell)}(\lh(\nu_\ell)): \ell<n^*\rangle:
\varepsilon\in\S_3\}\] 
forms a $\Delta$-system of sequences.

For $\varepsilon\in\S_3$ and $\nu\in T$ denote
\[b^\varepsilon_\nu\stackrel{\rm def}{=}\bigcap\{y_{\eta_{i(\varepsilon,
\ell)}\rest (\lh(\nu)+1)}: \ell<n^*, \nu\vartriangleleft\eta_{i
(\varepsilon,\ell)}\}\in\ba_\nu.\]

\begin{claim}
For each $\varepsilon\in\S_3$, $\nu\in T$ the element
$b^\varepsilon_\nu$ (of the algebra $\ba_\nu$) is non-zero.
\end{claim}

\noindent{\em Why?}\ \ \ Because of the definition of $\ba^{\red}$ and the
fact that $z_\varepsilon\neq\bz$: 
\[b^\varepsilon_\nu=\bz\ \ \ \Rightarrow\ \ \
\bigcap\{x_{\eta_{i(\varepsilon,\ell)}}: \ell<n^*, \nu\vartriangleleft\eta_{i  
(\varepsilon,\ell)}\}=\bz \ \ \ \Rightarrow\ \ \ z_\varepsilon=\bz.\]
\medskip

Since for each $\ell<n^*$ the algebra $\ba_{\nu_\ell}$ satisfies the
$\theta$--Knaster property we find a set $\S_4\in [\S_3]^{\textstyle \theta}$
such that for each $\ell<n^*$ and $\varepsilon_1,\varepsilon_2\in\S_4$ we have
\[\varepsilon_1\neq\varepsilon_2\ \ \ \Rightarrow\ \ \
b^{\varepsilon_1}_{\nu_\ell}\cap b^{\varepsilon_2}_{\nu_\ell}\neq\bz
\quad\mbox{ in }\ba_{\nu_\ell}.\]

Now we may finish by proving the following claim.

\begin{claim}
For each $\varepsilon_1,\varepsilon_2\in\S_4$
\[\ba^{\red}\models z_{\varepsilon_1}\cap z_{\varepsilon_2}\neq\bz.\]
\end{claim}

\noindent{\em Why?}\ \ \ Since $z_{\varepsilon_1}\cap z_{\varepsilon_2}$ is
just the intersection of generators it is enough to show that (remember the
definition of $\ba^{\red}$):
\begin{description}
\item[$(\otimes)$] for each $\varepsilon_1,\varepsilon_2\in\S_4$ and for every
$\nu\in T$ 
\[\ba_\nu\models\bigcap\{y_{\eta_i\rest(\lh(\nu)+1)}:
i\in\{i(\varepsilon_1,\ell), i(\varepsilon_2,\ell): \ell<n^*\} 
\mbox{ and } \nu\vartriangleleft\eta_i\}\neq\bz.\]
\end{description}
If $\nu=\nu_\ell$, $\ell<n^*$ then the intersection is
$b^{\varepsilon_1}_{\nu_\ell}\cap b^{\varepsilon_2}_{\nu_\ell}$ which by the
choice of the set $\S_4$ is not zero. So suppose that $\nu\notin\{\nu_\ell:
\ell<n^*\}$. Put 
\[u_\nu\stackrel{\rm def}{=}\{i: \nu\vartriangleleft\eta_i\mbox{ and for some
}\ell<n^* \mbox{ either } i=i(\varepsilon_1,\ell) \mbox{ or }
i=i(\varepsilon_2,\ell)\}.\]   
If 
\[\{\eta_i(\lh(\nu)): i\in u_\nu\}\subseteq\{\eta_{i(\varepsilon_2,\ell)}
(\lh(\nu)):\ell<n^*\ \&\ \nu\vartriangleleft\eta_{i(\varepsilon_2,\ell)}\}\]
then we are done as $b^{\varepsilon_2}_\nu\neq\bz$. So there is $\ell_1<n^*$
such that $\nu\vartriangleleft\eta_{i(\varepsilon_1,\ell_1)}$ and
\[\eta_{i(\varepsilon_1,\ell_1)}\rest
(\lh(\nu)+1)\notin\{\eta_{i(\varepsilon_2,\ell)}\rest (\lh(\nu)+1):
\ell<n^*\ \&\ \nu\vartriangleleft\eta_{i(\varepsilon_2,\ell)}\}.\]
Similarly we may assume that there is $\ell_2<n^*$ such that
$\nu\vartriangleleft\eta_{i(\varepsilon_2,\ell_2)}$ and
\[\eta_{i(\varepsilon_2,\ell_2)}\rest
(\lh(\nu)+1)\notin\{\eta_{i(\varepsilon_1,\ell)}\rest (\lh(\nu)+1):
\ell<n^*\ \&\ \nu\vartriangleleft\eta_{i(\varepsilon_1,\ell)}\}.\]
Because of the symmetry we may assume that $\varepsilon_1<\varepsilon_2$. Then
\[\nu=\eta_{i(\varepsilon_2,\ell_2)}\rest\lh(\nu)\in
A_{\varepsilon_1+1}\subseteq B_{\varepsilon_2}\]
and hence $\zeta(\varepsilon_2,\ell_2)\geq\lh(\nu)$. By the choice of $\S_2$
(remember $\varepsilon_1,\varepsilon_2\in\S_4\subseteq\S_2$), we get 
$\nu\trianglelefteq\nu_{\ell_2}$. But we have assumed that
$\nu\neq\nu_{\ell_2}$, so $\nu\vartriangleleft\nu_{\ell_2}$. Hence (once again
due to $\varepsilon_1,\varepsilon_2\in\S_2$) 
\[\eta_{i(\varepsilon_2,\ell_2)}\rest (\lh(\nu)+1)=
\eta_{i(\varepsilon_1,\ell_2)}\rest (\lh(\nu)+1)=\nu_{\ell_2}\rest
(\lh(\nu)+1)\]  
which contradicts the choice of $\ell_2$. 

The claim and so the lemma are proved. \QED 

\begin{remark}
{\em
We can strengthen ``$\theta$-Knaster'' in the assumption and conclusion of
\ref{consknaster} in various ways. For example we may have that ``intersection
of any $n$ members of the final set is non-zero''.  
}
\end{remark}

\begin{definition}
Let $(\ba,\bar{y})$ be a $\lambda$-marked Boolean algebra,
$\kappa\leq\lambda$. We say that 
\begin{enumerate}
\item $(\ba,\bar{y})$ satisfies the $\kappa$-Knaster property if $\ba$
satisfies the definition of the $\kappa$-Knaster property (see \ref{knaster})
with restriction to subsequences of $\bar{y}$.
\item $(\ba,\bar{y})$ is $(\kappa,{\rm not}\lambda)$--Knaster if
\begin{description}
\item[(a)] the algebra $\ba$ has the $\kappa$-Knaster property but
\item[(b)] the sequence $\bar{y}$ witnesses that the $\lambda$-Knaster
property fails for $\ba$. 
\end{description}
\end{enumerate}
\end{definition}

\begin{conclusion}
\label{knacon}
Assume that $\mu$ is a strong limit singular cardinal, $\lambda=2^\mu=\mu^+$
and $\theta=(2^{\cf(\mu)})^+$. 

\noindent{\em Then} there exists a $\lambda$-marked Boolean algebra
$(\ba,\bar{y})$ which is $(\theta,{\rm not}\lambda)$--Knaster.
\end{conclusion}

\Proof Choose cardinals $\mu_i^0,\mu_i<\mu$ (for $i<\cf(\mu)$) such that 
\begin{description}
\item[$(\alpha)$] $\cf(\mu)<\mu^0_0$
\item[$(\beta)$]  $\prod\limits_{j<i}\mu_j<\mu^0_i$, $\mu_i=(2^{\mu^0_i})^+$
\item[$(\gamma)$] the sequences $\langle\mu_i: i<\cf(\mu)\rangle$,
$\langle\mu_i^0: i<\cf(\mu)\rangle$ are increasing cofinal in $\mu$
\end{description}
(possible as $\mu$ is strong limit singular). By proposition \ref{algide} we
find $\mu_i$-marked Boolean algebras $(\ba_i,\bar{y}^i)$ and
$(\mu^0_i)^+$-complete ideals $I_i$ on $\mu_i$ (for $i<\delta$) such that 
\begin{description}
\item[(a)] if $X\subseteq\mu_i$, $X\notin I_i$ then $(\exists\alpha,\beta \in
X)(\ba_i\models y^i_\alpha\cap y^i_\beta=\bz)$
\item[(b)] the algebra $\ba_i$ has the $(2^{\cf(\mu)})^+$--Knaster property.
\end{description}
Let $T=\bigcup\limits_{i<\cf(\mu)}\prod\limits_{j<i}\mu_j$ and for $\nu\in
T_i$ ($i<\cf(\mu)$) let $I_\nu=I_i$, $\ba_\nu=\ba_i$, $\bar{y}_\nu=\bar{y}^i$
and $\lambda_\nu=\mu_i$. Now we may apply proposition \ref{stronger} to
$\mu$, $\langle\mu^0_i: i<\cf(\mu)\rangle$ and $(T,\bar{\lambda},\bar{I})$ to
find a stronger $J^{\bd}_{\cf(\mu)}$-cofinal sequence $\bar{\eta}$ for
$(T,\bar{\lambda},\bar{I})$ of the length $\lambda$. Consider the
$(\cf(\mu),\mu,\lambda)$-constructor $\C=(T,\bar{\lambda},\bar{\eta},\langle
(\ba_\nu,\bar{y}_\nu): \nu\in T\rangle)$. By {\bf (b)} above we may apply
lemma \ref{consknaster} to get that the algebra $\ba^{\red}(\C)$
satisfies the $(2^{\cf(\mu)})^+$--Knaster condition. Finally we use
proposition \ref{notknaster} (and {\bf (a)} above) to conclude that
$(\ba^{\red}(\C),\langle x^{\red}_\alpha:\alpha<\lambda\rangle)$ is
$(\theta,{\rm not}\lambda)$--Knaster. \QED

\begin{proposition}
\label{clospcf}
Assume that:\\
$\kappa$ is a regular cardinal such that
$(\forall\alpha<\kappa)(|\alpha|^{|\delta|}<\kappa)$, $\bar{\lambda}=
\langle\lambda_i: i<\delta\rangle$ is an increasing sequence of regular
cardinals such that $\kappa\leq\lambda_0$,
$\prod\limits_{j<i}\lambda_j<\lambda_i$ (or just $\max\pcf\{\lambda_j: j<i\}<
\lambda_i$) for $i<\delta$ and $\lambda\in\pcf\{\lambda_i: i<\delta\}$.
Further suppose that for each $i<\delta$ there exists a $\lambda_i$-marked
Boolean algebra which is $(\kappa,{\rm not}\lambda_i)$--Knaster. 

\noindent Then there exists a $\lambda$-marked Boolean algebra which is
$(\kappa,{\rm not}\lambda)$--Knaster. 
\end{proposition}

\Proof If $\lambda=\lambda_i$ for some $i<\delta$ then there is nothing to do.
If $\lambda<\lambda_i$ for some $i<\delta$ then let $\alpha<\delta$ be the
maximal limit ordinal such that $(\forall i<\alpha)(\lambda_i<\lambda)$ (it
neccessarily exists) . Now we may replace $\langle_i:i<\delta\rangle$ by
$\langle \lambda_i: i<\alpha\rangle$. Thus we may assume that $(\forall
i<\delta)(\lambda_i<\lambda)$. Further we  may assume that 
\[\lambda=\max\pcf\{\lambda_i: i<\delta\}\]
(by \cite{Sh:g}, I, 1.8). Now, due to \cite{Sh:g}, II, 3.5, p.65, we find a
sequence $\bar{\eta}\subseteq\prod\limits_{i<\delta}\lambda_i$ and an ideal
$J$ on $\delta$ such that    
\begin{enumerate}
\item $J\supseteq J^{\bd}_\delta$ and $\lambda=\tcf(\prod\limits_{i<\delta}
\lambda_i/J)$ 

\noindent (naturally: $J=\{a\subseteq\delta: \max\pcf\{\lambda_i: i\in a\}<
\lambda\}$) 
\item $\bar{\eta}=\langle\eta_\varepsilon:\varepsilon<\lambda\rangle$ is
${<}_J$-increasing cofinal in $\prod\limits_{i<\delta}\lambda_i/J$
\item for each $i<\delta$
\[|\{\eta_\varepsilon\rest i: \varepsilon<\lambda\}|<\lambda_i.\]
\end{enumerate}
Let $T=\bigcup\limits_{i<\delta}\prod\limits_{j<i}\lambda_j$ and for
$\nu\in T_i$ ($i<\delta$) let $\lambda_\nu=\lambda_i$,
$I_\nu=J^{\bd}_{\lambda_i}$.  

\noindent It follows from the choice of $\bar{\eta}, J$ above and our
assumptions that we may apply proposition \ref{getcofinal} and hence
$\bar{\eta}$ is $J$-cofinal for $(T,\bar{\lambda},\bar{I})$. For $\nu\in T$
let $(\ba_\nu,\bar{y}_\nu)$ be a $\lambda_\nu$-marked $(\kappa,{\rm
not}\lambda_\nu)$--Knaster Boolean algebra (exists by our assumptions). Now we
may finish using \ref{consknaster} and \ref{notknaster} for
$\C=(T,\bar{\lambda},\bar{\eta},\langle(\ba_\eta,\bar{y}_\eta):\eta\in 
T\rangle)$, $\bar{I}$ and $J$ (note the assumption {\bf (b)} of
\ref{notknaster} is satisfied as $I_\eta=J^{\bd}_{\lambda_\eta}$; remember the
choice of $(\ba_\eta,\bar{y}_\eta)$). \QED

\begin{remark}
{\em
Note that from cardinal arithmetic hypothesis $\cf(\mu)=\chi$, $\chi^{{<}\chi}
<\chi<\mu$, $\mu^+=\lambda<2^\chi$ alone we cannot hope to build a
counterexample. This is because of section 4 of \cite{Sh:93}, particularly
lemma 4.13 there. It was shown in that paper that if $\chi^{{<}\chi}<\chi_1=
\chi_1^{{<}\chi_1}$ then there is a $\chi^+$-cc $\chi$-complete forcing notion
$\p$ of size $\chi_1$ such that
\[\begin{array}{ll}
\forces_{\p}&\mbox{``if }|\ba|<\chi_1, \ba\models \chi\mbox{-cc}\\
\ &\mbox{then }\ba^+\mbox{ is the union of }\chi\mbox{ ultrafilters''}.\\
\end{array}\]
More on this see in section 8.

\noindent So the centrality of $\lambda\in\Reg\cap (\mu,2^\mu]$, $\mu$ strong
limit singular, is very natural.
}
\end{remark}

\section{The main result}

\begin{proposition}
\label{notcc}
Suppose that $\C$ is a $(\delta,\mu,\lambda)$-constructor. Then the free
product $\ba^{\red}(\C)*\ba^{\green}(\C)$ fails the $\lambda$-cc (so
$\cel(\ba^{\red}(\C)*\ba^{\green}(\C))\geq\lambda$). 
\end{proposition}

\Proof Look at the elements $(x_i,x_i)\in\ba^{\red}*\ba^{\green}$ for
$i<\lambda$. It follows directly from the definition of the algebras that for
each $i<j<\lambda$:
\[\mbox{either }\quad \ba^{\red}\models x_i^{\red}\cap x_j^{\red}=\bz
\quad\mbox{ or }\quad \ba^{\green}\models x_i^{\green}\cap x_j^{\green}=\bz.\]
Consequently the sequence $\langle (x_i,x_i):i<\lambda\rangle$ witnesses the
assertion of the proposition. \QED

\begin{proposition}
\label{finite}
Suppose that $n<\omega$ and for $\ell\leq n$:
\begin{enumerate}
\item $\chi_\ell,\lambda_\ell$ are regular cardinals, $\chi_\ell<\lambda_\ell
<\chi_{\ell+1}$ 
\item $(\ba_\ell,\bar{y}_\ell, I_\ell)$ is a $(\lambda_\ell,\chi_\ell)$-well
marked Boolean algebra (see definition \ref{well1}), $\bar{y}_\ell= \langle
y^\ell_i: i<\lambda_\ell\rangle$ 
\item $\ba$ is the Boolean algebra freely generated by $\{y_\eta:
\eta\in\prod\limits_{\ell\leq n}\lambda_\ell\}$ {\em except that} 
\begin{quotation}
\noindent {\em if} $\eta_{i_0},\ldots,\eta_{i_{k-1}}\in\prod\limits_{\ell\leq
n}\lambda_\ell$, $\eta_{i_0}\rest \ell=\eta_{i_1}\rest\ell=\ldots=
\eta_{i_{n-1}}\rest\ell$ and\\
$\ba_\ell\models\bigcap\limits_{m<k} y^\ell_{\eta_{i_m}(\ell)}=\bz$ 

\noindent {\em then} $\bigcap\limits_{m<k} y_{\eta_{i_m}}=\bz$.
\end{quotation}
[Compare to the definition of the algebras $\ba^{\red}(\C)$.]
\item $I=\{B\subseteq\prod\limits_{\ell\leq n}\lambda_\ell: \neg
(\exists^{I_0}\gamma_0)\ldots(\exists^{I_n}\gamma_n)(\langle\gamma_0,\ldots,
\gamma_n\rangle\in B)\}$.
\end{enumerate}
{\em Then:}
\begin{description}
\item[(a)] if all the algebras $\ba_\ell$ (for $\ell\leq n$) satisfy the
$\theta$-Knaster property, $\theta$ is a regular uncountable cardinal

then $\ba$ has the $\theta$-Knaster property;
\item[(b)] $I$ is a $\chi_0$-complete ideal on $\prod\limits_{\ell\leq
n}\lambda_i$;
\item[(c)] if $Y\subseteq (\prod\limits_{\ell\leq n}\lambda_\ell)^n$ is such
that 
\[(\exists^{I}\eta_0)\ldots(\exists^{I}\eta_n)(\langle\eta_0,\ldots,
\eta_n\rangle\in Y)\]
then there are $\langle\eta_0^\prime,\ldots,\eta_n^\prime\rangle,
\langle\eta_0^{\prime\prime},\ldots,\eta_n^{\prime\prime}\rangle \in Y$ such
that for all $\ell\leq n$
\[\ba\models y_{\eta^\prime_\ell}\cap y_{\eta^{\prime\prime}_\ell}=\bz.\]
\end{description}
\end{proposition}

\Proof {\bf (a)}\ \ \ The proof that the algebra $\ba$ satisfies
$\theta$-Knaster condition is exactly the same as that of \ref{consknaster}
(actually it is a special case of that).

\noindent{\bf (b)}\ \ \ Should be clear.

\noindent{\bf (c)}\ \ \ For $\ell,m\leq n$ put
\[\chi^m_\ell=\chi_\ell, \lambda^m_\ell=\lambda_\ell,\ I^m_\ell= I_\ell,\
P^m_\ell= \{\{\alpha,\beta\}\subseteq\lambda_\ell: \ba_\ell\models
y^\ell_\alpha\cap y^\ell_\beta=\bz\}, B=Y.\] 
It is easy to check that the assumptions of proposition \ref{fubini} are
satisfied. Applying it we find sets $X_0,\ldots,X_n$ satisfying the respective
versions of clauses {\bf (a)}--{\bf (d)} there. Note that our choice of the
sets $P^m_\ell$ and clauses {\bf (b)}, {\bf (c)} of \ref{fubini} imply that
\begin{quotation}
\noindent $X_m=\{\nu_m^\prime,\nu_m^{\prime\prime}\}\subseteq\prod
\limits_{\ell\leq n-m}\lambda_\ell$

\noindent $\nu_m^\prime\rest (n-m)=\nu_m^{\prime\prime}\rest (n-m)$

\noindent $\ba_{n-m}\models y^{n-m}_{\nu_m^{\prime}(n-m)}\cap
y^{n-m}_{\nu_m^{\prime\prime}(n-m)}=\bz$
\end{quotation}
Look at the sequences $\langle\nu^\prime_0,\ldots,\nu^\prime_n\rangle$,
$\langle\nu^{\prime\prime}_0,\ldots,\nu^{\prime\prime}_n\rangle$. By the
clause {\bf (d)} of \ref{fubini} we find $\langle\eta^\prime_0,\ldots,
\eta^\prime_m\rangle,\langle\nu^{\prime\prime}_0,\ldots,
\nu^{\prime\prime}_n\rangle\in Y$ such that for each $m\leq n$
\[\nu_m^\prime\trianglelefteq\eta^\prime_m,\quad \nu_m^{\prime\prime}
\trianglelefteq\eta^{\prime\prime}_m.\]
Now, the properties of $\nu_m^\prime$, $\nu^{\prime\prime}_m$ and the
definition of the algebra $\ba$ imply that for each $m\leq n$:
\[\ba\models y_{\eta^\prime_m}\cap y_{\eta_m^{\prime\prime}}=\bz,\]
finishing the proof. \QED 

\begin{lemma}
\label{green}
Assume that $\lambda$ is a regular cardinal, $|\delta|<\lambda$, $J$ is an
ideal on $\delta$ extending $J^{\bd}_\delta$,
$\C=(T,\bar{\lambda},\bar{\eta},\langle(\ba_\eta,\bar{y}_\eta): \eta\in
T\rangle)$ is a $(\delta,\mu,\lambda)$-constructor and  $\bar{I}$ is such that
$(T,\bar{\lambda},\bar{I})\in \K^{\id}_{\delta,\mu}$. Suppose that
$\bar{\eta}=\langle \eta_\alpha: \alpha<\lambda\rangle$ is a stronger
(or big) $J$-cofinal in $(T,\bar{\lambda},\bar{I})$ sequence such that
\[(\forall i<\delta)(|\{\eta_\alpha\rest i: \alpha<\lambda\}|<\lambda).\]
Further assume that
\begin{description}
\item[$(\circleddash)$] for every $n<\omega$ for a $J$-positive set of
$i<\delta$ we have:
\begin{quotation}
\noindent {\em if} $\eta_0,\ldots,\eta_n\in T_i$ are pairwise distinct and the
set $Y\subseteq\prod\limits_{\ell\leq n}\lambda_{\eta_\ell}$ is such that
\[(\exists^{I_{\eta_0}}\gamma_0)\ldots(\exists^{I_{\eta_n}}\gamma_n)(\langle
\gamma_0,\ldots,\gamma_n\rangle\in Y)\]

\noindent{\em then} for some $\gamma_\ell^\prime,\gamma_\ell^{\prime\prime}<
\lambda_{\eta_\ell}$ (for $\ell\leq n$) we have
\[\langle\gamma_\ell^\prime:\ell\leq n\rangle, \langle\gamma_\ell^{\prime
\prime}:\ell\leq n\rangle\in Y\quad\mbox{ and for all }\ell\leq n\]
\[\ba_{\eta_{\ell}}\models y_{\eta_\ell\hat{\ }\langle\gamma^\prime_\ell
\rangle}\cap y_{\eta_\ell\hat{\ }\langle\gamma^{\prime\prime}_\ell\rangle}=
\bz.\] 
\end{quotation}
\end{description}
{\em Then} the Boolean algebra $\ba^{\green}(\C)$ satisfies $\lambda$-cc.
\end{lemma}

\Proof Suppose that $\langle z_\alpha: \alpha<\lambda\rangle\subseteq
\ba^{\green}\setminus\{\bz\}$. By the standard cleaning (compare the first part
of the proof of \ref{consknaster}) we may assume that there are $n^*\in\omega$
and a function $\varepsilon:\lambda\times n^*\longrightarrow\lambda$ such that 
\begin{enumerate}
\item $z_\alpha=\bigcap\limits_{\ell<n^*}x_{\varepsilon(\alpha,\ell)}$ (in
$\ba^{\green}$)
\item $\varepsilon(\alpha,0)<\varepsilon(\alpha,1)<\ldots<
\varepsilon(\alpha,n^*-1)$ 
\item $\langle\langle\varepsilon(\alpha,\ell):\ell<n^*\rangle:
\alpha<\lambda\rangle$ forms a $\Delta$-system of sequences with the kernel
$m^*$, i.e. $(\forall \ell<m^*)(\varepsilon(\alpha,\ell)=\varepsilon(\ell))$
and 
\[(\forall \ell\in [m^*,n^*))(\forall\alpha<\lambda)(\varepsilon(\alpha,\ell)
\notin\{\varepsilon(\beta,k):(\beta,k)\neq (\alpha,\ell)\})\]
\item there is $i^*<\delta$ such that for each $\alpha<\lambda$ there is no
repetition in the sequence $\langle\eta_{\varepsilon(\alpha,\ell)}\rest i^*:
\ell<n^*\rangle$.
\end{enumerate}
Since $|\{\eta_\alpha\rest i: \alpha<\lambda\}|<\lambda$ (for $i<\delta$) and
$|\delta|<\lambda$ we may additionally require that 
\begin{description}
\item[$\hat{(*)}$] for each $i<\delta$, for every $\alpha<\lambda$ we have
\[(\exists^\lambda \beta<\lambda)(\forall \ell<n^*)(\eta_{\varepsilon(\alpha,
\ell)}\rest (i+1)=\eta_{\varepsilon(\beta,\ell)}\rest (i+1))\]
and
\item[$\hat{(**)}$] for each $\alpha<\beta<\lambda$, $\ell<n^*$
\[\eta_{\varepsilon(\alpha,\ell)}\rest i^*=\eta_{\varepsilon(\beta,\ell)}
\rest i^*.\]
\end{description}
\medskip

\noindent{\bf Remark}: \ \ \ \ Note that the claim below is like an
$(n^*-m^*)$-place version of \ref{notknaster}. Having an $(n^*-m^*)$-ary
version is extra for the construction but it also costs.

\begin{claim}
\label{claimgreen}
{\em Assume that:}

\noindent $\C=(T,\bar{\lambda},\bar{\eta},\langle(\ba_\eta,\bar{y}_\eta): 
\eta\in T\rangle$ is a $(\delta,\mu,\lambda)$-constructor, $\lambda$ a regular
cardinal, $\delta<\lambda$, $\bar{I}$ is such that
$(T,\bar{\lambda},\bar{I})\in \K^{\id}_{\delta,\mu}$, $J$ is an ideal on
$\delta$ extending $J^{\bd}_\delta$ and the sequence $\bar{\eta}$ is stronger 
$J$-cofinal in $(T,\bar{\lambda},\bar{I})$.\\
Further suppose that $\varepsilon:\lambda\times n^*\longrightarrow$, $m^*,n^*$
and $i^*<\delta$ are as above (after the reduction, but the property
$\hat{(**)}$ is not needed). 

\noindent {\em Then}
\begin{description}
\item[$(\boxtimes)$] for every large enough $\alpha<\lambda$ the set:
\[\hspace{-1cm}\begin{array}{l}
Z_\alpha\stackrel{\rm def}{=}\{i<\delta:
\neg(\exists^{I_{\eta_{\varepsilon(\alpha,m^*)}\rest i}}  
\gamma_{m^*})(\exists^{I_{\eta_{\varepsilon(\alpha,m^*+1)}\rest i}}
\gamma_{m^*+1})\ldots\\
\ \ldots(\exists^{I_{\eta_{\varepsilon(\alpha,n^*-1)}\rest i}}
\gamma_{n^*-1})(\exists^\lambda \beta)(\forall \ell\in [m^*,n^*))
(\eta_{\varepsilon(\beta,\ell)}\rest (i{+}1)=(\eta_{\varepsilon(\alpha,\ell)}
\rest i)\hat{\ }\langle\gamma_\ell\rangle)\}\\
\end{array}
\]
is in the ideal $J$.
\end{description}
\end{claim}

\noindent{\em Why?}\ \ \ For each $i<\delta$, $i\geq i^*$ and distinct
sequences $\nu_{m^*},\ldots,\nu_{n^*-1}\in T_i$ define
\[\begin{array}{ll}
B_{\langle\nu_\ell:\ell\in [m^*,n^*)\rangle}\stackrel{\rm def}{=}
\{\bar{\gamma}: &\bar{\gamma}=\langle\gamma_\ell: \ell\in [m^*,n^*)\rangle
\mbox{ and}\\
\ &\mbox{for arbitrarly large }\alpha<\lambda\mbox{ for all }m^*\leq\ell<n^*\\ 
\ &\nu_\ell\hat{\ }\langle\gamma_\ell\rangle\vartriangleleft
\eta_{\varepsilon(\alpha,\ell)}\}.\\
\end{array}\]
We will call a sequence $\langle\nu_\ell: \ell\in [m^*,n^*)\rangle$ {\em a
success} if
\[(\exists^{I_{\nu_{m^*}}}\gamma_{m^*})\ldots(\exists^{I_{\nu_{n^*-1}}}
\gamma_{n^*-1})(\langle\gamma_\ell: \ell\in [m^*,n^*)\rangle\in
B_{\langle\nu_\ell\in [m^*,n^*)\rangle}).\]
Using this notion we may reformulate $(\boxtimes)$ (which we have to prove) to
\begin{description}
\item[$(\boxtimes^*)$] for every large enough $\alpha<\lambda$, for
$J$-majority of $i<\delta$, $i>i^*$ the sequence
$\langle\eta_{\varepsilon(\alpha,\ell)} \rest i: \ell\in [m^*,n^*)\rangle$ is
a success. 
\end{description}

\noindent To show $(\boxtimes^*)$ note that if a sequence
$\langle\nu_\ell:\ell\in [m^*,n^*)\rangle$ is not a success then there are
functions $f^k_{\langle\nu_\ell: \ell\in [m^*,n^*)\rangle}$ (for $m^*\leq
k<n^*$) such that
\[f^k_{\langle\nu_\ell: \ell\in [m^*,n^*)\rangle}:\prod_{\ell=m^*}^{k-1}
\lambda_{\nu_\ell}\ \longrightarrow\ I_{\nu_k}\quad\mbox{ and}\]
\[\begin{array}{ll}
\mbox{if }&\langle\gamma_\ell:\ell\in [m^*,n^*)\rangle\in B_{\langle\nu_\ell:
\ell\in [m^*,n^*)\rangle}\\
\mbox{then }&(\exists k\in [m^*,n^*))(\gamma_k\in f^k_{\langle\nu_\ell:\ell\in
[m^*,n^*)\rangle}(\gamma_{m^*},\ldots,\gamma_{k-1})).\\
\end{array}\] 
If $\langle\nu_\ell: \ell\in [m^*,n^*)\rangle$ is a success then we declare
that $f^k_{\langle\nu_\ell: \ell\in [m^*,n^*)\rangle}$ is constantly equal to
$\emptyset$. 

\noindent Now we may finish the proof of the claim applying clause {\bf (b)} of
definition \ref{cofinal}(5) to $n^*-1$ and functions $F_0,\ldots,F_{n^*-1}$
such that for $k\in [m^*,n^*)$
\[F_k(\nu_0\hat{\ }\langle\gamma_0\rangle,\ldots, \nu_{k-1}\hat{\ }\langle
\gamma_{k-1}\rangle,\nu_k,\ldots,\nu_{n^*-1}\rangle)=
f^k_{\langle\nu_\ell: \ell\in [m^*,n^*)\rangle}(\gamma_{m^*},\ldots,
\gamma_{k-1}).\]
This gives us a suitable $\alpha^*<\lambda$. Suppose $\varepsilon(\alpha,m^*)
\geq\alpha^*$. Then for $J$-majority of $i<\delta$ for each $k\in [m^*,n^*)$
we have 
\medskip

if 
\[F_m(\eta_{\varepsilon(\alpha,0)}\rest (i+1),\ldots,
\eta_{\varepsilon(\alpha,k-1)}\rest (i+1),\eta_{\varepsilon(\alpha,k)}\rest i,
\ldots,\eta_{\varepsilon(\alpha,n^*-1)}\rest i )\in
I_{\eta_{\varepsilon(\alpha,k)\rest i}}\]
then 
\[\eta_{\varepsilon(\alpha,k)}\rest (i+1)\notin
F_m(\eta_{\varepsilon(\alpha,0)}\rest (i+1),\ldots,
\eta_{\varepsilon(\alpha,k-1)}\rest (i+1),\eta_{\varepsilon(\alpha,k)}\rest i
,\ldots,\eta_{\varepsilon(\alpha,n^*-1)}\rest i ).\]
\smallskip

\noindent But the choice of the functions $F_k$ implies that thus for
$J$-majority of $i<\delta$, for each $k\in [m^*,n^*)$
\[\eta_{\varepsilon(\alpha,k)}(i)\notin
f^k_{\langle\eta_{\varepsilon(\alpha,\ell)}\rest i: \ell\in [m^*,n^*)\rangle}
(\eta_{\varepsilon(\alpha,m^*)}(i),\ldots,\eta_{\varepsilon(\alpha,
k-1)}(i)).\]
Now the definition of the function $f^k_{\langle\nu_\ell: \ell\in
[m^*,n^*)\rangle}$ works:

if for some relevant $i<\delta$ above the sequence
$\langle\eta_{\varepsilon(\alpha,\ell)}\rest i: \ell\in [m^*,n^*)\rangle$ is
not a success then $\langle\eta_{\varepsilon(\alpha,\ell)}(i): \ell\in
[m^*,n^*)\rangle\notin B_{\langle\eta_{\varepsilon(\alpha,\ell)}\rest i:
\ell\in [m^*,n^*)\rangle}$ and this contradicts $\hat{*}$ before.

The claim is proved.
\medskip

Let $\alpha^*$ be such that for each $\alpha\geq\alpha^*$ we have $Z_\alpha\in
J$. Choose $i\in\delta\setminus Z_{\alpha^*}$ such that the clause
$(\circleddash)$ applies for $n^*-m^*$ and $i$. Let
\[Y\stackrel{\rm def}{=}\{\langle\gamma_{m^*},\ldots,\gamma_{n^*-1}\rangle\!:
(\exists^\lambda\beta)(\forall \ell{\in}[m^*,n^*))(\eta_{\varepsilon
(\beta,\ell)}\rest (i{+}1){=} (\eta_{\varepsilon(\alpha^*,\ell)}\rest i)\hat{\
}\langle\gamma_\ell\rangle)\}.\]
The definition of $Z_{\alpha^*}$ (and the choice of $i$) imply that the
assumption $(\circleddash)$ applies to the set $Y$ and we get
$\gamma_\ell^\prime,\gamma_\ell^{\prime\prime}<
\lambda_{\eta_{\varepsilon(\alpha^*,\ell)}\rest i}$ (for $m^*\leq\ell<n^*$)
such that 
\[\langle\gamma^\prime_\ell:m^*\leq\ell<n^*\rangle,\langle
\gamma^{\prime\prime}_\ell:m^*\leq\ell<n^*\rangle\in Y\quad\mbox{ and}\]
\[\ba_{\eta_{\varepsilon(\alpha^*,\ell)}\rest i}\models y_{\eta_{\varepsilon(
\alpha^*,\ell)}\rest i\hat{\ }\langle\gamma^\prime_\ell\rangle} \cap
y_{\eta_{\varepsilon(\alpha^*,\ell)}\rest i\hat{\ }\langle
\gamma^{\prime\prime}_\ell\rangle}=\bz\quad\mbox{for }m^*\leq\ell< n^*.\]
Now choose $\alpha<\beta<\lambda$ such that for $m^*\leq\ell<n^*$
\[\eta_{\varepsilon(\alpha^*,\ell)}\rest i\hat{\ }\langle\gamma_\ell^\prime
\rangle = \eta_{\varepsilon(\alpha,\ell)}\rest (i+1),\quad
\eta_{\varepsilon(\alpha^*,\ell)}\rest i\hat{\ }\langle
\gamma_\ell^{\prime\prime}\rangle=\eta_{\varepsilon(\beta,\ell)}\rest (i+1)\]
(possible by the choice of $Y$ and $\gamma^\prime_\ell,
\gamma^{\prime\prime}_\ell$). The definition of the algebra $\ba^{\green}(\C)$
and the choice of $\gamma^\prime_\ell,\gamma^{\prime\prime}_\ell$ imply that
for $m^*\leq\ell<n^*$
\[\ba^{\green}(\C)\models x_{\varepsilon(\alpha,\ell)}\cap
x_{\varepsilon(\beta,\ell)}\neq \bz.\]
If $\ell\neq m$ then 
\[\ba^{\green}(\C)\models x_{\varepsilon(\alpha,\ell)}\cap
x_{\varepsilon(\beta,m)}\neq \bz\]
by the conditions $\hat{(**)}$ and 4) of the preliminary cleaning (and the
definition of $\ba^{\green}(\C)$, remember $z_\alpha\neq\bz$). Finally,
remembering that $\varepsilon(\alpha,\ell)=\varepsilon(\beta,\ell)$ for
$\ell<m^*$, $z_\alpha\neq\bz$ and $z_\beta\neq\bz$, we may conclude that
\[\ba^{\green}(\C)\models \bigcap_{\ell<n^*} x_{\varepsilon(\alpha,\ell)}\cap 
\bigcap_{\ell<n^*} x_{\varepsilon(\beta,\ell)}\neq \bz\]
finishing the proof. \QED

\begin{theorem}
\label{main}
{\em If} $\mu$ is a strong limit singular cardinal, $\lambda\stackrel{\rm
def}{=}2^\mu=\mu^+$

\noindent {\em then} there are Boolean algebras $\ba_1,\ba_2$ such that
the algebra $\ba_1$ satisfies the $\lambda$-cc, the algebra $\ba_2$ has the
$(2^{\cf(\mu)})^+$--Knaster property but the free product $\ba_1*\ba_2$ does
not satisfy the $\lambda$-cc.   
\end{theorem}

\Proof Let $\delta=\cf(\mu)$ and let $h:\delta\longrightarrow\omega$ be a
function such that
\[(\forall n\in\omega)(\exists^\delta i)(h(i)=n).\]
Choose an increasing sequence $\langle\mu_i:i<\delta\rangle$ of regular
cardinals such that $\mu=\sum\limits_{i<\delta}\mu_i$. Next, by induction on
$i<\delta$ choose $\lambda_i$, $\chi_i$, $(\ba_i,\bar{y}_i)$ and $I_i$ such
that
\begin{enumerate}
\item $\lambda_i,\chi_i$ are regular cardinals below $\mu$
\item $\lambda_i>\chi_i\geq\prod\limits_{j<i}\lambda_j + \mu_i$
\item $I_i$ is a $\chi_i^+$-complete ideal on $\lambda_i$ (containing all
singletons) 
\item $(\ba_i,\bar{y}_i)$ is a $\lambda_i$-marked Boolean algebra such that

\noindent {\em if} $n=h(i)$ and the set $Y\subseteq (\lambda_i)^{n+1}$ is such
that  
\[(\exists^{I_i}\gamma_0)\ldots(\exists^{I_i}\gamma_n)(\langle
\gamma_0,\ldots,\gamma_n\rangle\in Y)\]

\noindent{\em then} for some $\gamma_\ell^\prime,\gamma_\ell^{\prime\prime}<
\lambda_i$ (for $\ell\leq n$) we have
\[\langle\gamma_\ell^\prime:\ell\leq n\rangle, \langle\gamma_\ell^{\prime
\prime}:\ell\leq n\rangle\in Y\quad\mbox{ and for all }\ell\leq n\]
\[\ba_{i}\models y^i_{\gamma^\prime_\ell}\cap y^i_{\gamma^{\prime
\prime}_\ell}=\bz.\]
\item Each algebra $\ba_i$ satisfies the $(2^{|\delta|})^+$--Knaster
condition. 
\end{enumerate}
Arriving at the stage $i$ of the construction first we put $\chi_i=
(\prod\limits_{j<i}\lambda_j+\mu_i)^+$. Next we define inductively
$\chi_{i,k},\lambda_{i,k}$ for $k\leq h(i)$ such that
\[\chi_{i,0}=\chi_i,\quad \lambda_{i,k}=(2^{\chi_{i,k}})^+,\quad
\chi_{i,k+1}=(\lambda_{i,k})^+.\] 
By \ref{algide}, for each $k\leq h(i)$ we find a
$(\lambda_{i,k},\chi_{i,k}^+)$--well marked Boolean algebra
$(\ba_{i,k},\bar{y}_{i,k},I_{i,k})$ such that $\ba_{i,k}$ has the
$(2^\delta)^+$--Knaster property (compare \ref{well2}). Let
$\lambda_i=\lambda_{i,h(i)}$. Proposition \ref{finite} applied to
$\langle(\ba_{i,k},\bar{y}_{i,k},I_{i,k}): k\leq h(i)\rangle$ provides a
$\lambda_i$-marked Boolean algebra $(\ba_i,\bar{y}_i)$ and a
$\chi_i^+$-complete ideal $I_i$ on $\lambda_i$ such that the requirements 4,5
above are satisfied.
\medskip

Now put $T=\bigcup\limits_{j<\delta}\prod\limits_{i<j}\lambda_i$ and for
$\eta\in T$:
\[\ba_\eta=\ba_{\lh(\eta)},\ \bar{y}_\eta=\bar{y}_{\lh(\eta)},\
I_\eta=I_{\lh(\eta)}.\]
By \ref{stronger} we find a stronger $J^{\bd}_\delta$-cofinal sequence
$\bar{\eta}=\langle \eta_\alpha: \alpha<\lambda\rangle$ for
$(T,\bar{\lambda},\bar{I})$. Take the $(\delta,\mu,\mu^+)$-conctructor $\C$
determined by these parameters. Look at the algebras $\ba_2=\ba^{\red}(\C)$,
$\ba_1=\ba^{\green}(\C)$. Applying \ref{notcc} we get that $\ba_1 *\ba_2$
fails the $\lambda$-cc. The choice of the function $h$ and the requirement 4
above allow us to apply \ref{green} to conclude that the algebra $\ba_2$
satisfies $\lambda$-cc. Finally, by \ref{consknaster}, we have that $\ba_1$
has the $(2^\delta)^+$--Knaster property. \QED

\begin{remark}
{\em
\begin{enumerate}
\item We shall later give results not using $2^\mu=\mu^+$ but still not in ZFC
\item Applying the methods of \cite{Sh:576} one can the consistency of: 
{\em for some $\mu$ strong limit singular there is no example for
$\lambda=\mu^+$}.
\item If we want ``for no regular $\lambda\in [\mu,2^\mu]$'' more is needed, we
expect the consistency, but it is harder (not speaking of ``for all $\mu$'')
\item Remark 1) above shows that $2^\mu>\mu^+$ is not enough for the negative
result.  
\end{enumerate}
}
\end{remark}

\section{Toward improvements}

\begin{definition}
\label{super}
Let $(T,\bar{\lambda},\bar{I})\in\K_{\mu,\delta}^{\id}$ and let $J$ be an
ideal on $\delta$ (including $J^{\bd}_\delta$, as usual). We say that a
sequence $\bar{\eta}=\langle\eta_\alpha: \alpha<\lambda\rangle$ of
$\delta$-branches through $T$ is {\em super $J$-cofinal for}
$(T,\bar{\lambda},\bar{I})$ if 
\begin{description}
\item[(a)] $\eta_\alpha\neq\eta_\beta$ for distinct $\alpha,\beta<\lambda$
\item[(b)] for every function $F$ there is $\alpha^*<\lambda$ such that  

\noindent {\em if} $\alpha_0<\ldots<\alpha_n<\lambda$, $\alpha^*\leq\alpha_n$

\noindent {\em then} the set
\[\hspace{-1.3cm}\begin{array}{ll}
\{i<\delta:&\mbox{(ii)$^*$ }F(\eta_{\alpha_0},\ldots,\eta_{\alpha_{n-1}},
\eta_{\alpha_n}\rest i)\in I_{\eta_{\alpha_n}\rest i}\\ 
\ &\ \ \mbox{(and well defined) but }\\
\ &\ \ \eta_{\alpha_n}\rest(i{+}1)\in F(\eta_{\alpha_0},\ldots, 
\eta_{\alpha_{n-1}},\eta_{\alpha_n}\rest i)\}\\
\end{array}\] 
is in the ideal $J$.
\end{description}
\end{definition}

\begin{remark}
{\em
\begin{enumerate}
\item The main difference between the definition of super $J$-cofinal sequence
and those in \ref{cofinal} is the fact that here the values of the function
$F$ depend on $\eta_{\alpha_\ell}$ (for $\ell<n$), not on the restrictions of 
these sequences as it was in earlier notions.
\item ``super$^*$ $J$--cofinal'' is defined by adding
``$\alpha^*\leq\alpha_0$'' (compare \ref{cofinal}(10)). 
\end{enumerate}
}
\end{remark}

\begin{proposition}
\label{supimpstr}
Suppose that $(T,\bar{\lambda},\bar{I})\in \K^{\id}_{\mu,\delta}$ is such that 
for each $\nu\in T_i$, $i<\delta$ the ideal $I_\nu$ is $|T_i|^+$-complete. Let
$J\supseteq J^{\bd}_\delta$ be an ideal on $\delta$. Then every super
$J$-cofinal sequence is stronger$^*$ $J$-cofinal.
\end{proposition}

\Proof Assume that $\bar{\eta}=\langle\eta_\alpha:\alpha<\lambda\rangle
\subseteq\lim_\delta(T)$ is super $J$-cofinal for $(T,\bar{\lambda},\bar{I})$.
Let $n<\omega$ and let $F_0,\ldots,F_{n-1}$ be functions. For each $\ell\leq
n$ we define an $(\ell+1)$--place function $F^*_\ell$ such that

if $\alpha_0<\alpha_1<\ldots<\alpha_{\ell-1}<\lambda$, $\rho\in
T_i$, $i<\delta$

then 
\[\begin{array}{ll}
F^*_\ell(\eta_{\alpha_0},\ldots,\eta_{\alpha_{\ell-1}},\rho)=
&\bigcup\{F_\ell(\eta_{\alpha_0}\rest(i{+}1),\ldots,\eta_{\alpha_{\ell-1}}\rest
(i{+}1), \rho,\nu_{\ell+1},\ldots,\nu_n):\\
\ &\nu_{\ell+1},\ldots,\nu_n\in T_i\ \&\\
\ &F_\ell(\eta_{\alpha_0}\rest(i{+}1),\ldots,\eta_{\alpha_{\ell-1}}\rest
(i{+}1), \rho,\nu_{\ell+1},\ldots,\nu_n)\in I_\rho\\
\ &\mbox{(and well defined)}\}.\\
  \end{array}\]
As the ideals $I_\rho$ (for $\rho\in T_i$) are $|T_i|^+$-complete we know that
\[F^*_\ell(\eta_{\alpha_0},\ldots,\eta_{\alpha_{\ell-1}},\rho)\in I_\rho.\]
Applying \ref{super}(b) to the functions $F^*_\ell$ ($\ell<n$) we choose
$\alpha^*_\ell<\lambda$ such that

if $\alpha_0<\ldots<\alpha_\ell<\lambda$, $\alpha^*_\ell\leq\alpha_\ell$

then the set
\[\begin{array}{ll}
B^*_\ell\stackrel{\rm def}{=}\{i<\delta: &F^*_\ell(\eta_{\alpha_0},\ldots,
\eta_{\alpha_{\ell-1}},\eta_{\alpha_\ell}\rest i)\in I_{\eta_{\alpha_\ell}
\rest i}\ \mbox{ but}\\
\ &\eta_{\alpha_\ell}\rest (i+1)\in F^*_\ell(\eta_{\alpha_0},\ldots,
\eta_{\alpha_{\ell-1}},\eta_{\alpha_\ell}\rest i)\}\\ 
  \end{array}\]

is in the ideal $J$.

\noindent Put $\alpha^*=\max\{\alpha^*_\ell:\ell\leq n\}$. We want to show
that this $\alpha^*$ works for the condition \ref{cofinal}(6)(b) (version for
``stronger$^*$''). So suppose that $m\leq n$, $\alpha^*\leq\alpha_0<\alpha_1<
\ldots<\alpha_n<\lambda$. Let
\[\begin{array}{ll}
B_m\stackrel{\rm def}{=}&\{i<\delta: F_m(\eta_{\alpha_0}\rest(i{+}1),\ldots,
\eta_{\alpha_{m-1}}\rest (i{+}1),\eta_{\alpha_m}\rest i,\ldots,\eta_{\alpha_n}
\rest i)\in I_{\eta_{\alpha_m}\rest i}\ \&\\
\ &\eta_{\alpha_m}\rest (i{+}1)\in F_m(\eta_{\alpha_0}\rest(i{+}1),\ldots,
\eta_{\alpha_{m-1}}\rest (i{+}1),\eta_{\alpha_m}\rest i,\ldots,\eta_{\alpha_n}
\rest i)\}.\\
  \end{array}\]
Note that if $i\in B_m$ then, as $\alpha^*_m\leq\alpha^*\leq\alpha_m$,
\[\begin{array}{rl}
\eta_{\alpha_m}\rest (i{+}1)\in &F_m(\eta_{\alpha_0}\rest(i{+}1),\ldots,
\eta_{\alpha_{m-1}}\rest (i{+}1),\eta_{\alpha_m}\rest i,\ldots,\eta_{\alpha_n}
\rest i)\subseteq\\
\subseteq &F^*_m(\eta_{\alpha_0},\ldots,\eta_{\alpha_{m-1}},\eta_{\alpha_m}
\rest i)\in I_{\eta_{\alpha_m}\rest i}. 
  \end{array}\]
Hence we conclude that $B_m\subseteq B^*_m$ and therefore $B_m\in J$,
what finishes the proof of the proposition. \QED

\begin{proposition}
\label{supimpfub}
Assume that $(T,\bar{\lambda},\bar{I})\in\K_{\mu,\delta}^{\id}$, each ideal
$I_\eta$ (for $\eta\in T_i$, $i<\delta$) is $(|\delta|+|T_i|)^+$-complete and
$J\supseteq J^{\bd}_\delta$ is an ideal on $\delta$. Further suppose that a
sequence $\bar{\eta}=\langle\eta_\alpha:\alpha<\lambda\rangle$ is super
$J$-cofinal for $(T,\bar{\lambda},\bar{I})$, $\lambda$ is a regular cardinal
greater than $|T|$ and a sequence $\langle\alpha_{\varepsilon,\ell}:
\varepsilon<\lambda, \ell<n\rangle\subseteq \lambda$ is with no repitition and 
such that 
\[\alpha_{\varepsilon,0}<\alpha_{\varepsilon,1}<\ldots<\alpha_{\varepsilon,
n-1}\quad\quad\mbox{ for all }\varepsilon<\lambda.\]
{\em Then} for every $\varepsilon<\lambda$ large enough there is $a\in J$ such
that 
\begin{description}
\item[$(\boxdot)$] if $i_\ell\in\delta\setminus a$ (for $\ell<n$), $i_0\geq
i_1\geq\ldots\geq i_{n-1}$ then
\[\hspace{-1cm}\begin{array}{l}
(\exists^{I_{\eta_{\alpha_{\varepsilon,0}}\rest i_0}}\gamma_0)\ldots
(\exists^{I_{\eta_{\alpha_{\varepsilon,n-1}}\rest i_{n-1}}}\gamma_{n-1})\\
(\exists^\lambda \zeta<\lambda)(\forall\ell<n)(\eta_{\alpha_{\zeta,\ell}}\rest
(i_\ell+1) = \eta_{\alpha_{\varepsilon,\ell}}\rest i_\ell\hat{\
}\langle\gamma_\ell\rangle).\\
\end{array}\] 
\end{description}
\end{proposition}

\Proof This is very similar to claim \ref{claimgreen}. First choose
$\varepsilon_0<\lambda$ such that for each $\varepsilon\in
[\varepsilon_0,\lambda)$ and for every $i_0,\ldots,i_{n-1}<\delta$ we have  
\[(\exists^\lambda \zeta<\lambda)(\forall\ell<n)(\eta_{\alpha_{\zeta,\ell}}
\rest(i_\ell+1)=\eta_{\alpha_{\varepsilon,\ell}}\rest (i_\ell+1))\]
(possible as $|T|<\cf(\lambda)=\lambda$).

Now, for $\bar{\imath}=\langle i_\ell: \ell<n\rangle\subseteq\delta$ and
$\bar{\nu}=\langle\nu_\ell: \ell<n\rangle$ such that $i_0\geq i_1\geq\ldots
\geq i_{n-1}$, $\nu_\ell\in T_{i_\ell}$ and $k<n$ we define a function
$f^k_{\bar{\imath},\bar{\nu}}: \prod\limits_{\ell<k} \lambda_{\nu_\ell}
\longrightarrow I_{\nu_k}$ (with a convention that
$f^0_{\bar{\imath},\bar{\nu}}$ is supposed to be a $0$-place function, i.e. a
constant) as follows

Let
\[B_{\bar{\imath},\bar{\nu}}\stackrel{\rm def}{=}\{\langle\gamma_\ell:
\ell<n\rangle\in\prod_{\ell<n}\lambda_{\nu_\ell}: (\exists^\lambda \zeta<
\lambda)(\forall \ell<n)(\eta_{\alpha_{\zeta,\ell}}\rest (i_\ell +1)=\nu_\ell
\hat{\ }\langle\gamma_\ell\rangle)\}.\]
If
\begin{description}
\item[$(\blacklozenge_{\bar{\imath},\bar{\nu}})$] $\neg(\exists^{I_{\nu_0}}
\gamma_0)\ldots(\exists^{I_{\nu_{n-1}}}\gamma_{n-1})(\langle\gamma_0,\ldots,
\gamma_{n-1}\rangle\in B_{\bar{\imath},\bar{\nu}})$
\end{description}
then $f^0_{\bar{\imath},\bar{\nu}},\ldots,f^{n-1}_{\bar{\imath},\bar{\nu}}$
are such that
\begin{description}
\item[$(\lozenge)$] if $\langle\gamma_0,\ldots,\gamma_{n-1}\rangle\in
B_{\bar{\imath},\bar{\nu}}$ 

then $(\exists k<n)(\gamma_k\in f^k_{\bar{\imath},\bar{\nu}}(\gamma_0,\ldots,
\gamma_{k-1}))$. 
\end{description}
Otherwise (i.e. if not $(\blacklozenge_{\bar{\imath},\bar{\nu}})$ the
functions $f^k_{\bar{\imath},\bar{\nu}}$ are constantly equal to $\emptyset$
(for $k<n$). 

Next, for $k<n$, choose functions $F_k$ such that if $\eta_0,\ldots,\eta_k\in
\lim_\delta(T)$, $i<\delta$ then 
\[\begin{array}{ll}
F_k(\eta_0,\ldots,\eta_{k-1},\eta_k\rest i)= &\ \\
\bigcup\{f^k_{\bar{\imath},\bar{\nu}}(\eta_0(i_0),\ldots,\eta_{k-1}(i_{k-1})):
&\bar{\imath}=\langle i_\ell:\ell<n\rangle,\ \bar{\nu}=\langle\nu_\ell:
\ell<n\rangle,\\  
\ &\delta>i_0\geq\ldots\geq i_k=i\geq i_{k+1}\geq\ldots\geq i_{n-1}\\
\ &\nu_\ell=\eta_\ell\rest i_\ell\ \mbox{ for }\ell\leq k\ \mbox{ and}\\
\ &\nu_\ell\in T_{i_\ell}\ \mbox{ for }k<\ell<n\}.
  \end{array}\]
Note that $F_k(\eta_0,\ldots,\eta_{k-1},\eta_k\rest i)$ is a union of at most
$|\delta|+|T_i|$ sets from the ideal $I_{\eta_k\rest i}$ and hence
$F_k(\eta_0,\ldots,\eta_{k-1},\eta_k\rest i)\in I_{\eta_k\rest i}$ (for each
$\eta_0,\ldots,\eta_k\in \lim_\delta(T)$, $i<\delta$). Thus, using the super
$J$-cofinality of $\bar{\eta}$ we find $\alpha^*<\lambda$ such that
\begin{quotation}
\noindent if $\alpha^*\leq\alpha_<\ldots<\alpha_n<\lambda$

\noindent then the set
\[\{i<\delta: (\exists k<n)(\eta_{\alpha_k}(i)\in F_k(\eta_{\alpha_0},\ldots,
\eta_{\alpha_{k-1}},\eta_{\alpha_k}))\}\]
is in the ideal $J$.
\end{quotation}
Let $\varepsilon_1>\varepsilon_0$ be such that for every $\varepsilon\in
[\varepsilon_1,\lambda)$ we have $\alpha^*<\alpha_{\varepsilon,0}<\ldots<
\alpha_{\varepsilon,n-1}$. 

Suppose now that $\varepsilon_1<\varepsilon<\lambda$. By the choice of
$\alpha^*$ we know that the set
\[a\stackrel{\rm def}{=}\{i<\delta: (\exists\ell<n)(\eta_{\alpha_{\varepsilon,
\ell}}(i)\in F_\ell(\eta_{\alpha_{\varepsilon,0}},\ldots\eta_{
\alpha_{\varepsilon,\ell-1}},\eta_{\alpha_{\varepsilon,\ell}}\rest i))\}\]
is in the ideal $J$. We are going to show that the assertion $(\boxdot)$ holds
for $\varepsilon$ and $a$.
\smallskip

\noindent Suppose that $\bar{\imath}=\langle
i_\ell:\ell<n\rangle\subseteq\delta\setminus a$, $i_0\geq i_1\geq\ldots\geq
i_{n-1}$. Let $\bar{\nu}=\langle \nu_\ell: \ell<n\rangle$,
$\nu_\ell=\eta_{\alpha_{\varepsilon,\ell}}\rest i_\ell$. If the condition
$(\blacklozenge_{\bar{\imath},\bar{\nu}})$ fails then we are done. So assume
that it holds true. By the choice of the set $a$ (and $\alpha^*$) we have
\[(\forall\ell<n)(\eta_{\alpha_{\varepsilon,\ell}}(i_\ell)\notin F_\ell
(\eta_{\alpha_{\varepsilon,0}},\ldots,\eta_{\alpha_{\varepsilon,\ell-1}},
\eta_{\alpha_{\varepsilon,\ell}}\rest i_l))\]
what, by the definition of $F_\ell$, implies that
\[(\forall\ell<n)(\eta_{\alpha_{\varepsilon,\ell}}(i_\ell)\notin
f^\ell_{\bar{\imath},\bar{\nu}}(\eta_{\alpha_{\varepsilon,0}}(i_0),\ldots,
\eta_{\alpha_{\varepsilon,\ell-1}}(i_{\ell-1}))).\]
By $(\lozenge)$ we conclude that
\[\langle\eta_{\alpha_{\varepsilon,0}}(i_0),\ldots,\eta_{\alpha_{\varepsilon,
n-1}}(i_{n-1})\rangle\notin B_{\bar{\imath},\bar{\nu}}\]
and hence, by the definition of $B_{\bar{\imath},\bar{\nu}}$,
\[\neg(\exists^\lambda \zeta)(\forall\ell<n)(\eta_{\alpha_{\zeta,\ell}}\rest
(i_\ell+1) =\eta_{\alpha_{\varepsilon,\ell}}\rest (i_\ell))\]
which contradicts the choice of $\varepsilon_0$ (remember $\varepsilon\geq
\varepsilon_1>\varepsilon_0$). \QED

\begin{definition}
We say that a $\lambda$-marked Boolean algebra $(\ba,\bar{y})$ {\em has
character $n$} if
\begin{quotation}
\noindent for every finite set $u\in [\lambda]^{{<}\omega}$ such that
$\ba\models\bigcap\limits_{\alpha\in u} y_\alpha=\bz$ there exist a subset
$v\subseteq u$ of size $|v|\leq n$ such that
$\ba\models\bigcap\limits_{\alpha\in v} y_\alpha=\bz$.
\end{quotation}
\end{definition}

\begin{proposition}
If a $\lambda$-marked Boolean algebra $(\ba,\bar{y})$ is $(\theta,{\rm
not}\lambda)$-Knaster (or other examples considered in the present paper) and
$(\ba,\bar{y})$ has character 2 then without loss of generality
$(\ba,\bar{y})$ is determined by a colouring on $\lambda$:
\begin{quotation}
\noindent if $c:[\lambda]^2\longrightarrow 2$ is such that 
\[c(\{\alpha,\beta\})=0\quad\mbox{ \rm iff }\quad\ba\models y_\alpha\cap
y_\beta=\bz\]

\noindent then the algebra $\ba$ is freely generated by $\{y_\alpha:
\alpha<\lambda\}$ except that
\begin{quotation}
{\em if} $c(\{\alpha,\beta\})=0$

{\em then} $y_\alpha\cap y_\beta=0$. \QED
\end{quotation}
\end{quotation}
\end{proposition}

\begin{remark}
{\em 
These are nice examples.
}
\end{remark}

\begin{proposition}
In all our results (like: \ref{algide} or \ref{consknaster}), the marked
Boolean algebra $(\ba,\bar{y})$ which we get is actually of character 2 as
long as any $(\ba_\eta,\bar{y}_\eta)$ appearing in the assumptions (if any) is
like that.  

\noindent Then automatically the $\theta$--Knaster property of the marked
Boolean algebra $(\ba,\bar{y})$ implies a stronger condition: 

\noindent if $Z\in [\lh(\bar{y})]^\theta$ then there is a set $Y\in
[Z]^\theta$ such that $\{y_i: i\in Y\}$ generates a filter in $\ba$. \QED 
\end{proposition}

\begin{proposition}
\label{ultfub}
Let $(T,\bar{\lambda},\bar{I})\in\K_{\mu,\delta}^{\id}$ be such that for each
$\eta\in T$ the filter $(I_\eta)^c$ (dual to $I_\eta$) is an ultrafilter on
$\suc_T(\eta)$, and let $J$ be an ideal on $\delta$ (extending
$J^{\bd}_\delta$). {\em If:}
\begin{description}
\item[(a)] $\C=(T,\bar{\lambda},\bar{\eta},\langle(\ba_\eta,\bar{y}_\eta): 
\eta\in T\rangle)$ is a $(\delta,\mu,\lambda)$-constructor, the sequence
$\bar{\eta}$ is stronger $J$-cofinal for $(T,\bar{\lambda},\bar{I})$,
$|T|<\cf(\lambda)=\lambda$, 
\item[(b)] the sequence $\langle\alpha_{\varepsilon,\ell}:\varepsilon<\lambda,
\ell<n\rangle\subseteq\lambda$ is with no repetition,
\item[(c)] for each distinct $\eta,\nu\in T$ either the ideal $I_\eta$ is
$(2^{\lambda_\nu})^+$--complete (which, of course, implies $\lambda_\eta>
2^{\lambda_\nu}$) or the ideal $I_\nu$ is $(2^{\lambda_\eta})^+$--complete (it
is enough if this holds true for $\eta,\nu$ such that $\lh(\eta)=\lh(\nu)$
\end{description}
{\em then} for every large enough $\varepsilon<\lambda$ for $J$-almost all
$i<\delta$ there are sets $X_\ell\in (I_{\eta_{\alpha_{\varepsilon,\ell}}\rest
i})^+$ (for $\ell<n$) such that
\[(\forall\gamma_0\in X_0)\ldots(\forall\gamma_{n-1}\in X_{n-1})
(\exists^\lambda\zeta<\lambda)(\forall\ell<n)
(\eta_{\alpha_{\varepsilon,\ell}}\rest i\hat{\ }\langle\gamma_\ell\rangle
\vartriangleleft\eta_{\alpha_{\zeta,\ell}}).\] 
\end{proposition}

\noindent{\bf Remark~\ref{ultfub}.A}\ \ \ We can replace {\em stronger} by
{\em big} and then omit being an ultrafilter.
\medskip

\Proof First note that we may slightly re-enumerate are sequence $\langle
\alpha_{\varepsilon,\ell}: \varepsilon<\lambda,\ell<n\rangle$ and we may
assume that for each $\varepsilon<\lambda$
\[\alpha_{\varepsilon,0}<\alpha_{\varepsilon,1}<\ldots<\alpha_{\varepsilon,
n-1}.\]
Now, since $|T|<\cf(\lambda)=\lambda$ we may apply claim \ref{claimgreen} to
\[\langle\langle\alpha_{\varepsilon,\ell}:\ell<n\rangle:\varepsilon_0\leq
\varepsilon<\lambda\rangle\]
(we need to take $\varepsilon_0$ large enough to get the condition $\hat{(*)}$
of the proof of \ref{green}). Consequently we may conclude that there is
$\varepsilon_1<\lambda$ such that for every $\varepsilon\in
[\varepsilon_1,\lambda)$ 
\begin{description}
\item[$(\boxtimes_\varepsilon)$] for $J$-majority of $i<\delta$ we have 
\[\hspace{-0.7cm} (\exists^{I_{\eta_{\alpha_{\varepsilon,0}}\rest i}}\gamma_0)
\ldots (\exists^{I_{\eta_{\alpha_{\varepsilon,n-1}}\rest i}}\gamma_{n-1})
(\exists^\lambda\zeta{<}\lambda)(\forall\ell{<}n)(\eta_{\alpha_{\zeta,\ell}}
\rest (i{+}1) =\eta_{\alpha_{\varepsilon,\ell}}\rest i\hat{\
}\langle\gamma_\ell \rangle).\]
\end{description}
Now we would like to apply \ref{prefub}. We can not do this directly as we do
not know if the cardinals $\lambda_{\eta_{\varepsilon,\ell}\rest i}$ are
decresing (with $\ell$). However the following claim helps us.

\begin{claim}
\label{transp}
Suppose that $\lambda_0<\lambda_1$ are cardinals and $I_0,I_1$ are maximal
ideals on $\lambda_0,\lambda_1$ respectively. Assume that the ideal $I_1$ is
$(\lambda_0)^+$--complete and $\varphi(x,y)$ is a formula. Then
\[(\exists^{I_0}\gamma_0)(\exists^{I_1}\gamma_1)\varphi(\gamma_0,\gamma_1)\ \
\ \Rightarrow\ \ \ (\exists^{I_1}\gamma_1)(\exists^{I_0}\gamma_0)
\varphi(\gamma_0,\gamma_1).\]  
\end{claim}

\noindent{\em Why?}\ \ \ First note that if $I$ is a maximal ideal then the
quantifiers $\exists^I$ and $\forall^I$ are equivalent. Suppose now that 
\[(\exists^{I_0}\gamma_0)(\exists^{I_1}\gamma_1)\varphi(\gamma_0,\gamma_1).\]
This implies (as $I_0,I_1$ are maximal) that
\[(\forall^{I_0}\gamma_0)(\forall^{I_1}\gamma_1)\varphi(\gamma_0,\gamma_1).\]
Thus we have a set $a\in I_0$ and for each $\gamma\in\lambda_0\setminus a$ we
have a set $b_\gamma\in I_1$ such that
\[(\forall\gamma_0\in\lambda_0\setminus a)(\forall\gamma_1\in\lambda_1
\setminus b_{\gamma_0})\varphi(\gamma_0,\gamma_1).\] 
Let $b=\bigcup\limits_{\gamma\in\lambda_0\setminus a} b_\gamma$. As $I_1$ is
$(\lambda_0)^+$--complete the set $b$ is in $I_1$. Clearly
\[(\forall\gamma_1\in\lambda_1\setminus b)(\forall\gamma_0\in\lambda\setminus
a)\varphi(\gamma_0,\gamma_1)\] 
which implies $(\exists^{I_1}\gamma_1)(\exists^{I_0}\gamma_0)\varphi
(\gamma_0,\gamma_1)$, finishing the proof of the claim.
\medskip

Now fix $\varepsilon>\varepsilon_1$ ($\varepsilon_1$ as chosen earlier). Take
$i^*<\delta$ such that the elements of $\langle\eta_{\alpha_{\varepsilon,
\ell}}\rest i:\ell<n\rangle$ are pairwise distinct. Suppose that $i\in
[i^*,\delta)$ is such that the formula of $(\boxtimes_\varepsilon)$ holds
true. Let $\{k_\ell: \ell<n\}$ be an enumeration of $n$ such that 
\[\lambda_{\eta_{\alpha_{\varepsilon,k_0}}\rest i}>
\lambda_{\eta_{\alpha_{\varepsilon,k_1}}\rest i}>\ldots>
\lambda_{\eta_{\alpha_{\varepsilon,k_{n-1}}}\rest i}.\] 
(Note that by the assumption {\bf (c)} we know that all the
$\lambda_{\eta_{\alpha_{\varepsilon,k_\ell}}\rest i}$ are distinct, remember
the choice of $i^*$.) Applying claim \ref{transp} we conclude that
\[(\exists^{I_{\eta_{\alpha_{\varepsilon,k_0}}\rest i}}\gamma_{k_0})\ldots
(\exists^{I_{\eta_{\alpha_{\varepsilon,k_{n-1}}}\rest i}}\gamma_{k_{n-1}})
(\exists^\lambda\zeta{<}\lambda)(\forall\ell{<}n)(\eta_{\alpha_{\varepsilon,
\ell}} \rest i\hat{\ }\langle\gamma_\ell\rangle=\eta_{\alpha_{\zeta,\ell}}\rest
(i{+}1)).\]
But now we are able to use \ref{prefub} and we get that there are sets
$X_{k_\ell}\subseteq\lambda_{\eta_{\alpha_{\varepsilon,k_\ell}}\rest i}$,
$X_{k_\ell}\notin I_{\eta_{\alpha_{\varepsilon,k_\ell}}\rest i}$ (for
$\ell<n$) such that 
\[\prod_{\ell<n}X_\ell\subseteq\{\langle\gamma_0,\ldots,\gamma_{n-1}
\rangle: (\exists^\lambda\zeta{<}\lambda)(\forall\ell{<}n)
(\eta_{\alpha_{\varepsilon,\ell}}\rest i\hat{\ }\langle\gamma_\ell\rangle=
\eta_{\alpha_{\zeta,\ell}}\rest (i{+}1))\}\]
what is exactly what we need. \QED
\medskip

If we assume less completeness of the ideals $I_\eta$ in \ref{ultfub} then
still we may say something.

\begin{proposition}
Let $\langle\sigma_i: i<\delta\rangle$ be a sequence of cardinals. Suppose
that $T,\bar{\lambda},\bar{I},\bar{\eta},J,\lambda,\mu,\delta$ and
$\langle\alpha_{\varepsilon,\ell}:\varepsilon<\lambda,\ell<n\rangle$ are as in
\ref{ultfub} but with condition {\bf (c)} replaced by
\begin{description}
\item[(c)$^-_{\langle\sigma_i: i<\delta\rangle}$] if $\eta,\nu\in T_i$,
$\eta\neq\nu$, $i<\delta$ then either the ideal $I_\eta$ is
$((\lambda_\nu)^{\sigma_i})^+$--complete or the ideal $I_\nu$ is
$((\lambda_\eta)^{\sigma_i})^+$--complete. 
\end{description}
Then for every large enough $\varepsilon<\lambda$ for $J$-almost all
$i<\delta$ there are sets $X_\ell\in
[\lambda_{\eta_{\alpha_{\varepsilon,\ell}} \rest i}]^{\sigma_i}$ (for
$\ell<n$) such that
\[(\forall\gamma_0\in X_0)\ldots(\forall\gamma_{n-1}X_{n-1})(\exists^\lambda
\zeta<\lambda)(\forall\ell<n)(\eta_{\alpha_{\varepsilon,\ell}}\rest i\hat{\ }
\langle\gamma_\ell\rangle\vartriangleleft\eta_{\alpha_{\zeta,\ell}}).\] 
\end{proposition}

\Proof The proof goes exactly as the one of \ref{ultfub}, but instead of
\ref{prefub} we use \ref{weakfub}. \QED

\begin{remark}
{\em
\begin{enumerate}
\item Note that in the situation as in \ref{ultfub}, we usually have that
``$J$--cofinal'' implies ``stronger $J$--cofinal'' (see \ref{cofimpstrong},
\ref{strimpstr}).
\item The first assumption of \ref{ultfub} (ultrafilters) coupled with our
normal completeness demands is a very heavy condition, but it has rewards.
\item A natural context here is when $\langle\mu_i: i\leq\kappa\rangle$ is a
strictly increasing continues sequence of cardinals such that each $\mu_{i+1}$
is compact and $\mu=\mu_\kappa$. Then every $\mu_{i+1}$-complete filter can be
extended to an $\mu_{i+1}$-complete ultrafilter. Moreover $2^\mu=\mu^+$
follows by Solovay \cite{So74}.

{\em If} for some function $f$ from cardinals to cardinals, for each $\chi$
there is an algebra $\ba_\chi$ of cardinality $f(\chi)$ which cannot be
decomposed into $\leq\mu$ sets $X_i$ each with some property ${\bf
Pr}(\ba_\chi, X_i)$ and if each $\mu_i$ if $f$-inaccessible

\noindent {\em then} we can find $T,\bar{I},\bar{\lambda}$ as in \ref{ultfub}
and such that $\eta\in T_i\ \ \Rightarrow\ \ \mu_i<\chi_\eta<\lambda_\eta<
\mu_{i+1}$ and for $\eta\in T_i$ there is an algebra $\ba_\eta$ with universe
$\lambda_\eta$ and the ideal $I_\eta$ is $\chi_\eta$--complete, 
\[\mbox{if } X\subseteq \ba_\eta\mbox{ and } {\bf Pr}(\ba_\eta,X)\mbox{ then }
X\in I_\eta\]
(compare \ref{algide}) and $\lambda_\eta<\lambda_\nu\ \ \Rightarrow\ \
(2^{\lambda_\eta})^+<\chi_\nu$. Now choosing cofinal $\bar{\eta}$ we may
proceed as in earlier arguments. 
\item It seems to be good for building nice examples, however we did not find
the right question yet.
\item Central to our proofs is an assumption that 
\[\mbox{``}\langle\alpha_{\zeta,\ell}: \zeta<\lambda,\ell<n\rangle\subseteq
\lambda\mbox{ is a sequence with no repetition'',}\]
i.e. we deal with $\lambda$ disjoint $n$-tuples. This is natural as the
examples constructed here are generated from $\{x_i: i<\lambda\}$ by finitary
functions. One may ask what happens if we admit functions with, say,
$\aleph_0$ places? We can still try to get for $\mu$ as above that:
\begin{description}
\item[$(\boxtimes)$] there is $h:[\mu^+]^2\longrightarrow 2$ such that 

{\em if} $\langle u_\varepsilon:\varepsilon<\lambda\rangle$ are pairwise
disjoint, $u_\varepsilon=\{\alpha_{\varepsilon,\ell}:\ell<\ell^*\}$ is the
increasing (with $\ell$) enumeration, $\ell^*<\mu$ ($\ell^*$ infinite), for a
sequence $\langle\nu_\ell:\ell<\ell^*\rangle\subseteq T_i$ 
\[\hspace{-1cm}\begin{array}{l}
B_{\langle\nu_\ell:\ell<\ell^*\rangle}\stackrel{\rm def}{=}\\
\{\langle\eta_{\alpha_{\varepsilon,\ell}}(i): \ell<\ell^*\rangle:
(\exists^\lambda\zeta<\lambda)(\forall\ell<\ell^*)(\eta_{\alpha_{\varepsilon,
\ell}}\rest(i+1)=\eta_{\alpha_{\zeta,\ell}}\rest (i+1))\},\\
\end{array}\]
for some $i^*<\delta$ there are no repetitions in $\langle
\eta_{\alpha_{\varepsilon,\ell}}\rest i^*: \ell<\ell^*\rangle$ and $h\rest
[u_\varepsilon]^2\equiv 1$ (for each $\varepsilon<\lambda$)

{\em then} there are $\alpha<\beta$ (really a large set of these) such that
\[h\rest [u_\alpha\cup u_\beta]^2\equiv 1.\]
\end{description}
The point is that we can deal with functions with infinitely many variables.
Looking at previous proofs, ``in stronger'' we can get (for $\mu$ strong
limit singular etc):
\begin{quotation}
\noindent for $\alpha$ large enough

\noindent for $i<\delta=\cf(\mu)$ large enough

\noindent $\ldots\ldots\ldots$

\noindent we can defeat
\[\hspace{-0.5cm}(\ldots\ldots(\forall^{I_{\eta_{\alpha_{\varepsilon,\ell}}
\rest i}}\gamma_\ell)\ldots\ldots)(\langle\gamma_\ell: \ell<\ell^*\rangle\in
B_{\langle\eta_{\alpha_{\varepsilon,\ell}}\rest i: \ell<\ell^*\rangle})\]
but the duality of quantifiers fails, so the conclusion is only that
\[\hspace{-1cm}(\forall^J i<\delta)[\neg(\ldots
(\forall^{I_{\eta_{\alpha_{\varepsilon,\ell}}\rest i}} 
\gamma_\ell)\ldots)_{\ell<\ell^*}(\langle \eta_{\alpha_{\varepsilon,\ell}}(i):
\ell<\ell^*\rangle\notin B_{\langle\eta_{\alpha_{\varepsilon,\ell}}\rest i:
\ell<\ell^*\rangle})].\] 
\end{quotation}
\item (no ultrafilters) If $I\supseteq J^{\bd}_\eta$, $\delta$ is a regular
cardinal, $\lambda_\eta=\lambda_{\lh(\eta)}$ and for each $u\in [T_i]^{{<}
|\delta|\chi}$, $i<\delta$ the free product $\mathop{\bigstar}\limits_{\eta\in
u}\ba_\eta$ satisfies the $\lambda$-cc then we can show that the algebra
$\ba^{\red}_{<\chi}$ satisfies the $\lambda$-cc too, 

where for a cardinal $\chi$ the algebra $\ba^{\red}_{<\chi}$ is the Boolean
algebra freely generated by
\[\hspace{-0.5cm}\begin{array}{ll}
\{\bigcap\limits_{\alpha\in u} x_\alpha^{\gt(\alpha)}:& \gt:u\longrightarrow
2, u\in [\lambda]^{<\delta}, h\rest [u\cap\gt^{-1}[1]]^2\equiv 1\mbox{ and}\\
\ &|u|<\chi\mbox{ and}\\
\ &(\exists i<\delta)(\mbox{the mapping }\alpha\mapsto\eta_\alpha(i)\mbox{ is
one-to-one (for $\alpha\in u$)})\\
\ &(\exists i<\delta)(\exists\alpha\in u)(\forall j\in (i,\delta))(\forall
\beta\in u)(f_\alpha(j)\leq f_\beta(j))\}.\\
\end{array}\]
[Note that if $\chi\leq\cf(\delta)$ simpler.]
\end{enumerate}
}
\end{remark}

\[*\quad*\quad*\quad*\quad*\]

Now we will deal with an additional demand that the algebra $\ba^{\red}$
satisfies $|\delta|^+$-cc (or even has the $|\delta|^+$--Knaster property).
Note that the demand of $|\delta|$-cc does not seem to be reasonable: if
every $\bar{y}_\eta$ has two disjoint members (and every node $t\in T$ is an
initial segment of a branch through $T$) then we can find $\delta$
branches which, if in $\{\eta_\alpha: \alpha<\lambda\}$, give $\delta$
pairwise disjoint elements. {\em Moreover}:
\begin{quotation}
\noindent for each $\nu\in T_\ell$ let $A_\nu=\{\eta_\alpha(i):
\eta_\alpha\rest i=\nu\}$ and
\[a_\alpha=\{i<\delta: (\exists \beta\in A_{\eta_\alpha\rest i})
(\ba_{\eta_\alpha\rest i}\models y_{\eta_\alpha(i)}\cap y_\beta=\bz)\}.\]
So if $\ba^{\red}\models\sigma$-cc then
$(\forall\alpha<\lambda)(|a_\alpha|<\sigma)$.
\end{quotation}

\begin{definition}
\label{free}
Let $(T,\bar{\lambda})\in \K_{\mu,\delta}$ and let $\bar{\eta}=\langle
\eta_\alpha:\alpha<\lambda\rangle\subseteq{\lim}_\delta(T)$. We say that
$\bar{\eta}$ is {\em hereditary $\theta$-free} if for every $Y\in
[\lambda]^\theta$ there are $Z\in [Y]^\theta$ and $i<\delta$ such that
\[(\forall \alpha,\beta\in Z)(\alpha\neq\beta\ \ \ \Rightarrow\ \ \
[\eta_\alpha\rest i=\eta_\beta\rest i\ \&\ \eta_\alpha(i)\neq\eta_\beta(i)]).\]
\end{definition}

\begin{proposition}
\label{morekna}
Assume that $\C=(T,\bar{\lambda},\bar{\eta},\langle (\ba_\eta,\bar{y}_\eta):
\eta\in T\rangle)$ is a $(\delta,\mu,\lambda)$-constructor.
If $\bar{\eta}$ is hereditary $\theta$-free, each algebra $\ba_\eta$ has
the $\theta$-Knaster property and $\theta$ is regular then the algebra
$\ba^{\red}(\C)$ has the $\theta$-Knaster property.
\end{proposition}

\Proof The same as for \ref{consknaster}. Note that the proof there
shows actually that,
\begin{quotation}
\noindent if $(\forall\alpha<\theta)(|\alpha|^{|\delta|}<\theta=\cf(\theta))$

\noindent then $\bar{\eta}$ is $\theta$-hereditary free. \QED
\end{quotation}

\begin{proposition}
\label{getfree}
Assume that $(T,\bar{\lambda})\in \K_{\mu,\delta}$, $\bar{\eta}=\langle
\eta_\alpha:\alpha<\lambda\rangle\subseteq{\lim}_\delta(T)$, $\lambda$ is a
regular cardinal. Further suppose that 
\begin{description}
\item[(a)] $(\forall\alpha<\theta)(|\alpha|^{{<}\delta}<\theta =\cf(\theta))$,
$\delta<\theta$, $J$ is an ideal on $\delta$ extending $J^{\bd}_\delta$ and
\item[(b)] the sequence $\bar{\eta}$ is $<_J$-increasing and one of the
following conditions is satisfied:
\begin{description}
\item[$(\alpha)$] $\bar{\eta}$ is $<_J$-cofinal in $\prod\limits_{i<\delta}
\lambda_i/J$, $\lambda_i$ are regular cardinals above $\theta$ (at least for
$J$-majority of $i<\delta$), $\{\alpha<\lambda: \cf(\alpha)=\theta\}\in
I[\lambda]$ and $\lambda_\eta=\lambda_{\lh(\eta)}$;  
\item[$(\beta)$] there are a sequence $\langle C_\alpha:\alpha<\lambda\rangle$
of subsets of $\lambda$, a closed unbounded subset $E$ of $\lambda$ and
$i^*<\delta$ such that
\begin{enumerate}
\item $C_\alpha\subseteq\alpha$, $\otp(C_\alpha)\leq\theta$,  
\item if $\beta\in C_\alpha$ then $C_\beta=C_\alpha\cap\beta$ and $\eta_\beta
\rest [i^*,\delta)<\eta_\alpha\rest [i^*,\delta)$,
\item if $\alpha\in E$ and $\cf(\alpha)=\theta$ then $\alpha=\sup(C_\alpha)$.
\end{enumerate}
\end{description}
\end{description}
Then there is $A\in [\lambda]^\lambda$ such that the restriction $\bar{\eta}
\rest A$ is $\theta$-hereditary free. 
\end{proposition}

\Proof First let us assume that the case $(\beta)$ of the clause {\bf (b)} of
the assumptions holds.

\begin{claim}
\label{claimfree}
Suppose that $Y\in [E]^\theta$. Then
\begin{enumerate}
\item $(\exists Z\in [Y]^\theta)(\exists i^\otimes)(\mbox{the sequence
}\langle f_{\beta_\varepsilon}(i^\otimes): \varepsilon\in Z\rangle\mbox{ is
strictly increasing})$. 
\item If additionally $J=J^{\bd}_\delta$ then 
\[\hspace{-0.8cm}(\exists Z\in [Y]^\theta)(\exists i^\otimes<\delta)(\mbox{the
sequence }\langle\eta_\beta\rest [i^\otimes,\delta):\beta\in Z\rangle\mbox{ is
strictly increasing}).\] 
\end{enumerate}
\end{claim}

\noindent{\em Why?}\ \ \ Suppose $Y\in [E]^\theta$. Without loss of generality
we may assume that $\otp(Y)=\theta$. Let $\alpha=\sup(Y)$. So $\alpha\in E$,
$\cf(\alpha)=\theta$ and hence $C_\alpha$ is unbounded in $\alpha$. Let
$C_\alpha=\langle\alpha_\varepsilon:\varepsilon<\theta\rangle$ be the
increasing enumeration. Clearly the set 
\[A\stackrel{\rm def}{=}\{\varepsilon<\theta:[\alpha_\varepsilon,
\alpha_{\varepsilon+1})\cap Y\neq\emptyset\}\]
is unbounded in $\theta$. For $\varepsilon\in A$ choose $\beta_\varepsilon\in
[\alpha_\varepsilon, \alpha_{\varepsilon+1})\cap Y$. Then
\[(\exists a_\varepsilon\in J)(\eta_{\alpha_\varepsilon}\rest (\delta\setminus
a_\varepsilon)\leq \eta_{\beta_\varepsilon}\rest (\delta\setminus
a_\varepsilon)<\eta_{\alpha_{\varepsilon+1}}\rest (\delta\setminus 
a_\varepsilon)).\]  
Now choose $i_\varepsilon\in\delta\setminus a_\varepsilon$, $i_\varepsilon>
i^*$ and find $B\in [A]^\theta$ such that 
\[\varepsilon\in B\ \ \ \Rightarrow\ \ \ i_\varepsilon=i^\otimes.\]
Easily, by the assumption $(\beta)(2)$, this $i^\otimes$ and
$Z=\{\beta_\varepsilon: \varepsilon\in B\}$ are as required in
\ref{claimfree}(1). 

\noindent If additionally we know that $J=J^{\bd}_\delta$ then for some $B\in
[A]^\theta$ we have 
\[(\exists i^\otimes\in [i^*,\delta))(\varepsilon\in B\ \ \Rightarrow\ \
a_\varepsilon\subseteq i^\otimes)\]
and hence the sequence $\langle f_{\beta_\varepsilon}\rest [i^\otimes,\delta):
\varepsilon\in B\rangle$ is as required in \ref{claimfree}(2) (remember
$(\beta)(2)$). 
\medskip

But now, using $i^\otimes$ given by \ref{claimfree} we may deal with the
sequence $\langle f_{\beta_\varepsilon}\rest (i^\otimes +1): \varepsilon\in
B\rangle$ and using the old proof (see \ref{consknaster}) on the tree
$\bigcup\limits_{i\leq i^\otimes} T_i$ (note that we may apply the assumption
{\bf (a)} to arguments like there) we may get the desired conclusion. This
finishes the case when $(\beta)$ of {\bf (b)} holds true.
\bigskip

\noindent Now, assume that the case $(\alpha)$ of the clause {\bf (b)} of
the assumptions holds. We reduce this case to the previous one (using
cofinality). 

\noindent Take $\bar{C}$, $E$ witnessing that the set $\{\alpha<\lambda:
\cf(\alpha)=\theta\}$ is in $I[\lambda]$ and build a $<_J$-increasing sequence
$\bar{\eta}'=\langle\eta_\alpha^\prime:\alpha<\lambda\rangle\subseteq
\prod\limits_{i<\delta}\lambda_i$ such that $\eta_\alpha^\prime>\eta_\alpha$
and $\bar{\eta}'$ satisfies the clause $(\beta)$ of {\bf (b)}  for $\bar{C}$,
$E$. [The construction of $\eta_\alpha^\prime$ is by induction on
$\alpha<\lambda$. Suppose that we have defined $\eta_\beta^\prime$ for
$\beta<\alpha$. Now, at the stage $\alpha$ of the construction, we first
choose $\eta^0_\alpha\in\prod\limits_{i<\delta} \lambda_i$ such that
\[(\forall\beta<\alpha)(\eta_\beta^\prime<_J \eta^0_\alpha).\]
This is possible since the condition $(\alpha)$ implies that $\lambda=
\tcf(\prod_{i<\delta}\lambda_i/J)$, $\alpha<\lambda$. Now we put for
$i<\delta$: 
\[\eta_\alpha^\prime(i)=\max\big\{\eta_\alpha^0(i),\eta_\alpha(i)+1,
\sup\{\eta_\gamma^\prime(i)+1: \gamma\in C_\alpha\}\big\}.\]
Now one can check that this $\bar{\eta}^\prime$ is as required.]
\smallskip

Now we use the fact that $\bar{\eta}$ is cofinal. The set 
\[E'=\{\gamma\in E:(\forall\alpha<\gamma)(\exists\beta<\gamma)
(\eta_\alpha^\prime<_J \eta_\beta)\}\]
is a club of $\lambda$. Look at $\bar{\eta}\rest E'$. Suppose that $Y\in
[E']^\theta$. Without loss of generality we may assume that $\otp(Y)=\theta$
and let $\alpha=\sup(Y)$. By induction on $\varepsilon<\theta$ choose
$\alpha_\varepsilon<\beta_\varepsilon<\gamma_\varepsilon$ such that
\begin{quotation}
$\beta_\varepsilon\in Y$, $\alpha_\varepsilon\in C_\alpha$,
$\gamma_\varepsilon\in C_\alpha$, $\eta_{\alpha_\varepsilon}^\prime<_J
\eta_{\beta_\varepsilon}<_J \eta_{\gamma_\varepsilon}^\prime$ and

if $\zeta<\varepsilon$ then $\gamma_\zeta<\alpha_\varepsilon$.
\end{quotation}
Next choose $i_\varepsilon>i^*$ such that
\[\eta_{\alpha_\varepsilon}^\prime(i_\varepsilon)< \eta_{\beta_\varepsilon}
(i_\varepsilon)<\eta_{\gamma_\varepsilon}^\prime(i_\varepsilon).\]
We may assume that $i_\varepsilon=i^\otimes$ for all $\varepsilon<\theta$.
Now, as $\bar{\eta}^\prime$ obeys $\bar{C}$, we have
\[\zeta<\varepsilon\ \ \ \Rightarrow\ \ \ \eta_{\gamma_\zeta}^\prime
(i^\otimes)<\eta_{\alpha_\varepsilon}^\prime (i^\otimes)\]
and hence we conclude that the sequence $\langle\eta_{\beta_\varepsilon}
(i^\otimes): \varepsilon<\theta\rangle$ is strictly increasing. 
Now we may finish the proof like earlier. \QED

\begin{conclusion}
\label{moreconclusion}
If $\mu$ is a strong limit singular cardinal, $2^\mu=\mu^+=\lambda$ and
$\neg(\exists 0^{\#})$ or at least 
\[\{\delta<\mu^+: \cf(\delta)=(2^{<\cf(\mu)})^+\}\in I[\lambda]\]
then there is a $(\cf(\mu),\mu,\lambda)$-constructor $\C$ such that 
the algebra $\ba^{\red}(\C)$ has the $(2^{<\cf(\mu)})^+$--Knaster property,
its counterpart $\ba^{\green}(\C)$ is $\lambda$--cc and the free product is
not $\lambda$--cc. 

\noindent [Note that if {\rm GCH} holds then $(2^{<\cf(\mu)})^+=(\cf(\mu))^+$
so the problem is closed then.] 
 \end{conclusion}

\Proof Like \ref{main} using \ref{getfree}, \ref{morekna} instead of
\ref{stronger}, \ref{consknaster}. \QED

\section{The use of pcf}
Assuming that $2^{{<}\kappa}$ is much larger than $\kappa=\cf(\kappa)$
($=\cf(\mu)<\mu$) we may still want to have examples with the $(\kappa^+,{\rm
not}\lambda)$--Knaster property and the non-multiplicativity. Here
\ref{moreconclusion} does not help if GCH holds on an end segment of the
cardinals (and $\neg(\exists 0^\#)$). We try to remedy this.

It is done inductively. So \ref{startind} uses $\cf(\mu)=\aleph_0$ just to
start the induction. We can phrase (a part of) it without this assumption but
in applications we use it for $\cf(\mu)=\aleph_0$. Also \ref{startind}(b)
really needs this condition (otherwise we would have to assume that
$(\forall\alpha<\theta)(|\alpha|^{{<}\delta}<\mu)$). This result says that, if
$\cf(\mu)=\aleph_0$, then we have gotten the $\theta$--Knaster property for
{\em every} regular cardinal $\theta\in\mu\setminus \kappa^+$. 

\begin{definition}
\begin{enumerate}
\item By $\K_{\wmk}$ we will denote the class of all tuples $(\theta,\lambda,
\chi,J)$ such that $\theta<\lambda$, $\chi$ are regular cardinals, $J$ is a
$\chi$--complete ideal on $\lambda$ and there is a $(\lambda,\chi)$--well
marked Boolean algebra $(\ba,\bar{y},J)$ (see \ref{well1}) such that the
algebra $\ba$ satisfies the $\theta$--Knaster property ($\wmk$ stays for
``{\bf w}ell {\bf m}arked {\bf K}naster''). 

When we write $(\theta,\lambda)\in\K_{\wmk}$ we really mean $(\theta,\lambda,
\lambda,J^{\bd}_\lambda)\in \K_{\wmk}$ (what means just that there exists a
$(\theta,{\rm \lambda})$--Knaster marked Boolean algebra). 

\item By $\K_{\swmk}$ ($\swmk$ is for ``{\bf s}equence {\bf m}arked {\bf
K}naster'') we will denote the class of all triples $(\theta,\lambda,\chi)$ of
cardinals such that $\theta<\lambda$ are regular and there is a sequence
$\langle(\ba_\alpha,\bar{y}^\alpha):\alpha<\chi \rangle$ of $\lambda$--marked
Boolean algebras such that (for $\alpha<\chi$) the algebras $\ba_\alpha$ have
the $\theta$--Knaster property, $\bar{y}^\alpha=\langle y^\alpha_i:
i<\lambda\rangle$ and  
\begin{quotation}
\noindent {\em if} $n<\omega$, $\alpha_0<\ldots<\alpha_{n-1}<\chi$ and 
$\beta_{\varepsilon,\ell}<\lambda$ for $\varepsilon<\lambda$, $\ell<n$ are
such that $(\forall\varepsilon_1<\varepsilon_2<\lambda)(\forall \ell<n)
(\beta_{\varepsilon_1,\ell}<\beta_{\varepsilon_2,\ell})$

\noindent{\em then} there are $\varepsilon_1<\varepsilon_2<\lambda$ such that
\[\ell<n\ \ \ \Rightarrow\ \ \ \ba_{\alpha_\ell}\models\mbox{``}
y^{\alpha_\ell}_{\beta_{\varepsilon_1,\ell}}\cap
y^{\alpha_\ell}_{\beta_{\varepsilon_2,\ell}}=\bz\mbox{''}.\] 
\end{quotation}
\end{enumerate}
\end{definition}

\begin{remark}
{\em
\begin{enumerate}
\item On some closure properties of $\K_{\wmk}^\theta\stackrel{\rm def}{=}
\{\lambda: (\theta,\lambda)\in\K_{\wmk}\}$ under pcf see \ref{clospcf}: if
$\lambda_i\in\K^\theta_{\wmk}$ (for $i<\delta$), $\lambda_i>\max\pcf
\{\lambda_j: j<i\}$ and $\lambda\in\pcf\{\lambda_i: i<\delta\}$ and
$(\forall\alpha<\theta)(|\alpha|^{|\delta|}<\theta)$ {\em then}
$\lambda\in\K^\theta_{\wmk}$. 
\item We can replace $\theta$ by a set $\Theta$ of such cardinals, no real
difference. And thus we may consider the class $\K^*_{\wmk}$ of all tuples
$(\Theta,\lambda,\chi,J)$ such that there exists a $(\lambda,\chi)$--well
marked Boolean algebra $(\ba,\bar{y},J)$ with
\[(\forall\theta\in\Theta)(\ba\mbox{ satisfies the }\theta\mbox{--Knaster
property}).\] 
\end{enumerate}
}
\end{remark}

\begin{proposition}
\label{startind}
Assume that $\mu$ is a strong limit singular cardinal, $\aleph_0=\cf(\mu)<\mu$
and $\lambda=2^\mu=\mu^+$.
\begin{description}
\item[(a)] If $(\forall\alpha<\theta)(|\alpha|^{\cf(\mu)}<\theta=\cf(\theta)<
\lambda)$

then $(\theta,\lambda)\in\K_{\wmk}$. Moreover $(\theta,\lambda,2^\lambda)\in
\K_{\swmk}$. 

\item[(b)] If $\cf(\mu)<\theta=\cf(\theta)<\mu$ and $\{\alpha<\lambda:
\cf(\alpha)=\theta\}\in I[\lambda]$

then $(\theta,\lambda)\in \K_{\wmk}$. Moreover $(\theta,\lambda,2^\lambda)\in
\K_{\swmk}$.
\end{description}
\end{proposition}

\Proof This is similar to previous proofs and the first parts of
\ref{startind}{\bf (a)}, {\bf (b)} follow from what we have done already:
\ref{startind}{\bf (a)} is an obvious modification of \ref{knacon};
\ref{startind} is similar, but based on \ref{morekna}, \ref{getfree} (and
\ref{stronger}, \ref{notknaster}) (see below). What we actually have to prove
are the ``moreover'' parts.  We will sketch the proof of it for clause {\bf
(b)} only, modifying the proof of \ref{main}.

As in \ref{main} we choose a function $h:\cf(\mu)\longrightarrow\omega$ such
that for each $n\in\omega$ the preimage $h^{-1}[\{n\}]$ is unbounded (in
$\cf(\mu)$). Next we take an increasing equence $\langle\mu_i:i<\cf(\mu)
\rangle$ of regular cardinals such that $\mu=\sum\limits_{i<\delta}\mu_i$.
Finally (like in \ref{main}) we construct $\lambda_i$, $\chi_i$,
$(\ba_i,\bar{y}_i)$ and $I_i$ such that for $i<\cf(\mu)$:
\begin{enumerate}
\item $\lambda_i,\chi_i<\mu$ are regular cardinals,
\item $\lambda_i>\chi_i\geq\prod\limits_{j<i}\lambda_j + \mu_i$,
$\chi_0> \theta + \mu_0$,
\item $I_i$ is a $\chi_i^+$-complete ideal on $\lambda_i$,
\item $(\ba_i,\bar{y}_i)$ is a $\lambda_i$-marked Boolean algebra such that

\noindent {\em if} $n=h(i)$ and the set $Y\subseteq (\lambda_i)^{n+1}$ is such
that  
\[(\exists^{I_i}\gamma_0)\ldots(\exists^{I_i}\gamma_n)(\langle
\gamma_0,\ldots,\gamma_n\rangle\in Y)\]

\noindent{\em then} for some $\gamma_\ell^\prime,\gamma_\ell^{\prime\prime}<
\lambda_i$ (for $\ell\leq n$) we have
\[\langle\gamma_\ell^\prime:\ell\leq n\rangle, \langle\gamma_\ell^{\prime
\prime}:\ell\leq n\rangle\in Y\quad\mbox{ and for all }\ell\leq n\]
\[\ba_{i}\models y^i_{\gamma^\prime_\ell}\cap y^i_{\gamma^{\prime
\prime}_\ell}=\bz,\]
\item each algebra $\ba_i$ satisfies the $\theta$--Knaster condition,
\item for $\xi<\lambda_i$ the set $[\xi,\lambda_i)$ is not in the ideal $I_i$.
\end{enumerate}
Note that the last requirement is new here. Though we cannot demand that the
ideals $I_i$ extend $I^{\bd}_{\lambda_i}$, the condition (6) above is
satisfied in our standard construction. Note that the ideal from \ref{algide}
has this property if $\lambda$ there is regular. Moreover it is preserved when
the (finite) products of ideals (as in \ref{finite}) are considered. Also, if
$I$ is an ideal on $\lambda$, $A_0\in I$ is such that $|\lambda\setminus A_0|$
is minimal and $A_1\in I^+$ is such that $|A_1|$ is minimal then we can use
either $I\rest A_0$ or $I\rest A_1$. All relevant information is preserved
then (in the first case the condition (6) holds in the second $J^{\bd}_\lambda
\subseteq I$ -- under suitable renaming).

\noindent Now we put $T=\bigcup\limits_{i<\delta}\prod\limits_{j<i}\lambda_j$, 
$\ba_\eta=\ba_{\lh(\eta)}$, $\bar{y}_\eta=\bar{y}_{\lh(\eta)}$, $I_\eta=
I_{\lh(\eta)}$. Applying \ref{stronger} we find a stronger
$J^{\bd}_\delta$-cofinal sequence $\bar{\eta}=\langle \eta_\alpha:
\alpha<\lambda\rangle$ for $(T,\bar{\lambda},\bar{I})$. Due to the requirement
(6) above we may additionally demand that $\bar{\eta}$ is
$<_{J^{\bd}_{\cf(\mu)}}$--increasing cofinal in $\prod\limits_{i<\cf(\mu)}
\lambda_i/J^{\bd}_{\cf(\mu)}$. Let $\langle B_\xi: \xi< 2^\lambda\rangle$ be a
sequence of pairwise almost disjoint elements of $[\lambda]^\lambda$ (i.e.
$|B_\xi\cap B_\zeta|<\lambda$ for distinct $\xi,\zeta<2^\lambda$). For each
$\xi< 2^\lambda$ we may apply \ref{getfree} (the version of {\bf
(b)}($\alpha$)) to the sequence $\langle\eta_\alpha:\alpha\in B_\xi\rangle$
and we find $A_\xi\in [B_\xi]^\lambda$ such that each sequence $\langle
\eta_\alpha:\alpha\in A_\xi\rangle$ is $\theta$--hereditary free. Let
\[\ba^*_\xi=\ba^{\red}(T,\bar{\lambda},\langle\eta_\alpha: \alpha\in
A_\xi\rangle,\langle(\ba_\eta,\bar{y}_\eta): \eta\in T\rangle),\quad\quad
\bar{x}_\xi= \langle x^{\red}_{\alpha}: \alpha\in A_\xi\rangle.\]
Of course each $\ba^*_\xi$ is a subalgebra of $\ba^{\red}(T,\bar{\lambda},
\bar{\eta},\langle(\ba_\eta,\bar{y}_\eta): \eta\in T\rangle)$ (generated by
$\bar{x}_\xi$). By \ref{morekna} and \ref{notknaster} we know that the marked
Boolean algebras $(\ba^*_\xi, \bar{x}_\xi)$ are $(\theta,{\rm
not}\lambda)$--Knaster. To show that they witness $(\theta,\lambda,2^\lambda)
\in \K_{\swmk}$ suppose that $n<\omega$, $\xi_0,\ldots,\xi_{n-1}< 2^\lambda$,
$\beta_{\varepsilon,\ell}<\lambda$ (for $\varepsilon<\lambda$, $\ell<n$) are
such that 
\[(\forall\varepsilon_1<\varepsilon_2<\lambda)(\forall\ell<n)
(\beta_{\varepsilon_1,\ell}< \beta_{\varepsilon_2,\ell})\]
and of course $\{\beta_{\varepsilon,\ell}: \varepsilon<\lambda\}\subseteq
A_{\xi_\ell}$. Since $A_{\xi_\ell}$ are almost disjoint we may assume that 
\[(\forall\varepsilon_1,\varepsilon_2<\lambda)(\forall\ell_1<\ell_2<n)
(\beta_{\varepsilon_1,\ell_1}\neq \beta_{\varepsilon_2,\ell_2}).\]
Further we may assume that we have $i^*<\cf(\mu)$ such that for each
$\varepsilon<\lambda$ the sequences $\eta_{\beta_{\varepsilon,\ell}}\rest i^*$
for $\ell< n$ are pairwise distinct.

\noindent By the choice of $\bar{\eta}$, $T$, $\bar{\lambda}$ etc we may apply
\ref{claimgreen} and conclude that for all sufficiently large
$\varepsilon<\lambda$ the set
\[\begin{array}{r}
Z_\varepsilon=\{i<\cf(\mu): \neg(\exists^{I_{\eta_{\beta_{\varepsilon,0}}\rest
i}}\gamma_{0})\ldots(\exists^{I_{\eta_{\beta_{\varepsilon,n-1}}\rest i}}
\gamma_{n-1})(\exists^\lambda\zeta)(\forall \ell<n)\\
(\eta_{\beta_{\varepsilon,\ell}}\rest (i+1)=(\eta_{\beta_{\varepsilon,\ell}}
\rest i)\hat{\ }\langle\gamma_\ell\rangle)\}\\
\end{array}
\]
is in the ideal $J^{\bd}_{\cf(\mu)}$. Take one such $\varepsilon$. Choosing
$i\in \cf(\mu)\setminus Z_\varepsilon$, $i>i^*$ such that $h(i)=n$ we may
follow exactly as in the last part of the proof of \ref{green} and we find
$\varepsilon_0,\varepsilon_1<\lambda$ such that for each $\ell<n$ 
\begin{quotation}
\noindent $\eta_{\beta_{\varepsilon_0,\ell}}\rest i=
\eta_{\beta_{\varepsilon_1,\ell}}\rest i$ but

\noindent $\ba_{\eta_{\beta_{\varepsilon_0,\ell}}\rest i}\models
y_{\eta_{\beta_{\varepsilon_0,\ell}}\rest (i+1)} \cap
y_{\eta_{\beta_{\varepsilon_1,\ell}}\rest (i+1)}=\bz$ 
\end{quotation}
which implies that
\[(\forall\ell<n)(\ba^*_{\xi_\ell}\models x^{\red}_{\beta_{\varepsilon_0,\ell}}
\cap x^{\red}_{\beta_{\varepsilon_1,\ell}}=\bz).\]
\QED

\begin{proposition}
\label{contind}
Assume that
\begin{description}
\item[(a)] $\langle\lambda_i: i<\delta\rangle$ is an increasing sequence of
regular cardinals such that $\delta<\lambda_0$,
$\lambda_i>\max\pcf\{\lambda_j: j<i\}$ (the last is our natural assumption);
\item[(b)] $\aleph_0<\theta=\cf(\theta)<\bigcup\limits_{i<\delta}\lambda_i$
(naturally we assume just $\cf(\theta)=\theta<\lambda_0$)
\item[(c)] $\lambda=\max\pcf\{\lambda_i: i<\delta\}$
\item[(d)] $(\theta,\lambda_i,\max\pcf\{\lambda_j:j<i\})\in\K_{\swmk}$
\item[(e)] for each $\tau\in\{\lambda\}\cup\bigcup\limits_{\alpha<\delta}
\pcf\{\lambda_i: i<\alpha\}$ we have 
\[\{\xi<\tau: \cf(\xi)=\theta\}\in I[\tau]\]
or at least 
\begin{quotation}
\noindent for some $\bar{f}^\tau=\langle f^\tau_\varepsilon: \varepsilon<
\tau\rangle$, $<_{J_{=\tau}}$--increasing cofinal in $\prod\limits_{i<\alpha}
\lambda_i/J_{=\tau}$ we have:
\[\gamma<\tau\ \&\ \cf(\gamma)=\theta\ \ \ \Rightarrow\ \ \
f^\tau_\gamma\mbox{ is good in }\bar{f}^\tau\]
\end{quotation}
(see \cite{Sh:345a}, \cite{Sh:355}, section 1 (see 1.6(1) and then
\cite{MgSh:204}), 
\item[(f)] $|\pcf\{\lambda_i:i<\delta\}|<\theta$ or at least for each
$\alpha<\delta$ we have $|\pcf\{\lambda_i: i<\alpha\}|<\theta$.
\end{description}
{\em Then} $(\theta,\lambda)\in\K_{\wmk}$. Moreover
$(\theta,\lambda,\chi)\in\K_{\swmk}$ provided there is an almost disjoint
family of size $\chi$ in $[\lambda]^\lambda$. We may get algebras $\ba^{\red}$,
$\ba^{\green}$ as in main constructions such that 
\[\ba^{\red}\models\theta\mbox{-Knaster},\quad\ba^{\green}\models\lambda
\mbox{-cc}\quad\mbox{ and }\ba^{\red}*\ba^{\green}\models\neg\lambda
\mbox{-cc.}\] 
\end{proposition}

\noindent{\bf Remark~\ref{contind}.A:}\ \ \ Continues also the proof of 3.5 of
\cite{Sh:g}. 
\medskip

\Proof The main difficulty of the proof will be to construct a hereditary
$\theta$--free $<_{J_{{<}\lambda}}$--increasing sequence $\bar{\eta}=\langle
\eta_\alpha:\alpha<\lambda\rangle \subseteq \prod\limits_{i<\delta}\lambda_i$.
This is done in the claim below. For the notation used there let us note that
if $\alpha\leq\delta$ is a limit ordinal, $\tau\in \pcf\{\lambda_i:
i<\alpha\}$ then $J_{=\tau}[\{\lambda_i: i<\alpha\}]=J^\alpha_\tau$ is the
ideal on $\alpha$ generated by 
\[J_{{<}\tau}[\{\lambda_i: i<\alpha\}]\cup \{\alpha\setminus {\frak b}_\tau
[\{\lambda_i: i<\alpha\}]\}.\]
So in particular $\tcf(\prod\limits_{\i<\alpha}\lambda_i/J^\alpha_\tau)=\tau$.

\begin{claim}
There exists a tree $T\subseteq \bigcup\limits_{i<\delta}\prod\limits_{j<i}
\lambda_j$ such that $\lim_\delta(T)$ is $\theta$--hereditary free (and
$<_{J_{{<}\lambda}}$--cofinal).\\
Moreover for each $\alpha<\delta$ the size of $T_\alpha$ is
$\leq\max\pcf\{\lambda_i: i<\alpha\}$. 
\end{claim}

\noindent{\em Why?}\ \ \ For a limit ordinal $\alpha\leq\delta$ and $\tau\in
\pcf\{\lambda_i: i\leq\alpha\}$ (if $\alpha=\delta$ then $\tau=\lambda$)
choose a $<_{J^\alpha_\tau}$--increasing sequence $\bar{f}^{\alpha,\tau}=
\langle f^{\alpha,\tau}_\zeta: \zeta<\tau\rangle\subseteq\prod
\limits_{i<\alpha} \lambda_i$ cofinal in $\prod\limits_{i<\alpha}\lambda_i
/J^\alpha_\tau$ and such that
\begin{description}
\item[$(\tilde{\otimes})$] {\em if} $\zeta<\tau$, $\cf(\zeta)=\theta$

{\em then} for some unbounded set $Y_\zeta\subseteq \zeta$ (for simplicity
consisting of successor ordinals) and a sequence $\bar{s}^\tau=\langle
s^\tau_\xi: \xi\in Y_\zeta\rangle\subseteq J^\alpha_\tau$ we have
\[[\xi_1,\xi_2\in Y_\zeta\ \&\ \xi_1<\xi_2\ \&\ i\in\alpha\setminus
(s^\tau_{\xi_1}\cup s^\tau_{\xi_2})]\quad \Rightarrow\quad
f^{\alpha,\tau}_{\xi_1}(i) < f^{\alpha,\tau}_{\xi_2}(i).\]
\end{description}
[Why we can demand $(\tilde{\otimes})$? If in the assumption {\bf (e)} the
first part is satisfied then we follow similarly to the proof of
\ref{getfree}, compare \cite{Sh:355}, 1.5A, 1.6, pp 51--52. If we are in the
case of ``at least'' then this is exactly the meaning of goodness.]
Further we may demand that the sequence $\bar{f}^{\alpha,\tau}$ is
$^b$continuous: 
\begin{description}
\item[$(\tilde{\oplus})$] {\em if} $|\delta|<\cf(\zeta)<\lambda_0$,
$\zeta<\tau$ 

{\em then} $f^{\alpha,\tau}_\zeta(i)=\min\{\bigcup\limits_{\xi\in C}
f^{\alpha,\tau}_\xi(i): C$ is a club of $\zeta\}$ 
\end{description}
[compare the proof of 3.4 of \cite{Sh:345a}, pp 25--26]. 

For a limit ordinal $\alpha<\delta$ we define
\[\begin{array}{ll}
T^0_\alpha=\{f{\in}\prod\limits_{i<\alpha}\lambda_i: &\mbox{\bf (a) }f=
\max\{f^{\alpha,\tau_\ell}_{\zeta_\ell}: \ell{<}n\}\mbox{ for some}\\
\ &\quad n<\omega, \tau_\ell\in\pcf\{\lambda_i: i<\alpha\},
\zeta_\ell<\tau_\ell\\ 
\ &\mbox{{\bf (b)} for every }\tau{\in}\pcf\{\lambda_i: i{<}\alpha\}\mbox{ if
} \tau=\lambda\mbox{ or }\alpha<\delta\mbox{ then}\\
\ &\mbox{there is } \zeta_f(\tau){<}\tau\mbox{ such that}\\
\ &\quad f^{\alpha,\tau}_{\zeta_f(\tau)}\leq f\ \&\ 
f^{\alpha,\tau}_{\zeta_f(\tau)} = f\mbox{ mod } J^\alpha_\tau\}.\\
 \end{array}\]
(Note that if $\alpha=\delta$ then there is only one value of $\tau_\ell,\tau$
which we consider here: $\lambda$.) Let
$T'\subseteq\bigcup\limits_{i\leq\delta}\prod_{j<i}\lambda_j$ be a tree such
that for $\gamma\leq\delta$: 
\[T_\gamma^\prime=\{f\in\prod_{i<\gamma}\lambda_i: f\rest\alpha\in
T^0_\alpha\mbox{ for each limit }\alpha\leq\gamma\}.\]
Let 
\[A=\{\zeta<\lambda: (\exists
f\in\prod_{i\leq\delta}\lambda_i)[f^{\delta,\lambda}_\zeta\leq f\ \&\
f^{\delta,\lambda}_\zeta=f \mbox{ mod}\; J^\alpha_\lambda\ \&\ (\forall i\leq
\delta)(f\rest i\in T^\prime_\alpha)]\}\]
and for each $\zeta\in A$ let $f^*_\zeta$ be a function witnessing it. Now let
$T\subseteq\bigcup\limits_{i\leq\delta}\prod_{j<i}\lambda_j$ be a tree such
that $T_\delta=\{f^*_\zeta: \zeta\in A\}$. 

\noindent By the definition, $T$ is a tree, but maybe it does not have enough
levels? Let $\chi$ be a large enough regular cardinal. Take an increasing
continuous sequence $\langle N_i: i\leq \theta\rangle$ of elementary submodels
of $({\cal H}(\chi),{\in},{<}^*)$ such that
\begin{quotation}
$|N_i|=\Upsilon=\theta+|\pcf\{\lambda_\alpha: \alpha<\delta\}|<\lambda_0$,

$\Upsilon+1\subseteq N_i\in N_{i+1}$,

all relevant things are in $N_0$.
\end{quotation}
Now we define $f^*\in\prod\limits_{\alpha<\delta}\lambda_\alpha$ by 
\[f^*(\alpha)=\sup(\bigcup_{i<\theta}N_i\cap\lambda_\alpha).\]
Similarly as in \cite{Sh:355}, pp 63--65, one proves that $f^*\rest\alpha\in
T^0_\alpha$ for each limit $\alpha\leq\delta$. Hence for some $\zeta$ we have
$f^*= f^{\delta,\lambda}_\zeta\ \mbox{mod}\; J^\delta_\lambda$ and thus
$\zeta\in A$. Consequently $A$ is unbounded in $\lambda$.

By induction on $\alpha\leq\delta$ we prove that 
\begin{description}
\item[$(\tilde{\circledcirc})$] {\em if} $f_\zeta\in T_\alpha$ (for
$\zeta<\theta$) are pairwise distinct 

{\em then} there are $Z\in [\theta]^\theta$ and $j<\alpha$ such that
\[(\forall \zeta_0,\zeta_1\in Z)(\zeta_0\neq\zeta_1\ \Rightarrow\
[f_{\zeta_0}\rest j=f_{\zeta_1}\rest j\ \&\ f_{\zeta_0}(j)\neq f_{\zeta_1}(j)]
).\] 
\end{description}
If $\alpha$ is a non-limit ordinal then this is trivial. So suppose that
$\alpha$ is limit, $\alpha<\delta$. Then for some
$\tau_{\zeta,\ell}\in\pcf\{\lambda_i: i<\alpha\}$,
$\xi_{\zeta,\ell}<\tau_{\zeta,\ell}$, $n_\zeta<\omega$ (for $\zeta<\theta$,
$\ell<n_\zeta$) we have  
\[f_\zeta=\max\{f^{\alpha,\tau_{\zeta,\ell}}_{\xi_{\zeta,\ell}}: \ell<
n_\zeta\}.\] 
As $\theta>|\pcf\{\lambda_\beta: \beta<\alpha\}|$ we may assume that
$n_\zeta=n^*$, $\tau_{\zeta,\ell}=\tau_\ell$ and for each $\ell<n^*$ the
sequence $\langle\xi_{\zeta,\ell}: \zeta<\theta\rangle$ is either constant or
strictly increasing. Now, the second case has to occur for some $\ell$ and we
may follow similarly to \ref{claimfree} and then apply the inductive
hypothesis. We are left with the case $\alpha=\delta$. So let
$f_\zeta=f^*_{\beta_\zeta}$ for $\zeta<\delta$ and we continue as before (with
$\lambda$ for $\tau_\ell$).

\noindent This ends the proof of the claim (note that the arguments showing
that all the $T^0_\alpha$ are not empty prove actually that the tree $T$ has
enough branches to satisfy our additional requirements). 
\medskip

Now let $T$ be a tree is in the claim above. Let
$\bar{\eta}=\langle\eta_\alpha: \alpha<\lambda\rangle\subseteq\lim_\delta(T)$
be the enumeration of $\{f^*_\zeta: \zeta\in A\}$ such that $\bar{\eta}$ is
$<_{J_{<\lambda}}$--increasing cofinal in
$\prod\limits_{i<\delta}\lambda_i/J_{<\lambda}$. By the assumption {\bf (d)} 
for each $\eta\in T$ we find a marked Boolean algebra
$(\ba_\eta,\bar{y}_\eta)$ such that for every $i<\delta$ the sequence $\langle
(\ba_\eta,\bar{y}_\eta): \eta\in T_i\rangle$ witnesses that
$(\theta,\lambda_i, |T_i|)\in\K_{\swmk}$. These parameters determine a
$(\delta,\mu,\lambda)$--constructor $\C$, so we have the respective Boolean
algebra $\ba^{\red}(\C)$ (and its counterpart $\ba^{\green}(\C)$). To show
that they have the required properties we follow exactly the proof that
$(\theta,\lambda,\chi)\in\K_{\swmk}$, so we will present this only. 

First note that by \ref{morekna} the algebra $\ba^{\red}(\C)$ has the
$\theta$--Knaster property. Now, let $\langle A_\zeta: \zeta<\chi\rangle
\subseteq [\lambda]^\lambda$ be such that
\[\zeta_1<\zeta_2<\chi\ \ \ \Rightarrow\ \ \ |A_{\zeta_1}\cap A_{\zeta_2}| <
\lambda.\]
Let $\bar{x}_\zeta=\langle x^{\red}_\xi: \xi\in A_\zeta\rangle$ and let
$\ba_\zeta$ be the subalgebra of $\ba^{\red}(\C)$ generated by
$\bar{x}_\zeta$. We want to show that the sequence $\langle
(\ba_\zeta,\bar{x}_\zeta): \zeta<\chi\rangle$ witnesses
$(\theta,\lambda,\chi)\in \K_{\swmk}$. For this suppose that $\zeta_0<\ldots<
\zeta_{n-1}<\chi$, $n<\omega$ and $\beta_{\varepsilon,\ell}\in A_{\zeta,\ell}$
are increasing with $\varepsilon$ (for $\varepsilon<\lambda$, $\ell<n$) and
without loss of generality with no repetition. We may assume that
\[(\forall\ell<n)(\forall\varepsilon<\lambda)(\beta_{\varepsilon,\ell}\notin
\bigcup_{m\neq\ell} A_{\zeta_m}).\]
Further we may assume that for some $i^*<\delta$ and pairwise distinct
$\eta_\ell\in T_{i^*}$ (for $\ell<n$) we have
\[(\forall\varepsilon<\lambda)(\forall\ell<n)(\eta_{\beta_{\varepsilon,\ell}}
\rest i^*=\eta_\ell).\]
Now we take $i\in [i^*,\delta)$ such that
\[(\forall\gamma<\lambda_i)(\exists^\lambda \varepsilon<\lambda)(\forall
\ell<n)(\eta_{\beta_{\varepsilon,\ell}}(i)>\gamma)\]
(remember that each $\langle\eta_{\beta_{\varepsilon,\ell}}:
\varepsilon<\lambda\rangle$ is $<_{J_{<\lambda}}$--cofinal). Since $|T_i|<
\lambda_i$ we can find $\nu_0,\ldots,\nu_{n-1}\in T_i$ such that
$\eta_\ell\trianglelefteq\nu_\ell$ and 
\[(\forall\gamma<\lambda_i)(\exists^\lambda \varepsilon<\lambda)(\forall
\ell<n)(\eta_{\beta_{\varepsilon,\ell}}\rest i=\nu_\ell\ \&\
\eta_{\beta_{\varepsilon,\ell}}(i)>\gamma).\]
Consequently we may choose a sequence $\langle\langle \gamma_{\xi,\ell}:\ell<
n\rangle: \xi<\lambda_i\rangle\subseteq\lambda_i$ such that
$\xi<\gamma_{\xi,\ell}$ and
\[(\forall\xi<\lambda_i)(\exists^\lambda\varepsilon<\lambda)(\forall \ell<n)
(\eta_{\beta_{\varepsilon,\ell}}\rest (i+1)=\nu_\ell\hat{\ }\langle
\gamma_{\xi,\ell}\rangle).\]
Now we use the choice of $(\ba_{\nu_\ell},\bar{y}_{\nu_\ell})$ (witnessing
$(\theta,\lambda_i,|T_i|)\in \K_{\swmk}$) and we find $\xi_1<\xi_2<\lambda_i$
such that
\[(\forall\ell<n)(\ba_{\nu_\ell}\models y^{\nu_\ell}_{\gamma_{\xi_1,\ell}}
\cap y^{\nu_\ell}_{\gamma_{\xi_2,\ell}}=\bz),\]
which allows us to find $\varepsilon_1<\varepsilon_2<\lambda$ such that for
each $\ell<n$ the intersection $x_{\beta_{\varepsilon_1,\ell}}\cap
x_{\beta_{\varepsilon_2,\ell}}$ is $\bz$. \QED

\begin{conclusion}
If $\langle\mu_i: i\leq\kappa\rangle$ is a strictly increasing continuous
sequence of strong limit singular cardinals such that  $\kappa<\mu_0$,
$2^{\mu_i}=\mu_i^+$, $\kappa<\theta=\cf(\theta)<\mu_0$ and
\[\{\alpha<\mu^+_i: \cf(\alpha)=\theta\}\in I[\mu_i^+]\]
{\em then} $(\theta,\mu_\kappa^+)\in\K_{\wmk}$ and we may construct the
respective Boolean algebras $\ba^{\red}$, $\ba^{\green}$. \QED
\end{conclusion}

\begin{proposition}
Suppose that we have Boolean algebras $\ba^{\red}$, $\ba^{\green}$ such that 
\begin{quotation}
\noindent $\ba^{\red}$ satisfies the $\theta$--Knaster condition

\noindent for each $n<\omega$ the free product $(\ba^{\green})^n$ satisfies the
$\lambda$--cc 

\noindent the free product $\ba^{\red}*\ba^{\green}$ fails the $\lambda$-cc.
\end{quotation}
Then $(\theta,\lambda,\chi)\in\K_{\swmk}$, where $\chi=\lambda^+$ (or even if
$\chi$ is such that there is an almost disjoint family ${\cal A}\subseteq
[\lambda]^\lambda$ of size $\chi$). 
\end{proposition}

\Proof We have $y_\alpha\in(\ba^{\red})^+$ and $z_\alpha\in(\ba^{\green})^+$
for $\alpha<\lambda$ such that if $\alpha<\beta<\lambda$ then
\[\mbox{either }\quad\ba^{\red} \models y_\alpha\cap y_\beta=\bz\quad\mbox{ or
}\quad\ba^{\green}\models z_\alpha\cap z_\beta=\bz.\]
Let $A_\zeta\in [\lambda]^\lambda$  (for $\zeta<\chi$) be pairwise
almost disjoint sets. We want to show that the sequence
\[\langle(\ba^{\red},\bar{y}\rest A_\zeta): \zeta<\chi\rangle\]
is a witness for $(\theta,\lambda,\chi)\in\K_{\swmk}$. So we are given
$\zeta_0<\zeta_1<\ldots<\zeta_{n-1}<\chi$ and sequences
$\langle\alpha_{\varepsilon,\ell}: \varepsilon<\lambda\rangle\subseteq
A_{\zeta_\ell}$. Then, for some $\varepsilon^*<\lambda$ we have
\[\varepsilon^*\leq\varepsilon<\lambda\ \ \ \Rightarrow\ \ \
\alpha_{\varepsilon,\ell}\notin \bigcup_{m\neq \ell} A_{\zeta,m}.\]
We should find $\varepsilon_1<\varepsilon_2$ such that for all $\ell<n$
\[\ba^{\red}\models y_{\alpha_{\varepsilon_1,\ell}}\cap
y_{\alpha_{\varepsilon_2,\ell}}=\bz.\]
For this it is enough to find $\varepsilon^*<\varepsilon_1<\varepsilon_2$ such
that for $\ell<n$ 
\[\ba^{\green}\models z_{\alpha_{\varepsilon_1,\ell}}\cap
z_{\alpha_{\varepsilon_2,\ell}}\neq\bz.\]
But this we easily get from the fact that the free product $(\ba^{\green})^n$
satisfies the $\lambda$-cc. \QED

\begin{comment}
The proofs that the algebra $\ba^{\green}$ satisfies the $\lambda$-cc (see
\ref{green}, \ref{contind}) give that actually for each $n<\omega$ the product
$(\ba^{\green})^n$ satisfies $\lambda$-cc. So it is reasonable to add it
(though not needed originally).  
\end{comment}

\begin{comment}
The ``$\bar{\eta}$ is (strong{-}) $J$-cofinal for $(T,\bar{\lambda},\bar{I})$''
has easy consequences for the existence of colourings.
\end{comment}

\begin{remark}
{\em 
For $\mu$ strong limit singular we may sometimes get a cofinal sequence of
length $\lambda\in (\mu,2^\mu]$ without $2^\mu=\mu^+$. By \cite{Sh:430},
section 5, 

\noindent{\em if}:
\begin{description}
\item[(a)] $I_i$ is a $\chi_i$--complete, $|I_i|=\tau_i$, $\chi_i$ regular
\item[(b)] $\chi_i\leq\tau_i\leq(\chi_i)^{+n^*}$, $n^*<\omega$
\item[(c)] $\tcf(\prod\limits_{i<\delta}(\chi_i)^{+\ell}/J)=\lambda$ for each
$\ell\leq n^*$
\end{description}
{\em then} 
\begin{description}
\item[$(\alpha)$] there is a cofinal sequence in
$\prod\limits_{i<\delta}({\cal P}(\lambda_i)/I_i)/J$, because
\item[$(\beta)$] it has the true cofinality.
\end{description}
So if for arbitrarily large $\chi$, $2^{\chi}=\chi^+$, $2^{\chi^+}=\chi^{++}$
then we have the ideal we want and maybe the $\pcf$ condition holds. so
combining this and \ref{proexa} below we get that there may be an example of
our kind not because of GCH reasons, but still requiring some cardinal
arithmetic assumptions.   
}
\end{remark}

\begin{proposition}
\label{proexa}
Suppose that $\langle\lambda_i: i<\delta\rangle$ is a strictly increasing
sequence of regular cardinals, $I_i$ is a $(\prod\limits_{j<i}
\lambda_j)^+$--complete ideal on $\lambda_i$ (so
$\prod\limits_{j<i}\lambda_j<\lambda_i$) and $(\ba_i,\bar{y}_i,I_i)$
is a $\lambda_i$--well marked Boolean algebra (for $i<\delta$). 
\begin{enumerate}
\item Assume that $\prod\limits_{i<\delta}(I_i,{\subseteq})/J$ has true
cofinality $\lambda$. {\em Then} there exists a $(\theta,{\rm
not}\lambda)$--Knaster marked Boolean algebra.
\item Suppose in addition that $h:\delta\longrightarrow\omega$ is a function
such that   
\[(\forall n<\omega)(h^{-1}[\{n\}]\in J^+)\]
and $I^{[h(i)]}_i$ are the product ideals on $(\lambda_i)^n$ (for
$i<\delta$): 
\[I^{[h(i)]}_i\stackrel{\rm def}{=}\{B{\subseteq}(\lambda_i)^n:
\neg(\exists^{I_i}\gamma_0)\ldots(\exists^{I_i}\gamma_{h(i)-1})(\langle
\gamma_\ell: \ell< h(i)\rangle\in B).\]
Assume that 
\[\lambda=\tcf(\prod_{i<\delta}(I^{[h(i)]}_i,{\subseteq})/J)\]
and that the $(\ba_i,\bar{y}_i,I_i)$ satisfy the following requirement:
\begin{description}
\item[$(\tilde{*})_{h(i)}$] if $B\subseteq (\dom(\bar{y}_i))^{h(i)}$ is such
that 
\[(\exists^{I_i}\gamma_0)\ldots(\exists^{I_i}\gamma_{h(i)})(\langle\gamma_\ell:
\ell\leq h(i)\rangle\in B)\]
then there are $\gamma_\ell^\prime,\gamma_\ell^{\prime\prime}<\lambda_i$ (for
$\ell\leq h(i)$) such that for each $\ell$
\[\ba_i\models y_{i,\gamma_\ell^\prime}\cap y_{i,\gamma_\ell^{\prime\prime}}=
\bz.\] 
\end{description}
{\em Then} we can conclude that
$((2^{|\delta|})^+,\lambda,\lambda^+)\in\K_{\swmk}$ and we have a pair of
algebras $(\ba^{\red},\ba^{\green})$ as in main theorem \ref{main}. 
\end{enumerate}
\end{proposition}

\Proof The main point here is that with our assumptions in hands we may
construct a sequence $\langle\eta_\alpha: \alpha<\lambda\rangle \subseteq
\prod\limits_{i<\delta}\lambda_i$ which is quite stronger $J$--cofinal: 
it satisfies the requirement of \ref{cofinal}(6){\bf (b)} weakened to the
demand that the set there is not in the dual filter $J^c$. Of course this is
still enough to carry our proofs and we may use such a sequence to build the
right examples. 

\noindent{\em 1).}\ \ \ Let $\langle\langle A^\alpha_i:i<\delta\rangle:
\alpha<\lambda\rangle$ witness the true cofinality. By induction on
$\alpha<\lambda$ choose $\gamma_\alpha<\lambda$ and $\eta_\alpha\in
\prod\limits_{i<\alpha}\lambda_i$ such that
\begin{quotation}
\noindent $\langle\{\eta_\beta(i)\}: i<\delta\rangle\in\prod\limits_{i<\delta}
I_i$,

\noindent if $\beta<\alpha$ then $\gamma_\beta<\gamma_\alpha$ and $(\forall^J
i)(\eta_\beta(i)\in A^{\gamma_\alpha}_i)$ and

\noindent $\eta_\alpha(i)\notin A^{\gamma_\alpha}_i$.
\end{quotation}
For $\alpha=0$ or $\alpha$ limit, first choose $\gamma_\alpha=
\sup\{\gamma_{\alpha_1}+1: \alpha_1<\alpha\}$ and then choose $\eta_\alpha(i)$
by induction on $i$.\\
For $\alpha=\alpha_1+1$ first note that 
\[\langle\{\eta_{\alpha_1}(i)\}:i<\delta\in\prod\limits_{i<\delta} I_i.\]
Hence for some $\gamma^0_\alpha<\lambda$ we have
\[(\forall^J i)(\eta_{\alpha_1}(i)\in A^{\gamma_\alpha}_i).\]
Let $\gamma_\alpha=\max\{\gamma_{\alpha_1},\gamma^0_\alpha\}$. Now choose
$\eta_\alpha(i)$ by induction on $i$.

As $I_i$ is $|T_i|^+$--complete, clearly $\langle\eta_\alpha:
\alpha<\lambda\rangle$ is $J$--cofinal for $(T,J,\bar{I})$ and
\ref{notknaster}, \ref{consknaster} give the conclusion. 

\noindent{\em 2).}\ \ \ The construction of $\bar{\eta}$ is in a sense similar
to the one in the proof of \ref{stronger}, but we use our cofinality
assumptions. We have a cofinal sequence in
$\prod_{i<\delta}(I^{[h(i)]}_i,{\subseteq})/J$: 
\[\langle\langle A^\alpha_i: i<\delta\rangle: \alpha<\lambda\rangle.\]
For each $A^\alpha_i$  we have ``Skolem functions'' $f^\alpha_{i,\ell}$ for
$\ell< h(i)$ (like in the proofs of \ref{claimgreen}, \ref{supimpfub}).

We define $\eta_\alpha$ by induction on $\alpha<\lambda$. 
In the exclusion list we put all substitutions by
$\eta_{\gamma_0}\rest i,\ldots,\eta_{\gamma_{\ell-1}}\rest i$ for
$\gamma_k<\alpha$ to $f^\alpha_{i,\ell}$: each time we obtain a 
set in the ideal $I_i$ and a member $\bar{A}$ of $\prod\limits_{i<\delta} I_i$
such that if $(\forall^J i)(\eta(i)\notin A_i)$, $\eta\in\prod
\limits_{i<\delta} \lambda_i$ then $\eta$ satisfies the demand. Eventually we
have $|\alpha|^{<\omega}$ such elements of $\prod\limits_{i<\delta} I_i$. Let
them be $\{\bar{B}^{\alpha,\xi}: \xi\leq |\alpha|+\aleph_0\}$. Then for some
$\gamma_\alpha$ 
\[(\forall \xi<|\alpha|+\aleph_0)(\forall^J i<\delta)(B^{\alpha,\xi}_i
\subseteq A^{\gamma_\alpha}_i)\]
and similarly 
\[(\forall \beta<\alpha)(\forall^J i<\delta)(\eta_\beta(i)\in
A^{\alpha_i}_i).\] 
Choose $\eta_\alpha\in\prod\limits_{i<\delta}(\lambda_i\setminus
A^{\gamma_\alpha}_i)$. \QED

\begin{remark}
{\em
One of the main tools used in this section are (variants of) the following
observation:
\begin{quotation}
{\em if} $(\ba,\bar{y})$ is a $\lambda$-marked Boolean algebra such that $\ba$
is $\theta$--Knaster and if $\varepsilon(\alpha,\ell)<\lambda$ (for
$\alpha<\lambda$, $\ell<n$) are pairwise distinct then for some
$\alpha<\beta<\lambda$, for each $\ell<n$ we have $\ba\models
y_{\varepsilon(\alpha,\ell)}\cap y_{\varepsilon(\beta,\ell)}=\bz$

{\em then} $(\theta,\lambda,\lambda^+)\in\K_{\swmk}$. 
\end{quotation}
}
\end{remark}

\begin{conrem}
{\em
If $\mu$ is a strong limit singular cardinal,
$\cf(\mu)\leq\theta=\cf(\theta)<\mu$ 
then, by the methods of \cite{Sh:576}, one may get consistency of 
\begin{quotation}
if an algebra $\ba$ satisfies the $\theta$-ccc

then it satisfies the $\mu^+$-Knaster condition.
\end{quotation}
}
\end{conrem}

One may formulate the following question now:

\begin{pms}
\label{pms}
Suppose that $\ba$ is a Boolean algebra satisfying the $\theta$-cc and
$\lambda$ is a regular cardinal between $\mu^+$ and $(2^\mu)^+$.

Does $\ba$ satisfy the $\lambda$-Knaster condition?
\end{pms}
There a reasonable amount of information on consistency of the negative answer
in the next section, though \ref{pms} is not fully answered there. But a real
problem is the following.

\begin{problem}
Assume $\lambda=\mu^+$, $\cf(\mu)=\theta$ and $\mu$ is a strong limit
cardinal. Suppose that an algebra $\ba_0$ satisfies the $\lambda$-cc and an
algebra $\ba_1$ satisfies the $\theta^+$-cc.

Does the free product $\ba_0*\ba_1$ satisfy the $\lambda$-cc? (is this
consistent? see \ref{moreconclusion}).
\end{problem}

\begin{problem}
Is it consistent that 
\begin{quotation}
\noindent each Boolean algebra with the $\aleph_1$--Knaster property has the
$\lambda$--Knaster property for every regular (uncountable) cardinal
$\lambda$? 
\end{quotation}
\end{problem}

\section{Some consistency results}

We had seen that without inner models with large cardinals we have a complete
picture, e.g.: 
\begin{description}
\item[$(\aleph)$] if $\theta=\cf(\theta)>\aleph_0$, $\ba$ is a Boolean algebra
satisfying the $\theta$--cc and $\lambda$ is a regular cardinal such that
\[(\forall\tau<\lambda)(\tau^{{<}\theta}<\theta)\]
{\em then} the algebra $\ba$ satisfies the $\lambda$-Knaster condition.
\item[$(\beth)$] if $\theta=\cf(\theta)>\aleph_0$, $\theta<\mu=\mu^{{<}\mu}<
\lambda=\cf(\lambda)<\chi=\chi^{\lambda}$

\noindent{\em then} there is a $\mu^+$-cc $\mu$-complete forcing notion $\p$ of
size $\chi$ such that
\[\forces_{\p}\mbox{``the }\theta\mbox{--cc implies the }\lambda\mbox{--Knaster
property''.}\] 
Moreover
\item[$(\beth)^+$] if $\mu=\mu^{<\theta}<\lambda=\cf(\lambda)\leq 2^{\mu}$
then the $\theta$--cc implies the $\lambda$--Knaster property.
\item[$(\gimel)$] if $\theta=\cf(\theta)<\mu$, $\mu$ is a strong limit
singular cardinal, $\cf(\mu)=\theta$

\noindent{\em then} the $\theta^+$--cc {\em does not} imply the
$\mu^+$--Knaster property (and even we have the product example).
\end{description}
In $(\gimel)$, if we allow $(2^\theta)$-cc we may get even better conclusion.
In this section we want to show, under a large cardinals hypothesis, the
consistency of failure. 

\begin{proposition}
Assume that $\kappa$ is a supercompact cardinal,
$\kappa<\lambda=\cf(\lambda)$. Let $\ba$ be a Boolean algebra which does not
have the $\lambda$--Knaster property. {\em Then}
\[(\exists\theta)(\aleph_0<\theta=\cf(\theta)<\kappa\ \&\ \ba\mbox{ does not
have the }\theta\mbox{--Knaster property}).\]
\end{proposition}

\Proof Since $\kappa$ is supercompact, for every second order formula $\psi$:
\begin{quotation}
\noindent if $M\models\psi$

\noindent then for some $N\prec M$, $|N|<\kappa$, $N\models \psi$
\end{quotation}
(see Kanamori and Magidor, \cite{KnMg78}). \QED

\begin{proposition}
\begin{enumerate}
\item If $\aleph_0<\lambda_0<\lambda_1$ are regular cardinals such that
\begin{description}
\item[$(*)_{\lambda_0,\lambda_1}$] for every $x\in \cH(\lambda_1^+)$ there is
$N\prec(\cH(\lambda_1^+), {\in})$ such that $x\in N$ and $N\cong
(\cH(\lambda_0^+), {\in})$
\end{description}
{\em then} if a Boolean algebra $\ba$ has the $\lambda_0$--Knaster property
then it has $\lambda_1$--Knaster property (and $\ba\models\lambda_0$--cc
implies $\ba\models\lambda_1$--cc). 
\item The condition $(*)_{\lambda_0,\lambda_1}$ above holds if for some
$\kappa_0,\kappa_1$, $\kappa_0<\lambda_0$, $\kappa_1<\lambda_1$ we have:
\begin{description}
\item[$(\oplus)$] there is an elementary embedding $j:\bV\longrightarrow M$
with the critical point $\kappa_0$ and such that $j(\kappa_0)=\kappa_1$,
$j(\lambda_0)=\lambda_1$ and $M^{\lambda_1}\subseteq M$.
\end{description}
\item If $\kappa_0$ is a $2$-huge cardinal (or actually less) and e.g.
$\lambda_0=\kappa_0^{+\omega+1}$ then for some
$\lambda_1=\kappa_1^{+\omega+1}$ the condition $(\oplus)$ above holds (we can
assume GCH).
\end{enumerate}
\end{proposition}

\Proof Just check. \QED

\begin{proposition}
\label{consist}
Assume that
\[\bV\models\mbox{`` GCH}+\mbox{ there is 2-huge cardinal}>\theta=\cf(\theta)
\mbox{''}\]
(can think of $\theta=\aleph_0$). Then there is a $\theta$--complete forcing
notion $\p$ such that in $\bV^{\p}$:
\begin{description}
\item[(a)] GCH holds
\item[(b)] if a Boolean algebra $\ba$ has the $\theta^+$--Knaster property
then it has the $\theta^{+\theta+1}$--Knaster property

\noindent (note that if $\aleph_\theta>\theta$ then
$\theta^{+\theta+1}=\aleph_{\theta+1}$). 
\end{description}
\end{proposition}

\Proof Similar to \cite{LMSh:198}. \QED
\medskip

Chasing arrows what we use is

\begin{proposition}
If $\bV\models$GCH (for simplicity), $\theta=\cf(\theta)=\cf(\mu)<\mu$, a
Boolean algebra $\ba$ does not satisfy the $\mu^+$--Knaster condition and
$\q=\Levy(\theta,\mu)$ 

\noindent then $\bV^{\q}\models$``$\ba$ does not have the $\theta^+$--Knaster
property''. \QED
\end{proposition}

\section{More on getting the Knaster property}

Our aim here is to get a $\ZFC$ result (under reasonable cardinal arithmetic
assumptions) which implies that our looking for $(\kappa,{\rm
not}\lambda)$--Knaster marked Boolean algebras near strong limit singular is
natural. Bellow we discuss the relevant background. The proof relays on $\pcf$
theory (but only by quoting a simply stated theorem) and seems to be a good
example of the applicability of $\pcf$.

\begin{theorem}
\label{th8.1}
Assume $\mu= \mu^{<\beth_\omega}$.
\begin{enumerate}
\item If a Boolean algebra $\ba$ of cardinality $\leq 2^\mu$ satisfies the
$\aleph_1$--cc then $\ba$ is $\mu$-linked (see below).

\item If $\ba$ is a Boolean algebra satisfying the $\aleph_1$--cc then $\ba$
has the $\lambda$-Knaster property {\em for every} regular cardinal
$\lambda\in (\mu,2^\mu].$ 
\end{enumerate}
\end{theorem}

Where

\begin{definition}
\label{def8.2}
\begin{enumerate}
\item A Boolean algebra $\ba$ is $\mu$-linked if $\ba\setminus\{\bz\}$ is
the union of $\leq \mu$ sets of pairwise compatible elements.

\item A Boolean algebra $\ba$ is $\mu$-centered if $\ba\setminus \{\bz\}$
is the union of $\leq \mu$ filters.
\end{enumerate}
\end{definition}

Of course we can replace the $\aleph_1$--cc, $\beth_\omega$ by the
$\kappa$--cc, $\beth_\omega(\kappa)$ (see more later). The proof is self
contained except relayence on a theorem quoted from \cite{Sh:460}.

Let us review some background. By \cite{Sh:92}, 3.1, if $\ba$ is a
$\kappa$-cc Boolean algebra of cardinality $\mu^+$ and $\mu=\mu^{<\kappa}$
then $\ba$ is $\mu$-centered. The proof did not work for $\ba$ of cardinality
$\mu^{++}$ even if $2^\mu\geq \mu^{++}$ by \cite{Sh:126}, point being we
consider three elements. But if $\mu=\mu^{<\mu} < \lambda^{<\lambda}$, for
some $\mu^+$--cc $\mu$-complete forcing notion $\p$ of cardinality $\lambda$,
in $\bV^{\p}$:
\begin{quotation}
\noindent if $\ba$ is a $\mu-\cc$ Boolean algebra of cardinality $<\lambda$
the $\ba$ is $\mu$-centered 
\end{quotation}
(follows from an appropriate axiom). Juhasz, Hajnal and Szentmiklossy
\cite{HaJuSz} continue this restricting themselves to $\mu$-linked. Then proof
can be carried for $\mu^{++}$, and they continue by induction. However as in
not few cases, the problem was for $\lambda^+$, when $\cf(\lambda)=\aleph_0$
so they assume 
\begin{description}
\item{$\otimes$}\ \ \ \ \ if $\lambda\in (\mu, 2^\mu)$,
$\cf(\lambda)=\aleph_0$ then $\lambda=\lambda^{\aleph_0}$ and
$\square_\lambda$  
\end{description}
(on the square see Jensen \cite{Jn}). This implies that if we start with
$\bV={\bf L}$ and force, then the assumption ($\otimes$) holds, so it is a
reasonable assumption. Also they prove the consistency of the failure of the
conclusion when $\otimes$ fails relaying on \cite{HJSh:249} (on a set system +
graph constructed there) and on colouring of graphs (section 2 of
\cite{HaJuSz}), possibly $2^{\aleph_0}=\aleph_1$,
$2^{\aleph_1}=\aleph_{\omega+1}$, $|\ba|=2^{\aleph_1}$, $\ba$ satisfies the
$\aleph_1$--cc but is not $\aleph_1$-linked, only $\aleph_2$-linked.

This gives the impression of essentially closing the issue, and so I
would have certainly thought some years ago, but this is not the case,
examplifying the danger of looking at specific cases. In fact, as we shall
note in the end, their consistency result is best possible under our knowledge
of relevant forcing methods. They use \cite{HJSh:249} to have ``many very
disjoint sets''(i.e. $\langle X_\alpha: \alpha\in S\rangle$, $S\subseteq
\{\delta< \aleph_{\omega+1}: \cf(\delta)=\aleph_1\}$, $X_\alpha \subseteq
\alpha=\sup X_\alpha$, and $\alpha\neq \beta \Rightarrow X_\alpha\cap X_\beta$
finite). 

On $\pcf$ see \cite{Sh:g}. Now, \cite{Sh:460} has half jokingly a strong
claim of proving $\GCH$ under reasonable reinterpretation. In particular
\cite{Sh:460} says there cannot be many strongly almost disjoint quite
large sets, so this blocks reasonable extensions of \cite{HaJuSz}. Now the
main theorem of \cite{Sh:460} enables us to carry the induction on
$\lambda\in (\mu, 2^{\mu}]$ as in \cite{Sh:92}, 3.1, \cite{HaJuSz}, 3.x.

\begin{proposition}
\label{claim8.3}
Suppose that:
\begin{description}
\item[(a)] $\lambda>\theta = \cf(\theta)\geq \kappa =\cf(\kappa)>
\aleph_0$
\item[(b)] there are a club $E$ of $\lambda$ and a sequence $\bar {\P}
=\langle {\P}_\alpha: \alpha\in E\rangle$ (with $\alpha\in E \Rightarrow
|\alpha| \mid \alpha$) such that
\begin{description}
\item[(i)] ${\P}_\alpha \subseteq [\alpha]^{<\kappa}$,
$|{\P}_\alpha|\leq |\alpha|$ and $\bar {\P}$ is increasing continuous,
\item[(ii)] if $X\subseteq \lambda$ has order type $\theta$, then for
some increasing $\langle \gamma_\varepsilon: \varepsilon< \kappa\rangle$
we have $\gamma_\varepsilon \in X$ and for each $\varepsilon< \kappa$
for some $\zeta\in (\varepsilon, \kappa)$ and $\alpha\leq
\min(E\setminus \gamma_\varepsilon)$ we have $\{ \gamma_\zeta:
\zeta<\varepsilon\} \in {\P}_\alpha$,
\end{description}
\item[(c)] $\ba$ is a Boolean algebra satisfying the $\kappa$--cc,
$|\ba|=\lambda$.
\end{description}
{\em Then} we can find a Boolean algebra ${\ba}'$ and a sequence
$\langle{\ba}'_\alpha: \alpha\in E\rangle$ of subalgebras of $\ba'$ such that
\begin{description}
\item[$(\alpha)$] $\ba\subseteq {\ba}' \subseteq {\ba}^{\com}$ (the
completion)
\item[$(\beta)$] ${\ba}'=\bigcup\limits_{\alpha\in E} {\ba}'_\alpha$,
$|{\ba}'_\alpha|\leq |\alpha|\neq \aleph_0$, $\langle\ba'_\alpha:\alpha\in
E\rangle$ increasing continuous in $\alpha$
\item[$(\gamma)$] if $\alpha\in E$, $x\in {\ba}'\setminus \{\bz\}$
then for some $Y\subseteq {\ba}'_\alpha\setminus \{\bz\}$, $|Y|<\theta$ we
have
\begin{description}
\item[\ ] if $y\in Y$ then $y\cap x=0_{{\ba}'}$ and
\item[\ ] if $z\in {\ba}'_\alpha$ is such that $z\cap x= 0_{{\ba}'}$ then
$z\leq \sup Y'\in {\ba}'_\alpha$ for some $Y'\in [Y]^{<\kappa}$
\end{description}
\item[$(\delta)$] if either $(*)_1$ or $(*)_2$ (see below) holds then we can
add 
\begin{description}
\item[\ ] $Y$ generates the ideal $\{z\in {\ba}'_\alpha: z\cap
x=\bz_{{\ba}'}\}$
\end{description}
where
\begin{description}
\item[$(*)_1$] $(\forall \varepsilon<\theta)(|\varepsilon|^{<\kappa}<\theta)$
\item[$(*)_2$] in {\bf (b)} we can add:

if $X\subseteq \alpha$, $|X|<|\alpha|$ then for some $\tau$, $\tau^{<\kappa}
<\theta$ and $h: X \rightarrow \tau$ we have: if $Y\subseteq X$,
$h\restriction Y$ is constant then $Y\in {\P}_\alpha$. 
\end{description}
\end{description}
\end{proposition}

\Proof Let $\chi$ be a large enough regular cardinal. Let
$\ba=\{x_\varepsilon: \varepsilon<\lambda\}$, let ${\ba}^{\com}$ be the
completion of $\ba$. We choose by induction on $\alpha\in E$ an elementary
submodel $N_\alpha$ of $(H(\chi), \in, <^*_\chi)$ of cardinality $|\alpha|$,
increasing continuous in $\alpha$, such that $\ba$, $\langle x_\varepsilon:
\varepsilon< \lambda\rangle$, ${\ba}^\com$, $\bar {\P}$, $\lambda$, $\theta$,
$\kappa$ belong to $N_0$ and $\langle N_\zeta: \zeta\leq
\varepsilon\rangle \in N_{\varepsilon+1}$.

Note: if $\alpha\in \nacc (E)$ then $\alpha\in N_\alpha$ hence
${\P}_\alpha\subseteq N_\alpha$.

Let 
\[{\ba}'_\alpha\stackrel{\rm def}{=} N_\alpha \cap {\ba}^\com,\qquad
{\ba}'=\bigcup\limits_{\alpha\in E}{\ba}'_\alpha.\] 
We define by induction on $\alpha\in E$ a one-to-one function $g_\alpha$ from
${\ba}'_\alpha$ onto $\alpha$ such that
\[\beta\in \alpha\cap E\ \ \Rightarrow\ \ g_\beta \subseteq g_\alpha,\mbox{
and } g_\alpha\mbox{ is the }<^*_\chi\mbox{-first such }g,\] 
so $g_\alpha\in N_{\min(E\setminus (\alpha+1))}$. Let $g=
\bigcup\limits_{\alpha\in E} g_\alpha$. Thus $g$ is a one-to-one function 
from ${\ba}'$ onto $\lambda$. In the conclusion clauses $(\alpha)$, $(\beta)$
should be clear and let us prove clause $(\gamma)$. So let $\alpha\in E$,
$x\in {\ba}'\setminus \{\bz\}$. We define $J=\{z\in {\ba}'_\alpha:
{\ba}'\models$ ``$z\cap x = \bz$''$\}$. Then $J$ is an ideal of
${\ba}'_\alpha$. We now try to choose by induction on $\varepsilon<\theta$,
elements $y_\varepsilon\in J$ such that
\begin{description}
\item[(i)] $y_\varepsilon$ is a member of $J\setminus \{\bz_{\ba}\}$
\item[(ii)] there is no $u\in [\varepsilon]^{<\kappa}$ such that $y_\alpha\leq
\sup\limits_{\zeta\in u} y_\zeta \in {\ba}'_\alpha$ ($\sup$ - in the
complete Boolean algebra ${\ba}^{\com}$)
\item[(iii)] under {\bf (i)} + {\bf (ii)},  $g(y_\varepsilon)$ ($<\lambda$) is
minimal (hence under {\bf (i)} + {\bf (ii)}, $\beta_\varepsilon\stackrel{\rm
def}{=}\min\{ \beta\leq \alpha: y_\varepsilon \in {\ba}'_\beta\}$ is minimal).
\end{description}
If we are stuck for some $\varepsilon<\theta$, then for every $y\in J$ the
condition {\bf (ii)} fails (note that {\bf (iii)} does not change at this
point) i.e. there is a respective set $u$. So suppose $y_\varepsilon$
is defined for $\varepsilon<\theta$. Clearly 
\[\zeta<\varepsilon\ \ \Rightarrow\ \ g(y_\zeta)< g(y_\varepsilon)\]
and hence $\zeta<\varepsilon<\theta \Rightarrow \beta_\zeta \leq
\beta_\varepsilon$, and $\zeta<\varepsilon \Rightarrow y_\zeta \neq
y_\varepsilon$. Now apply clause {\bf (b)(ii)} of the assumption to the set
$X=\{\gamma_\varepsilon: \varepsilon< \theta\}$ to get a contradiction. \QED

\begin{proposition}
\label{claim8.4}
Suppose that 
\begin{description}
\item[(a)] $\lambda>\theta=\cf(\theta)\geq\kappa=\cf(\kappa)>\aleph_0$, and
$\mu=\mu^{<\theta}\leq \lambda\leq 2^{\mu}$ and 

$(\forall \lambda<\theta)[|\lambda|^{<\kappa}<\theta]$,
\item[(b)] as in \ref{claim8.3} and either $(*)_1$ or $(*)_2$ of clause
$(\delta)$ of \ref{claim8.3},
\item[(c)] $\ba$ is a $\kappa$-cc Boolean algebra of cardinality $\lambda$,
\item[(d)] every subalgebra ${\ba}'\subseteq {\ba}^\com$ of cardinality
$<\lambda$ is $\mu$-linked (see definition \ref{def8.2}(1)).
\end{description}
{\em Then} $\ba$ is $\mu$-linked.
\end{proposition}

\Proof Let $\langle {\ba}'_\alpha: \alpha \in E\rangle$ was as in the
conclusion of \ref{claim8.3}. Without loss of generality we may assume that
the set of elements ${\ba}'_\alpha$ is $\alpha$. Let for $\alpha\in E$,
$h_\alpha: {\ba}'_\alpha \setminus \{\bz\} \rightarrow  \mu$ be such that:
\[h_\alpha(x_1) = h_\alpha(x_2)\ \ \Rightarrow\ \ x_1 \cap x_2 \neq
\bz_{\ba}.\]
For each $x\in {\ba}'\setminus\ba'_{\min(E)}$ let $\alpha(x)=\max\{\alpha\in
E: x\notin {\ba}'_\alpha\}$ (well defined as ${\ba}'=\bigcup\limits_{\alpha\in
E}{\ba}'_\alpha$ and $\langle {\ba}'_\alpha: \alpha\in E\rangle$ is increasing
continuous), and let $Y_{x, \alpha}\subseteq {\ba}'_\alpha$ be such that
$|Y_{x,\alpha}|<\theta$ and
\[Y_x\subseteq J_x\stackrel{\rm def}{=}\{y\in {\ba}'_\alpha: y\cap
x=\bz_{\ba}\}\qquad \mbox{ and}\]
$Y_x$ is cofinal in $J_{x}$ ($Y_x$ exists by \ref{claim8.3}, see clause
$(\delta)$). 

Let us define $u^0_x=\{ 0, \alpha(x)\}$ and $Y^0_x$ the subalgebra of
${\ba}'$ generated by $\{x\}$, and $u^{n+1}_x= u^n_x \cup \{\alpha(y):
y\in Y^n_x\}$ and $Y^{n+1}_x$ be the subalgebra of ${\ba}'$ generated by
\[Y^n_x\cup\bigcup\{Y_{x_1, \alpha}: x_1\in Y^n_x \mbox{ and }\alpha\in 
u^n_x\}.\]
Finally let $Y^\omega_x= \bigcup\limits_{n<\omega} Y^n_x$. As $\theta$ is
regular, $|Y^n_x|<\theta$ and as in addition $\theta$ is uncountable,
$|Y^\omega_x|<\theta$. Let $u_x=\{\alpha(y): y\in Y^\omega_x\}$. We can find
$A_\zeta\subseteq {\ba}'\setminus \{\bz\}$ for $\zeta< \mu$ such that
${\ba}'\setminus \{\bz\}= \bigcup\limits_{\zeta< \mu} A_\zeta$ and
\begin{description}
\item[$(\tilde{\circledast})$] if $x_1$, $x_2\in A_\zeta$, then there are
one-to-one functions $f:Y^\omega_{x_1}\stackrel{\rm onto}{\longrightarrow}
Y^\omega_{x_2}$ and $g:u_{x_1}\stackrel{\rm onto}{\longrightarrow}u_{x_2}$
such that: 
\begin{description}
\item[(i)] $f$, $g$ preserve the order,
\item[(ii)] $f(x_1)=x_2$ and if $y\in Y^\omega_{x_1}$ then
$g(\alpha(y))=\alpha(f(y))$,  
\item[(iii)] if $\alpha\in u_{x_1}$, $y\in {\ba}'_\alpha\cap
Y^\omega_{x_1}$ then $h_\alpha(x_1)= h_{g(\alpha)}(f(x_1))$
\item[(iv)] $f$ is an isomorphism (of Boolean algebras)
\item[(v)] $g$ is the identity on $u_{x_1}\cap u_{x_2}$
\item[(vii)] $f$ is the identity on $Y^\omega_{x_1}\cap Y^\omega_{x_2}$
\end{description}
\end{description}
(Why? By \cite{EK} or use $\langle\eta_x: x\in {\ba}'\rangle$, $\eta_x\in
{}^{\mu}2$ with no repetitions.)

So it is enough to prove:
\[x_1, x_2\in A_\zeta \Rightarrow x_1\cap x_2 \neq 0_{\ba}.\]
Let $D_1$ be an ultrafilter of $Y^\omega_{x_1}$ to which $x_1$ belongs,
$D_2=: \{ f(y): y\in Y^\omega_{x_2}\}$ (an ultrafilter on
$Y^\omega_{x_2}$ to which $x_2$ belongs). It suffices to prove that for each
$\alpha\in E$, $D_1\cap {\ba}_\alpha$, $D_2\cap {\ba}_\alpha$ generate non
trivial filters on ${\ba}_\alpha$. We do it by induction on $\alpha$ (note if
$\alpha\leq \beta$ this holds for $\alpha$ provided it holds for $\beta$).
If $\alpha\in u_{x_1}\cap u_{x_2}$ use clause {\bf (iii)} of
$(\tilde{\circledast})$ -- note that this includes the case when $\alpha=0$.
For $\alpha\in \acc(E)$ it follows by the finiteness of the condition. In the
remaining case $\beta=\sup(E\cap\alpha)<\alpha$ and if $Y^\omega_{x_1}\cap
{\ba}'_\alpha\subseteq{\ba}'_\beta$, $Y^\omega_{x_2}\cap {\ba}'_\alpha
\subseteq {\ba}'_\beta$ this is trivial. So by symmetry we may assume that
$\alpha\in u_{x_1}\setminus u_{x_2}$ and use the definition of $Y_y$ for $y\in
B_\alpha\cap Y^\omega_{x_1}\setminus {\ba}_\beta$. \QED

\begin{proposition}
\label{claim8.5}
Assume $\mu=\mu^{<\beth_\omega(\kappa)}$. Then for every $\lambda\in (\mu,
2^\mu]$ of cardinality $>\mu$, for every large enough regular $\theta<
\beth_\omega(\kappa)$ clause (b) of \ref{claim8.3} holds.
\end{proposition}

\Proof By \cite{Sh:460}, for every $\tau\in [\mu, \lambda)$ for some
$\theta_\tau< \beth_\omega(\kappa)$, we have:
\begin{description}
\item[$(\tilde{\circleddash})$] there is $\P={\P}_\tau\subseteq [\tau]^{<
\beth_\omega(\kappa)}$ closed under subsets such that $|\P| \leq \tau$
and every $X\in [\tau]^{<\beth_\omega(\kappa)}$ is the union of $<\theta_\tau$
members of members of ${\P}_\tau$. 
\end{description}
Now as $\cf(\lambda)>\mu$ for some $n<\omega$, the set
\[\Theta=\{\tau: \mu<\tau<\lambda, \theta_\tau\leq \beth_\omega(\kappa)\}\]
is an unbounded subset of $\Card\cap (\mu, \lambda)$. Let $\theta\geq
(\beth_{n+1}(\kappa))$ be regular. Now choose $E$ club of $\lambda$ such
that $\alpha\in \nacc(E) \Rightarrow |\alpha|\in \Theta$ and choose
${\P}_\alpha\subseteq [\alpha]^{<\kappa}$ increasing continuous with
$\alpha\in E$, such that for $\alpha\in \nacc(E)$, for every $X\in
[\alpha]^\theta$, for some $h: X \longrightarrow \beth_n(\kappa)$, if
$Y\subseteq X$, $|Y|<\kappa$ and $h\restriction Y$ constant then $Y\in
{\P}_\alpha$.

Now suppose $X\subseteq \lambda$, $\otp(X)=\theta$, so let
$X=\{\gamma_\varepsilon: \varepsilon<\theta\}$, $\beta_\varepsilon$ increasing
with $\varepsilon$; let $\beta_\varepsilon=\min\{\alpha\in E:
\gamma_\varepsilon< \beta\}$, so $\zeta<\varepsilon \Rightarrow
\beta_\zeta\leq \beta_\varepsilon$ and $\beta_\varepsilon \in \nacc(E)$
and there is $h_\varepsilon: \{\zeta: \zeta<\varepsilon\} \rightarrow
\beth_n(\kappa)$ such that for every $j< \beth_n(\kappa)$,
\[u\in [\varepsilon]^{<\kappa}\ \&\ (h\restriction u\mbox{ constant})\ \
\Rightarrow\ \ \{ \gamma_\zeta: \zeta\in u\}\in {\P}_{\beta_\varepsilon}.\]
Applying the Erd\"os--Rado theorem (i.e. $\theta \rightarrow
(\beth_n(\kappa)^+)^2_{\beth_n(\kappa)}$) we get the desired result (the proof
is an overkill). \QED

\begin{mainconc}
\label{conc8.6}
Suppose that $\kappa$ is a regular uncountable cardinal,
$\mu=\mu^{\beth_\omega(\kappa)}$ and $\ba$ is a Boolean algebra satisfying 
the $\kappa$-cc.
\begin{enumerate}
\item If $|{\ba}|\leq 2^\mu$ then $\ba$ is $\mu$-linked
\item If $\lambda$ is regular $\in (\mu, 2^\mu]$ {\em then} $\ba$
satisfies the $\lambda$-Knaster condition.
\end{enumerate}
\end{mainconc}

\Proof 1) We prove this by induction on $\lambda=|\ba|$. If $|\ba|\leq \mu$
this is trivial and if $\cf(|\ba|)\leq \mu$ this follows easily by
the induction hypothesis. In other cases by \ref{claim8.5}, for some
$\theta^*< \beth_\omega(\kappa)$ for every regular $\theta\in (\theta^*,
\beth_\omega(\kappa))$, clause (b) of \ref{claim8.3} holds. Choose
$\theta=(\theta^{\aleph_0})^{++}$, so for this $\theta$ both clause (b)
of \ref{claim8.3} and $(*)_1$ of clause $(\delta)$ of \ref{claim8.3}
hold. Thus by claim \ref{claim8.4} we can prove the desired conclusion for
$\lambda=|\ba|$.
\medskip

\noindent 2) Follows from part 1). \QED

\begin{proposition}
\label{claim8.7}
\begin{enumerate}
\item In \ref{conc8.6} we can replace the assumption
$\mu=\mu^{\beth_\omega(\kappa)}$ by $\mu=\mu^{<\tau}$ if
\begin{description}
\item[$\otimes$] for every $\lambda\in (\mu, 2^\mu)$ of cardinality
$>\mu$ for some $\theta=\cf(\theta)\geq \kappa$ we have: clause (b) of
\ref{claim8.3} and $(*)_2$ of clause $(\delta)$ of \ref{claim8.3} hold.
\end{description}
\item If $\lambda^*\in (\mu, 2^\lambda)$ and we want to have the
conclusion of \ref{conc8.6}(1) with $|\ba|=\lambda^*$ and
\ref{conc8.6}(2) for $\lambda^*$-Knaster only {\em then} it suffice to
restrict ourselves in $\otimes$ to $\lambda\leq \lambda^*$. \QED
\end{enumerate}
\end{proposition}

\begin{proposition}
\label{fact8.8}
In \ref{claim8.3}, if $(\forall \varepsilon<
\theta)[|\varepsilon|^{<\kappa}<\theta]$ then we can weaken clause (ii)
of assumption (b) to
\begin{description}
\item[(ii)$'$] if $X\subseteq \lambda$ has order type $\theta$ then for
some $\langle \gamma_\varepsilon: \varepsilon< \kappa\rangle$ we have:
$\gamma_\varepsilon\in X$ and 
\[(\forall \varepsilon<\kappa)(\exists \gamma\in X)(\exists\alpha\leq
\min(E\setminus \gamma))(\{\gamma_\zeta: \zeta< \varepsilon\}\in
{\P}_\alpha).\] 
\end{description}
\end{proposition}

\Proof Let $X=\{j_\varepsilon: \varepsilon<\theta\}$ strictly increasing
with $\varepsilon$, and let $\beta_\varepsilon=\min(E\setminus
(j_\varepsilon+1))$, so $\zeta<\varepsilon \Rightarrow \beta_\zeta\leq
\beta_\varepsilon$. Let

\begin{eqnarray}
e\stackrel{\rm def}{=}\{\varepsilon < \theta: &\ & \varepsilon \mbox{ is a
limit ordinal and } \nonumber\\
\ &\ & \mbox{if }\varepsilon_1< \varepsilon\mbox{ and }u\in
[\varepsilon_1]^{<\kappa}\mbox{ and }\{j_\xi: \xi\in u]\in
\bigcup\limits_{\zeta< \theta}{\P}_{\beta_\zeta} \nonumber\\
\ &\ & \mbox{then }\{j_\varepsilon: \varepsilon\in u\}\in
\bigcup\limits_{\zeta< \varepsilon}{\P}_{\beta_\varepsilon}\}. \nonumber
\end{eqnarray}
Now $e$ is a club of $\theta$ as ($\theta$ is regular and) $(\forall
\varepsilon<\theta)[|\varepsilon|^{<\kappa}<\theta]$. So we can apply
clause {\bf (ii)}$'$ to $X'=:\{j_\varepsilon: \varepsilon \in e\}$, and get a
subset $\{\gamma_\varepsilon: \varepsilon< \kappa\}$ as there, it is as
required in clause {\bf (ii)}. \QED

\begin{proposition}
\label{claim8.9}
\begin{enumerate}
\item Assume $\lambda>\theta=\cf(\theta)\geq
\kappa=\cf(\kappa)>\aleph_0$. Then a sufficient condition for clause (b) +
$(\delta)(*)_1$ of claim \ref{claim8.3} is
\begin{description}
\item[$(\tilde{\otimes})$\ \ (a)] $\lambda>\theta=\cf(\theta)$
\item[\qquad\ (b)] for arbitrarily large $\alpha<\lambda$ for some regular
$\tau< \theta$ and $\lambda'< \lambda$, for every ${\frak a}\subseteq \Reg
\cap |\alpha|\setminus \theta$ for some $\langle {\frak b}_\varepsilon:
\varepsilon< \varepsilon^*<\tau\rangle$ we have ${\frak a}=
\bigcup\limits_{\varepsilon< \varepsilon^*} {\frak b}_\varepsilon$ and 
$[{\frak b}_\varepsilon]^{<\kappa} \subseteq J_{\leq \lambda'}[{\frak a}]$
for every $\varepsilon< \varepsilon^*$.
\item[\qquad\ (c)] $(\forall \varepsilon<\theta)[|\varepsilon|^{<\kappa}<
\theta]$ or for every $\lambda'\in [\mu, \lambda]$, $\square_{\{\delta<
\lambda': \cf(\delta)=\theta\}}$ 
\end{description}
\item Assume $\mu> \theta\geq \kappa=\cf(\kappa)> \aleph_0$, a
sufficient condition for clause (b) of \ref{claim8.3} to hold is:

for every $\lambda\in [\mu, 2^\mu]$ of cofinality $>\mu$ for some $\theta'\leq
\theta$, $\otimes_1$ holds (with $\theta'$ instead $\theta$).
\end{enumerate}
\end{proposition}

\Proof 1) By \cite{Sh:430}, \cite{Sh:420}, 2.6 or \cite{Sh:513}.

2) Follows. \QED

\begin{remark}
{\em 
So it is still possible that (assuming CH for simplicity)
\begin{description}
\item[$\otimes$] if $\mu=\mu^{\aleph_1}$, $\ba$ is $\ccc$ Boolean
algebra, $|\ba|\leq 2^\mu$ then $\ba$ is $\mu$-linked.
\end{description}
On the required assumption see \cite{Sh:410}, Hyp. 6.1(x).

\noindent Note that the assumptions of the form $\lambda\in I[\lambda]$ if
added save us a little on $\pcf$ hyp. (we mention it only in 8.x). But if we
are interested in the $[\kappa-\cc \Rightarrow \lambda$-Knaster], it can
be waived.
}
\end{remark}

\nocite{Ju}
\nocite{M1}
\nocite{M2}
\nocite{MgSh:204}
\nocite{Sh:93}
\nocite{Sh:345a}
\nocite{Sh:355}
\nocite{Sh:430}
\nocite{Sh:481}
\nocite{Sh:572}
\nocite{Sh:g}

\bibliographystyle{literal-unsrt}
\bibliography{listb,lista,listf,listx}

\shlhetal

\end{document}